\documentclass[a4paper]{article}%
\usepackage[pagewise]{lineno}
\usepackage{graphicx}
\usepackage{geometry}
\usepackage{amsmath}
\usepackage{amssymb}
\usepackage{amsfonts}
\usepackage{amsthm,amscd}
\usepackage[colorlinks=true]{hyperref}%
\providecommand{\U}[1]{\protect\rule{.1in}{.1in}}
\hypersetup{urlcolor=blue, citecolor=red}
\def\figurename{Figure}
\makeatletter
\renewcommand{\fnum@figure}[1]{\figurename~\thefigure.}
\makeatother
\def\tablename{Table}
\makeatletter
\renewcommand{\fnum@table}[1]{\tablename~\thetable.}
\makeatother
\def \bop {\noindent\textbf{Proof. }}
\def \eop {\hbox{}\nobreak\hfill
\vrule width 2mm height 2mm depth 0mm
\par \goodbreak \smallskip}
\newtheorem{theorem}{Theorem}[section]
\newtheorem{lemma}[theorem]{Lemma}

\theoremstyle{definition}
\newtheorem{definition}[theorem]{Definition}

\theoremstyle{remark}
\newtheorem{remark}[theorem]{Remark}
\numberwithin{equation}{section}

\begin{document}

\title{On Stochastic Maximum Principle: A Backward Stochastic Partial Differential
Equations Point of View}
\author{\textbf{Ishak Alia}
\and \textbf{Mohamed Sofiane Alia}}
\maketitle

\begin{abstract}
In this paper, we consider a class of stochastic control problems for
stochastic differential equations with random coefficients. The control domain
need not to be convex but the control process is not allowed to enter in
diffusion term. Moreover, the terminal cost involves a non linear term of the
expected value of terminal state. Our purpose is to derive a new version of
the Pontryagin's stochastic maximum principle by adopting an idea inspired
from the work of Peng [S. Peng, \textit{Maximum Principle for Stochastic
Optimal Control with Nonconvex Control Domain}, Lecture Notes in Control \&
Information Sciences, 114, (1990), pp. 724-732]. More specifically, we show
that if we combine the spike perturbation of the optimal control combined with
the stochastic Feynman-Kac representation of linear backward stochastic
partial differential equations (BSPDE, for short), a new version of the
stochastic maximum principle can be derived. We also investigate sufficient
conditions of optimality. In the last part of this paper, motivated by our
version of SMP, an interesting class of forward backward stochastic partial
differential equations is naturally introduced and the solvability of such
kind of equations is briefly presented.

\end{abstract}

\textbf{Keys words}: Stochastic control, linear degenerate backward stochastic
partial differential equation, optimal control, stochastic maximum principle.


\section{Introduction\label{section1}}

In this paper we consider a class of stochastic optimal control problems where
the state of the system under consideration is governed by the controlled SDE,%
\begin{equation}
\left\{
\begin{array}
[c]{l}%
dX\left(  t\right)  =b\left(  t,X\left(  t\right)  ,u\left(  t\right)
\right)  dt+\sum\limits_{j=1}^{d}\sigma^{j}\left(  t,X\left(  t\right)
\right)  dW_{j}\left(  t\right)  ,\\
X\left(  0\right)  =x_{0}\text{,}%
\end{array}
\right.  \label{Eq1}%
\end{equation}
and the objective of the controller is to minimize the following expected cost
functional,%
\begin{equation}
\mathbf{J}\left(  u\left(  \cdot\right)  \right)  =\mathbb{E}\left[  \int
_{0}^{T}f\left(  t,X\left(  t\right)  ,u\left(  t\right)  \right)  dt+h\left(
X\left(  T\right)  \right)  +G\left(  \mathbb{E}\left[  X\left(  T\right)
\right]  \right)  \right]  \text{,} \label{Eq2}%
\end{equation}
over the set of the admissible controls. Here $W\left(  \cdot\right)  =\left(
W_{1}\left(  \cdot\right)  ,...,W_{d}\left(  \cdot\right)  \right)  ^{\top}$
is a d-dimensional standard Brownian motion defined on some filtered
probability space $(\Omega,\mathcal{F},\left(  \mathcal{F}_{t}\right)
_{t\in\left[  0,T\right]  },\mathbb{P})$ satisfying usual conditions; the
coefficients $b$, $\sigma^{j}$, $f$, $h$ are sufficiently smooth $\left(
\mathcal{F}_{t}\right)  _{t\in\left[  0,T\right]  }$-progressively measurable
functions and $G$ is a deterministic measurable function; $u\left(
\cdot\right)  :\left[  0,T\right]  \times\Omega\rightarrow U$ represent the
control process, where $\left(  U,d\right)  $ is a separable metric space;
$X\left(  \cdot\right)  $ is the corresponding state process of $u\left(
\cdot\right)  $; and $x_{0}\in\mathbb{R}^{n}$ is regarded as the initial
state. A control process $\bar{u}\left(  \cdot\right)  $ that solves this
problem is called optimal.

One of the principal methods in solving continuous-time stochastic optimal
control problems is the Pontryagin's maximum principle approach, by which
necessary conditions for optimality are derived by considering a \textit{spike
perturbation} of the optimal control (i.e. by perturbing an optimal control on
a small time interval of length $\varepsilon>0$). Combining a sort of Taylor
expansion with respect to $\varepsilon$, with Euler's necessary condition of
optimality one obtains a kind of variational inequality. By the duality, the
stochastic maximum (or minimum) principle is obtained. It states that any
optimal control along with the corresponding optimal state must solve the so
called (extended) Hamiltonian system, which consists of a forward differential
equation and a linear backward differential equation called the adjoint
equation, plus a maximum (or minimum) condition of a function called the
Hamiltonian. The efficiency of the Pontryagin's maximum principle approach
lies in the fact that maximizing (or minimizing) the Hamiltonian is much more
manageable than the original control problem which is infinite-dimensional.

The Pontryagin's maximum principle approach was first performed for
deterministic problems; see e.g. Pontryagin et al. \cite{Pontryagin}.
Extension to stochastic diffusion control problems was first carried out by
Kushner \cite{Kuch}, followed by Haussmann (\cite{Hauss1}, \cite{Hauss2}),
Bismut \cite{Bis}, Elliott and Kohlmann \cite{Elliot} and Bensoussan
\cite{Ben}. However, at that time, the results were essentially obtained under
the assumption that the diffusion term is independent of the control variable.
Specially, in Bensoussan \cite{Ben} the stochastic maximum principle (SMP, for
short) is constructed for control problems, where the state is described by a
stochastic differential equation (\ref{Eq1}) (with $n=d$) and the objective of
the controller is to minimize the expected cost functional%
\[
\mathbf{J}_{0}\left(  u\left(  \cdot\right)  \right)  =\mathbb{E}\left[
\int_{0}^{T}f\left(  t,X\left(  t\right)  ,u\left(  t\right)  \right)
dt+h\left(  X\left(  T\right)  \right)  \right]
\]
over the set of admissible controls. The author constructs his SMP by
considering a \textit{spike perturbation} of the optimal control defined by
$u^{\varepsilon}\left(  t\right)  =v\mathbf{\chi}_{\left[  \tau,\tau
+\varepsilon\right)  }\left(  t\right)  +\bar{u}\left(  t\right)
\mathbf{\chi}_{\left[  0,T\right]  \backslash\left[  \tau,\tau+\varepsilon
\right)  }\left(  t\right)  ,$ for any $\tau\in\left[  0,T\right)  ,$ $v\in U$
and for any $\varepsilon\in\left[  0,T-\tau\right)  $. Performing a Taylor
expansion, the author obtained that%
\begin{align*}
0  &  \leq\frac{d\mathbf{J}_{0}}{d\varepsilon}\left(  u^{\varepsilon}\left(
\cdot\right)  \right)  |_{\varepsilon=0}\\
&  =\mathbb{E}\left[  h_{x}\left(  \bar{X}\left(  T\right)  \right)  ^{\top
}z\left(  T\right)  +\zeta\left(  T\right)  \right]  ,
\end{align*}
where $z\left(  \cdot\right)  $ is the unique solution of the so-called
variational equation:%
\[
\left\{
\begin{array}
[c]{l}%
dz\left(  t\right)  =b_{x}\left(  t,\bar{X}\left(  t\right)  ,\bar{u}\left(
t\right)  \right)  ^{\top}z\left(  t\right)  dt+\sum\limits_{i=1}^{d}%
\sigma_{x}^{j}\left(  t,\bar{X}\left(  t\right)  \right)  ^{\top}z\left(
t\right)  dW_{j}\left(  t\right) \\
z\left(  \tau\right)  =b\left(  \tau,\bar{X}\left(  \tau\right)  ,v\right)
-b\left(  \tau,\bar{X}\left(  \tau\right)  ,\bar{u}\left(  \tau\right)
\right)
\end{array}
\right.
\]
and $\zeta\left(  \cdot\right)  $ solves the following equation%
\[
\left\{
\begin{array}
[c]{l}%
\frac{d\zeta}{dt}\left(  t\right)  =f_{x}\left(  t,\bar{X}\left(  t\right)
,\bar{u}\left(  t\right)  \right)  ^{\top}z\left(  t\right) \\
\zeta\left(  \tau\right)  =f\left(  \tau,\bar{X}\left(  \tau\right)
,v\right)  -f\left(  \tau,\bar{X}\left(  \tau\right)  ,\bar{u}\left(
\tau\right)  \right)  .
\end{array}
\right.
\]
Then by introducing the first order adjoint processes $p\left(  \cdot\right)
$ and $q\left(  \cdot\right)  =\left(  q_{1}\left(  \cdot\right)
,...,q_{d}\left(  \cdot\right)  \right)  $ as the unique solution of the BSDE:%
\begin{equation}
\left\{
\begin{array}
[c]{l}%
dp\left(  t\right)  =-\left\{  \left\langle b_{x}\left(  t,\bar{X}\left(
t\right)  ,\bar{u}\left(  t\right)  \right)  ,p\left(  t\right)  \right\rangle
+\sum\limits_{j=1}^{d}\left\langle \sigma_{x}^{j}\left(  t,\bar{X}\left(
t\right)  \right)  ,q_{j}\left(  t\right)  \right\rangle \right. \\
\text{ \ \ \ \ \ \ \ \ \ \ \ \ \ }\left.  +f_{x}\left(  t,\bar{X}\left(
t\right)  ,\bar{u}\left(  t\right)  \right)  \right\}  dt+\sum\limits_{j=1}%
^{d}q_{j}\left(  t\right)  dW_{j}\left(  t\right)  \text{, }t\in\left[
0,T\right]  \text{,}\\
p\left(  T\right)  =h_{x}\left(  \bar{X}\left(  T\right)  \right)  \text{,}%
\end{array}
\right.  \label{Eq3}%
\end{equation}
and by applying It\^{o}'s formula to $t\rightarrow\left\langle p\left(
t\right)  ,z\left(  t\right)  \right\rangle $, the stochastic maximum
principle can be obtained. It sates that the optimal pair $\left(  \bar
{u}\left(  \cdot\right)  ,\bar{X}\left(  \cdot\right)  \right)  $ should
satisfy the following variational inequality%
\begin{align*}
0  &  \leq\left\langle b\left(  \tau,\bar{X}\left(  \tau\right)  ,v\right)
-b\left(  \tau,\bar{X}\left(  \tau\right)  ,\bar{u}\left(  \tau\right)
\right)  ,p\left(  \tau\right)  \right\rangle \\
&  +%
\begin{array}
[c]{c}%
f\left(  \tau,\bar{X}\left(  \tau\right)  ,v\right)  -f\left(  \tau,\bar
{X}\left(  \tau\right)  ,\bar{u}\left(  \tau\right)  \right)  \text{, a.s.,
}\forall v\in U\text{, a.e. }\tau\in\left[  0,T\right]  .
\end{array}
\end{align*}

In the general case when the diffusion coefficient is allowed to depend on the
control, the first version of the stochastic maximum principle was derived by
Peng \cite{Peng} in 1990. To overcome the difficulty arising from the
perturbed stochastic integral term, the author considered the second-order
term in the spike variation. In the form of his stochastic maximum principle
the first order adjoint BSDE was accompanied by a second order BSDE involving
the second order derivatives of the coefficients of the system, and the
Hamiltonian was extended accordingly. In the work of Peng \cite{Peng2}, the
author introduced an alternative approach. He showed that the first and second
order adjoint processes together with the variational inequality obtained in
\cite{Peng} can be derived in a natural way via random fields that solve a
linear system of BSPDEs. With the joint effort of many researchers in the last
30 years, there is a very extensive literature on different versions of Peng's
type SMP. Let us just mention a few: see Tang and Li \cite{Tang1} for systems
with jumps, see Bahlali and Mezerdi \cite{BahMez} for stochastic singular
control, see Buckdan et al. \cite{BucDjeLi} for stochastic systems of
mean-field type, see Yong \cite{Yong2} and Hu \cite{Hu} for general coupled
forward-backward stochastic differential equations (FBSDEs). More recently,
Agram and $\varnothing$ksendal \cite{Agram} presented an new approach based on
Hida-Malliavin calculus and white noise theory, which enabled them to derive
the SMP without involving the second adjoint BSDE.

The purpose of this paper is to suggest a new version of Pontryagin's SMP to
the class of control problems (\ref{Eq1})-(\ref{Eq2}). More specifically, by
combining the spike perturbation of the optimal control with stochastic
Feynman-Kac representations of linear degenerate BSPDEs (see e.g. [Ma and Yong
\cite{Yong1997}, Section 6]), we show that the optimal pair $\left(  \bar
{u}\left(  \cdot\right)  ,\bar{X}\left(  \cdot\right)  \right)  $ should
satisfy the following variational inequality%
\begin{align*}
0  &  \leq\left\langle b\left(  \tau,\bar{X}\left(  \tau\right)  ,u\right)
-b\left(  \tau,\bar{X}\left(  \tau\right)  ,\bar{u}\left(  \tau\right)
\right)  ,\bar{\theta}_{x}\left(  \tau,\bar{X}\left(  \tau\right)  \right)
\right\rangle \\
&  +\sum_{i=1}^{n}\left\langle b\left(  \tau,\bar{X}\left(  \tau\right)
,u\right)  -b\left(  \tau,\bar{X}\left(  \tau\right)  ,\bar{u}\left(
\tau\right)  \right)  ,G_{\bar{x}_{i}}\left(  \mathbb{E}\left[  \bar{X}\left(
T\right)  \right]  \right)  \bar{g}_{x}^{i}\left(  \tau,\bar{X}\left(
\tau\right)  \right)  \right\rangle \\
&  +%
\begin{array}
[c]{c}%
f\left(  \tau,\bar{X}\left(  \tau\right)  ,u\right)  -f\left(  \tau,\bar
{X}\left(  \tau\right)  ,\bar{u}\left(  \tau\right)  \right)  ,\text{ for all
}u\in U,\text{ a.s., a.e. }\tau\in\left[  0,T\right]  ,
\end{array}
\end{align*}
and the objective value of $\bar{u}\left(  \cdot\right)  $ is given by%
\[
\mathbf{J}\left(  \bar{u}\left(  \cdot\right)  \right)  =\bar{\theta}\left(
0,x_{0}\right)  +G\left(  \bar{g}\left(  0,x_{0}\right)  \right)  \text{,}%
\]
where the pair of random fields $\left(  \bar{\theta}\left(  \cdot
,\cdot\right)  ,\bar{\psi}\left(  \cdot,\cdot\right)  \right)  $ is the
classical solution of the following linear system of BSPDE:%
\begin{equation}
\left\{
\begin{array}
[c]{l}%
d\bar{\theta}\left(  t,x\right)  =-\left\{  \left\langle \bar{\theta}%
_{x}\left(  t,x\right)  ,b\left(  t,x,\bar{u}\left(  t\right)  \right)
\right\rangle +\frac{1}{2}\text{\textbf{tr}}\left[  \sigma\sigma^{\top}\left(
t,x\right)  \bar{\theta}_{xx}\left(  t,x\right)  \right]  \right. \\
\text{ \ \ \ \ \ \ \ \ \ \ \ \ \ \ \ }+\left.  \text{\textbf{tr}}\left[
\bar{\psi}_{x}\left(  t,x\right)  \sigma\left(  t,x\right)  \right]  +f\left(
t,x,\bar{u}\left(  t\right)  \right)  \right\}  dt\\
\text{ \ \ \ \ \ \ \ \ \ \ \ \ \ \ \ }+\bar{\psi}\left(  t,x\right)  ^{\top
}dW\left(  t\right)  \text{, }\left(  t,x\right)  \in\left[  0,T\right]
\times\mathbb{R}^{n}\text{,}\\
\bar{\theta}\left(  T,x\right)  =F\left(  x\right)  \text{, for }%
x\in\mathbb{R}^{n}%
\end{array}
\right.  \label{Eq4}%
\end{equation}
and%
\[
\bar{g}\left(  t,x\right)  \equiv\left(  \bar{g}^{1}\left(  t,x\right)
,...,\bar{g}^{n}\left(  t,x\right)  \right)  ^{\top}%
\]
such that for each for $1\leq i\leq n,$ $\left(  \bar{g}^{i}\left(
\cdot,\cdot\right)  ,\bar{\eta}^{i}\left(  \cdot,\cdot\right)  \right)  $ is
the classical solution of the following linear BSPDE:%
\begin{equation}
\left\{
\begin{array}
[c]{l}%
d\bar{g}^{i}\left(  t,x\right)  =-\left\{  \left\langle \bar{g}_{x}^{i}\left(
t,x\right)  ,b\left(  t,x,\bar{u}\left(  t\right)  \right)  \right\rangle
+\frac{1}{2}\text{\textbf{tr}}\left[  \sigma\sigma^{\top}\left(  t,x\right)
\bar{g}_{xx}^{i}\left(  t,x\right)  \right]  \right. \\
\text{ \ \ \ \ \ \ \ \ \ \ \ \ }%
\begin{array}
[c]{r}%
+\left.  \text{\textbf{tr}}\left[  \bar{\eta}_{x}^{i}\left(  t,x\right)
\sigma\left(  t,x\right)  \right]  \right\}  dt+\bar{\eta}^{i}\left(
t,x\right)  ^{\top}dW\left(  t\right)  \text{,}\\
\text{for }\left(  t,x\right)  \in\left[  0,T\right]  \times\mathbb{R}%
^{n}\text{,}%
\end{array}
\\
\bar{g}^{i}\left(  T,x\right)  =\bar{x}_{i}\text{, for }x\in\mathbb{R}%
^{n}\text{,}%
\end{array}
\right.  \label{Eq4*}%
\end{equation}
with $x_{i}$ denotes the i-th coordinate of $x\in\mathbb{R}^{n}$.

There exists many other situations where BSPDEs are involved to construct an
optimal solution of a stochastic optimal control problem. For instance, BSPDEs
appear as adjoint equations in the study of stochastic maximum principle for
stochastic parabolic PDEs (see e.g. \cite{Zhou2}) and as adjoint equations of
Duncan-Mortensen-Zakai filtering equations to formulate the stochastic maximum
principle for stochastic control problems with partial information (see e.g.
\cite{Ben2}, \cite{Naga}, \cite{Zhou}). A class of fully nonlinear BSPDEs, the
so-called backward stochastic Hamilton-Jacobi-Bellman equations, appear
naturally in the dynamic programming theory of controlled non-Markovian
processes \cite{Pengg}. More recently, Alia \cite{Alia} investigated
equilibrium solutions for a general class of time-inconsistent control
problems by using BSPDEs. However, to our best knowledge, the above-described
version of SMP for optimal controls seems to be new. Moreover, it permits us
to derive the optimal solution $\bar{u}\left(  \cdot\right)  $ as well as its
objective value\ $\mathbf{J}\left(  \bar{u}\left(  \cdot\right)  \right)  $ in
terms of the solutions of the BSPDEs (\ref{Eq4})-(\ref{Eq4*}); this is
different from the traditional SMP approach which does not provide the
objective value $\mathbf{J}\left(  \bar{u}\left(  \cdot\right)  \right)  $ in
a direct way.

In the latter part of the paper, we investigate a class of FBSPDEs which
naturally arises from the above-described version of SMP. These systems differ
from the classical FBSDEs (see e.g \cite{MY}), since they are consisting of a
forward SDE and a nonlinear BSPDE. Inspired by the idea of the classical
four-step scheme introduced by Ma et al. \cite{MPY}, we show that a solution
of the FBSPDEs can be constructed by solving a parabolic PDE. Under proper
conditions on the involved coefficients, a well-posedness result of the PDE is established.

The plan of the paper is as follows, in the second section, we give necessary
notations and some preliminaries on linear BSPDEs. In Section \ref{section3},
we formulate our stochastic optimal control problem. Section \ref{section4} is
devoted to the new version of the SMP. Section \ref{section5} is devoted to
the proof of the first main result in this paper. In Section \ref{section6},
we investigate sufficient conditions of optimality and we give a simple
example to illustrate our result. Finally, in Section \ref{section7}, a class
of FBSPDEs is briefly presented.

\section{Preliminaries\label{section2}}

\subsection{Notations}

Throughout this paper $(\Omega,\mathcal{F},\left(  \mathcal{F}_{t}\right)
_{t\in\left[  0,T\right]  },\mathbb{P})$ is a filtered probability space such
that $\mathcal{F}_{0}$ contains all $\mathbb{P}$-null sets, $\mathcal{F}%
_{T}=\mathcal{F}$ for an arbitrarily fixed finite time horizon $T>0,$ and
$\left(  \mathcal{F}_{t}\right)  _{t\in\left[  0,T\right]  }$ satisfies the
usual conditions. $\mathcal{F}_{t}$ stands for the information available up to
time $t$ and any decision made at time $t$ is based on this information. We
also assume that all processes and random variables are well defined and
adapted in this filtered probability space. Let $W\left(  \cdot\right)
=\left(  W_{1}\left(  \cdot\right)  ,...,W_{d}\left(  \cdot\right)  \right)  $
be a d-dimensional standard Brownian motion defined on $(\Omega,\mathcal{F}%
,\left(  \mathcal{F}_{t}\right)  _{t\in\left[  0,T\right]  },\mathbb{P}).$ For
simplicity, it is assumed that the filtration $\left(  \mathcal{F}_{t}\right)
_{t\in\left[  0,T\right]  }$, coincides with the one generated by the Brownian
motion; that is $\mathcal{F}_{t}=\sigma\left(  W\left(  r\right)  ;\text{
}0\leq r\leq t\right)  .$

We use $C^{\top}$ to denote the transpose of any vector or matrix $C$, $K>0$
is a generic constant which can be different from line to line and for a
function $f$, we denote by $f_{x}$ (resp. $f_{xx}$) the gradient or Jacobian
(resp. the Hessian) of $f$ with respect to the variable $x.$ We denote by
$\chi_{A}$ the indicator function of the set $A$. In addition, we use the
following notations for several sets and spaces of processes on the filtered
probability space, which will be used later:

\begin{enumerate}
\item[$\bullet$] $\bar{B}_{R}=\left\{  x\in\mathbb{R}^{n},\text{ such that
}\left\vert x\right\vert \leq R\right\}  $, for any $R>0$.

\item[$\bullet$] $\left(  U,d\right)  :$ a separable metric space.

\item[$\bullet$] $\mathcal{U}\left[  0,T\right]  :$ the set of $U$-valued
$\left(  \mathcal{F}_{t}\right)  _{t\in\left[  0,T\right]  }$-progressively
measurable processes $u\left(  \cdot\right)  $.
\end{enumerate}

For any $m\geq1$ and $p\geq2$ we denote by

\begin{enumerate}
\item[$\bullet$] $\mathbb{L}^{p}\left(  \Omega,\mathcal{F}_{t},\mathbb{P};%
\mathbb{R}
^{l}\right)  $: the set of $\mathbb{R}^{l}$-valued, $\mathcal{F}_{t}%
$-measurable random variables $\zeta,$ with%
\[
\left\Vert \zeta\right\Vert _{\mathbb{L}^{p}\left(  \Omega,\mathcal{F}%
_{t},\mathbb{P};%
\mathbb{R}
^{n}\right)  }=\left(  \mathbb{E}\left[  \left\vert \zeta\right\vert
^{p}\right]  \right)  ^{\frac{1}{p}}<\infty.
\]

\item[$\bullet$] $\mathcal{C}_{\mathcal{F}}^{p}\left(  0,T;\mathbb{R}%
^{l}\right)  $: the space of $\mathbb{R}^{l}$-valued, $\left(  \mathcal{F}%
_{t}\right)  _{t\in\left[  0,T\right]  }$-adapted continuous processes
$X\left(  \cdot\right)  $, with%
\[
\left\Vert X\left(  \cdot\right)  \right\Vert _{\mathcal{C}_{\mathcal{F}}%
^{p}\left(  0,T;\mathbb{R}^{l}\right)  }:=\left(  \mathbb{E}\left[
\sup\limits_{t\in\left[  0,T\right]  }\left\vert X\left(  t\right)
\right\vert ^{p}\right]  \right)  ^{\frac{1}{p}}<\infty\text{.}%
\]

\item[$\bullet$] $\mathcal{L}_{\mathcal{F}}^{p}\left(  0,T;\mathbb{R}%
^{l}\right)  $: the space of $\mathbb{R}^{l}$-valued, $\left(  \mathcal{F}%
_{t}\right)  _{t\in\left[  0,T\right]  }$-adapted processes $Z\left(
\cdot\right)  $, with%
\[
\left\Vert Z\left(  \cdot\right)  \right\Vert _{\mathcal{L}_{\mathcal{F}}%
^{p}\left(  0,T;\mathbb{R}^{l}\right)  }:=\left(  \mathbb{E}\left[  {\int
_{0}^{T}}\left\vert Z\left(  t\right)  \right\vert ^{p}dt\right]  \right)
^{\frac{1}{p}}<\infty\text{.}%
\]

\item[$\bullet$] $C^{m}\left(  \mathbb{R}^{n};\mathbb{R}^{l}\right)  $: the
set of functions from $\mathbb{R}^{n}$ to $\mathbb{R}^{l}$ that are
continuously differentiable up to order $m\geq1$.

\item[$\bullet$] $C_{b}^{m}\left(  \mathbb{R}^{n};\mathbb{R}^{l}\right)  $:
the set of those functions in $\mathcal{C}^{m}\left(  \mathbb{R}%
^{n};\mathbb{R}^{l}\right)  $ whose partial derivatives up to order $m$ are
uniformly bounded.
\end{enumerate}

Finally, for any infinite-dimensional Banach space $E$ with a norm $\left\Vert
\cdot\right\Vert _{E}$, we denote.

\begin{enumerate}
\item[$\bullet$] $\mathcal{C}_{\mathcal{F}}^{p}\left(  0,T;E\right)  $: the
set of all $\left(  \mathcal{F}_{t}\right)  _{t\in\left[  0,T\right]  }%
$-adapted $E$-valued continuous processes $\theta\left(  \cdot\right)  $ such
that:%
\[
\left\Vert \theta\left(  \cdot\right)  \right\Vert _{\mathcal{C}_{\mathcal{F}%
}^{p}\left(  0,T;E\right)  }:=\sup_{t\in\left[  0,T\right]  }\left(
\mathbb{E}\left[  \left\Vert \theta\left(  t\right)  \right\Vert _{E}%
^{p}\right]  \right)  ^{\frac{1}{p}}<\infty\text{.}%
\]

\item[$\bullet$] $\mathcal{L}_{\mathcal{F}}^{p}\left(  0,T;E\right)  $: the
set of all $\left(  \mathcal{F}_{t}\right)  _{t\in\left[  0,T\right]  }%
$-adapted $E$-valued processes $\psi\left(  \cdot\right)  $ such that:%
\[
\left\Vert \psi\left(  \cdot\right)  \right\Vert _{\mathcal{L}_{\mathcal{F}%
}^{p}\left(  0,T;E\right)  }:=\left(  \mathbb{E}\left[  {\int_{0}^{T}%
}\left\Vert \psi\left(  t\right)  \right\Vert _{E}^{p}dt\right]  \right)
^{\frac{1}{p}}<\infty\text{.}%
\]

\end{enumerate}

\subsection{Linear Backward Stochastic Partial Differential Equations}

In this paragraph, we review some well-known results on linear BSPDEs. Let
$b:\left[  0,T\right]  \times\mathbb{R}^{n}\times\Omega\rightarrow
\mathbb{R}^{n}$, $\sigma:\left[  0,T\right]  \times\mathbb{R}^{n}\times
\Omega\rightarrow\mathbb{R}^{n\times d}$, $l:\left[  0,T\right]
\times\mathbb{R}^{n}\times\Omega\rightarrow\mathbb{R}$ be three $\left(
\mathcal{F}_{t}\right)  _{t\in\left[  0,T\right]  }$-progressively measurable
functions and $F:\mathbb{R}^{n}\times\Omega\rightarrow\mathbb{R}$ be an
$\mathcal{F}_{T}$-measurable function. To simplify our notation, for
$\varphi=b$, $\sigma$, $l$ we write $\varphi\left(  t,x\right)  $ for
$\varphi\left(  t,x,\omega\right)  $ and $h\left(  x\right)  $ for $h\left(
x,\omega\right)  $.

Consider the following linear backward stochastic partial differential
equation in the unknown random fields $p\left(  t,x\right)  \in%
\mathbb{R}
$, $q\left(  t,x\right)  =\left(  q^{1}\left(  t,x\right)  ,...,q^{d}\left(
t,x\right)  \right)  ^{\top}\in%
\mathbb{R}
^{d}$:%
\begin{equation}
\left\{
\begin{array}
[c]{l}%
dp\left(  t,x\right)  =-\left\{  \left\langle p_{x}\left(  t,x\right)
,b\left(  t,x\right)  \right\rangle +\frac{1}{2}\text{\textbf{tr}}\left[
\sigma\left(  t,x\right)  \sigma\left(  t,x\right)  ^{\top}p_{xx}\left(
t,x\right)  \right]  \right. \\
\text{ \ \ \ \ \ \ \ \ \ \ \ \ \ \ \ }+\left.  \text{\textbf{tr}}\left[
\sigma\left(  t,x\right)  q_{x}\left(  t,x\right)  \right]  +l\left(
t,x\right)  \right\}  dt\\
\text{ \ \ \ \ \ \ \ \ \ \ \ \ \ \ \ }+\left\langle q\left(  t,x\right)
,dW\left(  t\right)  \right\rangle \text{, }\left(  t,x\right)  \in\left[
0,T\right]  \times\mathbb{R}^{n}\text{,}\\
p\left(  T,x\right)  =h\left(  x\right)  \text{, for }x\in\mathbb{R}%
^{n}\text{.}%
\end{array}
\right.  \label{Eq5}%
\end{equation}

A pair of random fields $\left(  p\left(  \cdot,\cdot\right)  ,q\left(
\cdot,\cdot\right)  \right)  $ is called a \textit{classical solution} of
(\ref{Eq5}) if (see e.g. \cite{Yong1999})%
\[
\left\{
\begin{array}
[c]{l}%
p\left(  \cdot,\cdot\right)  \in\mathcal{C}_{\mathcal{F}}^{2}\left(
0,T;C^{2}\left(  \bar{B}_{R};\mathbb{R}\right)  \right)  ,\\
q\left(  \cdot,\cdot\right)  \in\mathcal{L}_{\mathcal{F}}^{2}\left(
0,T;C^{1}\left(  \bar{B}_{R};\mathbb{R}^{d}\right)  \right)  ,
\end{array}
\right.  \forall R>0,
\]
such that the following holds for all $\left(  t,x\right)  \in\left[
0,T\right]  \times\mathbb{R}^{n}$, almost surely:%
\begin{align*}
p\left(  t,x\right)   &  =h\left(  x\right)  +\int_{t}^{T}\left\{
\left\langle p_{x}\left(  \tau,x\right)  ,b\left(  \tau,x\right)
\right\rangle \right. \\
&  +\frac{1}{2}\text{\textbf{tr}}\left[  \sigma\left(  \tau,x\right)
\sigma\left(  \tau,x\right)  ^{\top}p_{xx}\left(  \tau,x\right)  \right]
+\text{\textbf{tr}}\left[  \sigma\left(  \tau,x\right)  q_{x}\left(
\tau,x\right)  \right] \\
&  \left.  +l\left(  \tau,x\right)  \right\}  d\tau-\int_{t}^{T}\left\langle
q\left(  \tau,x\right)  ,dW\left(  \tau\right)  \right\rangle .
\end{align*}

The existence and uniqueness of a solution to linear or semi-linear degenerate
BSPDEs driven by a Brownian motion was first proved by Ma and Yong
(\cite{Yong1997}, \cite{Yong1999}), followed by Hu et al. \cite{Huetal2002},
Tang \cite{Tang}, Du et al. \cite{Kai}, Du and Zhang \cite{Kai1} and Ma et al.
\cite{Ma2012} in different frameworks. More recently, Chen and Tang
\cite{Chen} dealt with semi-linear backward stochastic integral partial
differential equations with jumps. In this work, we mostly focus on the work
of Tang \cite{Tang}, in which the author developed a probabilistic approach to
construct the adapted solution of semi-linear BSPDEs in terms of those of SDEs
and BSDEs.

Given an integer $m\geq1$, we consider the following assumption.

\begin{enumerate}
\item[\textbf{(A}$_{m}$\textbf{)}] The random fields $b$, $\sigma$, $f$ and
$h$ satisfy the following conditions:%
\[
\left\{
\begin{array}
[c]{l}%
b\in\mathcal{L}_{\mathcal{F}}^{\infty}\left(  0,T;C_{b}^{m}\left(
\mathbb{R}^{n};\mathbb{R}^{n}\right)  \right)  \text{,}\\
\sigma\in\mathcal{L}_{\mathcal{F}}^{\infty}\left(  0,T;C_{b}^{m}\left(
\mathbb{R}^{n};\mathbb{R}^{n\times d}\right)  \right)  \text{,}\\
f\in\mathcal{L}_{\mathcal{F}}^{\infty}\left(  0,T;C_{b}^{m}\left(
\mathbb{R}^{n};\mathbb{R}\right)  \right)  \text{,}\\
h\in\mathcal{L}_{\mathcal{F}_{T}}^{\infty}\left(  \Omega;C_{b}^{m}\left(
\mathbb{R}^{n};\mathbb{R}\right)  \right)  \text{.}%
\end{array}
\right.
\]

\end{enumerate}

The following theorem concerns the existence and uniqueness of a classical
solution to the BSPDEs (\ref{Eq5}).

\begin{theorem}
[\cite{Tang}]\label{result1}Let Assumption \textbf{(A}$_{m}$\textbf{)} be
satisfied with $m>2+\frac{n}{2}$. Then the BSPDE (\ref{Eq5}) admits a unique
adapted classical solution $\left(  p\left(  \cdot,\cdot\right)  ,q\left(
\cdot,\cdot\right)  \right)  .$ Moreover, the partial derivatives in $x$ of
$p\left(  \cdot,\cdot\right)  $ up to order $m$ are uniformly bounded.
\end{theorem}

The following theorem provides the stochastic Feynman--Kac representation of
the linear degenerate BSPDE (\ref{Eq5}).

\begin{theorem}
[Feynman--Kac representation]\label{result02}Suppose that Assumption
\textbf{(A}$_{m}$\textbf{)} is satisfied with $m>2+\frac{n}{2}$. Let $\left(
p\left(  \cdot,\cdot\right)  ,q\left(  \cdot,\cdot\right)  \right)  $ be the
classical solution of the BSPDEs (\ref{Eq5}), then we have the following
equality: For each $\left(  t,\zeta\right)  \in\left[  0,T\right]
\times\mathbb{L}^{2}\left(  \Omega,\mathcal{F}_{t},\mathbb{P};%
\mathbb{R}
^{n}\right)  $,%
\begin{equation}
p\left(  t,\zeta\right)  =\mathbb{E}_{t}\left[  \int_{t}^{T}l\left(
\tau,X^{t,\zeta}\left(  \tau\right)  \right)  d\tau+h\left(  X^{t,\zeta
}\left(  T\right)  \right)  \right]  \text{,} \label{Eq102}%
\end{equation}
where $\mathbb{E}_{t}\left[  \cdot\right]  =\mathbb{E}_{t}\left[
\mathbb{\cdot}\left\vert \mathcal{F}_{t}\right.  \right]  $ and $X^{t,\zeta
}\left(  \cdot\right)  \in\mathcal{C}_{\mathcal{F}}^{2}\left(  t,T;\mathbb{R}%
^{n}\right)  $ is the unique strong solution of the following SDE,%
\begin{equation}
\left\{
\begin{array}
[c]{l}%
dX^{t,\zeta}\left(  \tau\right)  =b\left(  \tau,X^{t,\zeta}\left(
\tau\right)  \right)  d\tau+\sigma\left(  \tau,X^{t,\zeta}\left(  \tau\right)
\right)  dW\left(  \tau\right)  \text{, }\tau\in\left[  t,T\right]  \text{,}\\
X^{t,\zeta}\left(  t\right)  =\zeta\text{.}%
\end{array}
\right.  \label{Eq100}%
\end{equation}

\end{theorem}

Before presenting a proof of the above theorem we first recall a version of an
It\^{o}-Wentzell formula to the composition of random fields and continuous
semimartingale processes; see e.g. Kunita (\cite{Kunita1}, \cite{Kunita2},
\cite{Kunita3}), Krylov \cite{Krylov} and Chen and Tang \cite{Chen}.

\begin{lemma}
[It\^{o}-Wentzell Formula]Let $X\left(  \cdot\right)  \in\mathcal{C}%
_{\mathcal{F}}^{2}\left(  0,T;\mathbb{R}^{n}\right)  $ be a process of the
form%
\[
dX\left(  t\right)  =b\left(  t\right)  dt+\sigma\left(  t\right)  dW\left(
t\right)  ,
\]
where $b\left(  \cdot\right)  \in\mathcal{L}_{\mathcal{F}}^{\infty}\left(
0,T;\mathbb{R}^{n}\right)  $ and $\sigma\left(  \cdot\right)  \in
\mathcal{L}_{\mathcal{F}}^{\infty}\left(  0,T;\mathbb{R}^{n\times d}\right)
$. Suppose that $V\left(  \cdot,\cdot\right)  \in\mathcal{C}_{\mathcal{F}}%
^{2}\left(  0,T;C^{2}\left(  \bar{B}_{R};\mathbb{R}\right)  \right)  $,
$\forall R>0$, is a semimartingale with spatial parameter $x\in\mathbb{R}$:%
\[
dV\left(  t,x\right)  =\Gamma\left(  t,x\right)  dt+\zeta\left(  t,x\right)
^{\top}dW\left(  t\right)  \text{, }\left(  t,x\right)  \in\left[  0,T\right]
\times\mathbb{R}\text{,}%
\]
with $\Gamma\left(  \cdot,\cdot\right)  \in\mathcal{L}_{\mathcal{F}}%
^{1}\left(  0,T;C\left(  \bar{B}_{R};\mathbb{R}\right)  \right)  $ and
$\zeta\left(  \cdot,\cdot\right)  \in\mathcal{L}_{\mathcal{F}}^{2}\left(
0,T;C^{1}\left(  \bar{B}_{R};\mathbb{R}^{d}\right)  \right)  $, $\forall R>0$.
Then the following holds for all $0\leq t\leq s\leq T$, almost surely,%
\begin{align*}
&  V\left(  s,X\left(  s\right)  \right)  -V\left(  t,X\left(  t\right)
\right) \\
&  =\int_{t}^{s}\left\{  \Gamma\left(  \tau,X\left(  \tau\right)  \right)
+\left\langle V_{x}\left(  \tau,X\left(  \tau\right)  \right)  ,b\left(
\tau\right)  \right\rangle \right. \\
&  \left.  +\text{\textbf{tr}}\left[  \zeta_{x}\left(  \tau,X\left(
\tau\right)  \right)  \sigma\left(  \tau\right)  \right]  +\frac{1}%
{2}\text{\textbf{tr}}\left[  \sigma\left(  \tau\right)  \sigma\left(
\tau\right)  ^{\top}V_{xx}\left(  \tau,X\left(  \tau\right)  \right)  \right]
\right\}  d\tau\\
&  +\int_{t}^{s}\left\{  \zeta\left(  \tau,X\left(  \tau\right)  \right)
^{\top}+V_{x}\left(  \tau,X\left(  \tau\right)  \right)  ^{\top}\sigma\left(
\tau\right)  \right\}  dW\left(  \tau\right)  \text{.}%
\end{align*}

\end{lemma}

\bop A proof of this lemma can be easily obtained by adapting the proof\ of
Lemma 3.1 in \cite{Chen}. We omit it.\eop

\textbf{Proof of Theorem \ref{result02}}. Let $\left(  p\left(  \cdot
,\cdot\right)  ,q\left(  \cdot,\cdot\right)  \right)  $ be the classical
solution of the BSPDEs (\ref{Eq5}) and $X^{t,x}\left(  \cdot\right)  $ is the
unique strong solution of the SDE (\ref{Eq100}). Define%
\[
\mathbf{J:=}\mathbb{E}_{t}\left[  \int_{t}^{T}l\left(  \tau,X^{t,\zeta}\left(
\tau\right)  \right)  d\tau+h\left(  X^{t,\zeta}\left(  T\right)  \right)
\right]  .
\]
By the terminal condition in (\ref{Eq5}) we have%
\begin{align}
\mathbf{J}  &  \mathbf{=}\mathbb{E}_{t}\left[  \int_{t}^{T}l\left(
\tau,X^{t,\zeta}\left(  \tau\right)  \right)  d\tau+p\left(  T,X^{t,\zeta
}\left(  T\right)  \right)  \right] \nonumber\\
&  =\mathbb{E}_{t}\left[  \int_{t}^{T}l\left(  \tau,X^{t,\zeta}\left(
\tau\right)  \right)  d\tau+p\left(  t,\zeta\right)  \right] \nonumber\\
&  +\mathbb{E}_{t}\left[  p\left(  T,X^{t,\zeta}\left(  T\right)  \right)
-p\left(  t,\zeta\right)  \right]  . \label{Eq101}%
\end{align}
Moreover, by applying It\^{o}-Wentzell formula to $p\left(  \tau,X^{t,\zeta
}\left(  \tau\right)  \right)  $ on time interval $\left[  t,T\right]  $, we
obtain that%
\[
\mathbb{E}_{t}\left[  p\left(  T,X^{t,\zeta}\left(  T\right)  \right)
-p\left(  t,\zeta\right)  \right]  =-\mathbb{E}_{t}\left[  \int_{t}%
^{T}l\left(  \tau,X^{t,\zeta}\left(  \tau\right)  \right)  d\tau\right]  .
\]
where $\mathbb{E}_{t}\left[  \cdot\right]  =\mathbb{E}\left[  \cdot
|\mathcal{F}_{t}\right]  $ is the conditional expectation with respect to
$\mathcal{F}_{t}$. Invoking this into (\ref{Eq101}), we obtain (\ref{Eq102}%
).\eop

\section{Formulation of the problem\label{section3}}

We consider a continuous-time, n-dimensional, controlled system%
\begin{equation}
\left\{
\begin{array}
[c]{l}%
dX\left(  t\right)  =b\left(  t,X\left(  t\right)  ,u\left(  t\right)
\right)  dt+\sigma\left(  t,X\left(  t\right)  \right)  dW\left(  t\right)
\text{, }t\in\left[  0,T\right]  \text{,}\\
X\left(  0\right)  =x_{0}\text{.}%
\end{array}
\right.  \label{Eq11}%
\end{equation}
Here $b:\left[  0,T\right]  \times\mathbb{R}^{n}\times U\times\Omega
\rightarrow\mathbb{R}^{n}$ and $\sigma:\left[  0,T\right]  \times
\mathbb{R}^{n}\times\Omega\rightarrow\mathbb{R}^{n\times d}$ are $\left(
\mathcal{F}_{t}\right)  _{t\in\left[  0,T\right]  }$-progressively measurable
functions; $u\left(  \cdot\right)  \in\mathcal{U}\left[  0,T\right]  $
represents an admissible control process; $X\left(  \cdot\right)
=X^{x_{0},u\left(  \cdot\right)  }\left(  \cdot\right)  $ is the controlled
state process; and $x_{0}\in\mathbb{R}$ is regarded as the initial state.

In order to evaluate the cost-performance of a control process $u\left(
\cdot\right)  $, we introduce the functional%
\begin{equation}
\mathbf{J}\left(  u\left(  \cdot\right)  \right)  :=\mathbb{E}\left[  \int
_{0}^{T}f\left(  t,X\left(  t\right)  ,u\left(  t\right)  \right)  dt+h\left(
X\left(  T\right)  \right)  +G\left(  \mathbb{E}\left[  X\left(  T\right)
\right]  \right)  \right]  \label{Eq12}%
\end{equation}
where $f:\left[  0,T\right]  \times\mathbb{R}^{n}\times U\times\Omega
\rightarrow%
\mathbb{R}
$ is an $\left(  \mathcal{F}_{t}\right)  _{t\in\left[  0,T\right]  }%
$-progressively measurable function, $h:%
\mathbb{R}
^{n}\times\Omega\rightarrow%
\mathbb{R}
$, is an $\mathcal{F}_{T}-$measurable function and $G:%
\mathbb{R}
^{n}\rightarrow%
\mathbb{R}
$ is a deterministic function.

Before going further, we introduce some notations. For any fixed control
$u\left(  \cdot\right)  \in\mathcal{U}\left[  0,T\right]  ,$ we define the
random fields $b^{u\left(  \cdot\right)  }:\left[  0,T\right]  \times
\mathbb{R}^{n}\times\Omega\rightarrow\mathbb{R}^{n}$ and $f^{u\left(
\cdot\right)  }:\left[  0,T\right]  \times\mathbb{R}^{n}\times\Omega
\rightarrow\mathbb{R}$ by%
\begin{align*}
b^{u\left(  \cdot\right)  }\left(  t,x\right)   &  :=b\left(  t,x,u\left(
t\right)  \right)  ,\\
f^{u\left(  \cdot\right)  }\left(  t,x\right)   &  :=f\left(  t,x,u\left(
t\right)  \right)  .
\end{align*}

Given an integer $m\geq1$, we shall use the following assumption.

\begin{enumerate}
\item[\textbf{(H}$_{m}$\textbf{)}] For each fixed $u\left(  \cdot\right)
\in\mathcal{U}\left[  0,T\right]  ,$ the random fields $b^{u\left(
\cdot\right)  }$, $\sigma$, $f^{u\left(  \cdot\right)  }$, $h$ and $G$ satisfy
the following conditions:%
\[
\left\{
\begin{array}
[c]{l}%
b^{u\left(  \cdot\right)  }\in\mathcal{L}_{\mathcal{F}}^{\infty}\left(
0,T;C_{b}^{m}\left(  \mathbb{R}^{n};\mathbb{R}^{n}\right)  \right)  \text{,}\\
\sigma\in\mathcal{L}_{\mathcal{F}}^{\infty}\left(  0,T;C_{b}^{m}\left(
\mathbb{R}^{n};\mathbb{R}^{n\times d}\right)  \right)  \text{,}\\
f^{u\left(  \cdot\right)  }\in\mathcal{L}_{\mathcal{F}}^{\infty}\left(
0,T;C_{b}^{m}\left(  \mathbb{R}^{n};\mathbb{R}\right)  \right)  \text{,}\\
F\in\mathcal{L}_{\mathcal{F}_{T}}^{\infty}\left(  \Omega;C_{b}^{m}\left(
\mathbb{R}^{n};\mathbb{R}\right)  \right)  \text{,}\\
G\in\mathcal{L}^{\infty}\left(  \Omega;C_{b}^{2}\left(  \mathbb{R}%
^{n};\mathbb{R}\right)  \right)
\end{array}
\right.
\]

\end{enumerate}

Note that \textbf{(H}$_{m}$\textbf{)} implies that the derivatives in $x$ of
$b^{u\left(  \cdot\right)  }$ and $\sigma$ up to order $m\geq1$ are uniformly
bounded by some positive constant. Under \textbf{(H}$_{m}$\textbf{)}, the
state equation (\ref{Eq11}) admits a unique strong solution $X\left(
\cdot\right)  \in\mathcal{C}_{\mathcal{F}}^{2}\left(  0,T;\mathbb{R}%
^{n}\right)  $ (see e.g. \cite{YongZhou}). Moreover, there exists a constant
$K>0$ such that%
\[
\mathbb{E}\left[  \sup_{0\leq t\leq T}\left\vert X\left(  t\right)
\right\vert ^{2}\right]  \leq K\left(  1+\left\vert x_{0}\right\vert
^{2}\right)  .
\]

Our stochastic optimal control problem can be stated as follows.

\subparagraph{Problem (S).}

\textit{Minimize (\ref{Eq12}) over} $\mathcal{U}\left[  0,T\right]  .$

Any $\bar{u}\left(  \cdot\right)  $ satisfying%
\[
\mathbf{J}\left(  \bar{u}\left(  \cdot\right)  \right)  =\inf_{u\left(
\cdot\right)  \in\mathcal{U}\left[  0,T\right]  }\mathbf{J}\left(  u\left(
\cdot\right)  \right)
\]
is called an optimal control. The corresponding $\bar{X}\left(  \cdot\right)
=X^{x_{0},\bar{u}\left(  \cdot\right)  }\left(  \cdot\right)  $ and $\left(
\bar{u}\left(  \cdot\right)  ,\bar{X}\left(  \cdot\right)  \right)  $ are
called an optimal state process and optimal pair, respectively.

\section{Necessary Condition For Optimality\label{section4}}

As well mentioned in Introduction Section, the first order Taylor expansion in
the spike variation together with the pair of adjoint processes $\left(
p\left(  \cdot\right)  ,q\left(  \cdot\right)  \right)  $ play a central role
in the traditional stochastic maximum principle approach \cite{Ben}. Moreover,
the adjoint equation that $\left(  p\left(  \cdot\right)  ,q\left(
\cdot\right)  \right)  $ satisfies is a standard BSDE having
finite-dimensional state process.

In this paper, inspired by Peng \cite{Peng2}, we follow an alternative
approach which permits us to derive a stochastic maximum principle of Problem
(S) without performing a Taylor expansion and without involving the adjoint
processes $\left(  p\left(  \cdot\right)  ,q\left(  \cdot\right)  \right)  $.
More specifically, we show that if we use the spike perturbation of the
optimal control combined with stochastic Feynman-Kac representations of linear
degenerate BSPDEs, we can obtain a new version of the stochastic maximum principle.

In the next, we introduce the backward stochastic partial differential
equations involved in the new version of the SMP.

For a given optimal control $\bar{u}\left(  \cdot\right)  \in\mathcal{U}%
\left[  0,T\right]  $, we define in $\left(  t,x\right)  \in\left[
0,T\right]  \times\mathbb{R}^{n},$ the following \textit{linear degenerate
backward stochastic partial differential equation} in the unknown random
fields $\bar{\theta}\left(  t,x\right)  \in\mathbb{R}$, $\bar{\psi}\left(
t,x\right)  \in\mathbb{R}^{d}$:%
\begin{equation}
\left\{
\begin{array}
[c]{l}%
d\bar{\theta}\left(  t,x\right)  =-\left\{  \left\langle \bar{\theta}%
_{x}\left(  t,x\right)  ,b^{\bar{u}\left(  \cdot\right)  }\left(  t,x\right)
\right\rangle +\frac{1}{2}\text{\textbf{tr}}\left[  \left(  \sigma\sigma
^{\top}\right)  ^{u\left(  \cdot\right)  }\left(  t,x\right)  \bar{\theta
}_{xx}\left(  t,x\right)  \right]  \right. \\
\text{ \ \ \ \ \ \ \ \ \ \ \ \ \ \ \ }+\left.  \text{\textbf{tr}}\left[
\bar{\psi}_{x}\left(  t,x\right)  \sigma\left(  t,x\right)  \right]
+f^{\bar{u}\left(  \cdot\right)  }\left(  t,x\right)  \right\}  dt\\
\text{ \ \ \ \ \ \ \ \ \ \ \ \ \ \ \ }+\bar{\psi}\left(  t,x\right)  ^{\top
}dW\left(  t\right)  \text{, }\left(  t,x\right)  \in\left[  0,T\right]
\times\mathbb{R}^{n}\text{,}\\
\bar{\theta}\left(  T,x\right)  =h\left(  x\right)  \text{, for }%
x\in\mathbb{R}^{n}\text{.}%
\end{array}
\right.  \label{Eq13}%
\end{equation}

Similarly, we define in $\left(  t,x\right)  \in\left[  0,T\right]
\times\mathbb{R}^{n}$ the following system of linear degenerate BSPDEs in the
unknown random fields $\bar{g}\left(  s,x\right)  =\left(  \bar{g}^{1}\left(
s,x\right)  ,...,\bar{g}^{n}\left(  s,x\right)  \right)  ^{\top}\in
\mathbb{R}^{n},$ $\bar{\eta}\left(  s,x\right)  =\left(  \bar{\eta}^{1}\left(
s,x\right)  ,...,\bar{\eta}^{n}\left(  s,x\right)  \right)  \in\mathbb{R}%
^{d\times n}$: For $1\leq i\leq n,$%
\begin{equation}
\left\{
\begin{array}
[c]{l}%
d\bar{g}^{i}\left(  t,x\right)  =-\left\{  \left\langle \bar{g}_{x}^{i}\left(
t,x\right)  ,b^{\bar{u}\left(  \cdot\right)  }\left(  t,x\right)
\right\rangle +\frac{1}{2}\text{\textbf{tr}}\left[  \left(  \sigma\sigma
^{\top}\right)  \left(  t,x\right)  \bar{g}_{xx}\left(  t,x\right)  \right]
\right. \\
\text{ \ \ \ \ \ \ \ \ \ \ \ \ }%
\begin{array}
[c]{r}%
+\left.  \text{\textbf{tr}}\left[  \sigma\left(  t,x\right)  \bar{\eta}%
_{x}^{i}\left(  t,x\right)  \right]  \right\}  dt+\bar{\eta}^{i}\left(
t,x\right)  ^{\top}dW\left(  t\right)  \text{,}\\
\text{for }\left(  t,x\right)  \in\left[  0,T\right]  \times\mathbb{R}%
^{n}\text{,}%
\end{array}
\\
\bar{g}^{i}\left(  T,x\right)  =x_{i}\text{, for }x\in\mathbb{R}^{n}.
\end{array}
\right.  \label{Eq13**}%
\end{equation}
where $\bar{\eta}^{i}\left(  s,x\right)  $ is the i-th column of $\bar{\eta
}\left(  s,x\right)  $ and $x_{i}$ denotes the i-th coordinate of
$x\in\mathbb{R}^{n}$.

According to Theorem \ref{result1}, under Assumption \textbf{(H}$_{m}%
$\textbf{),} with $m>2+\frac{n}{2}$, the BSPDE (\ref{Eq13}) admits a unique
adapted classical solution $\left(  \bar{\theta}\left(  \cdot,\cdot\right)
,\bar{\psi}\left(  \cdot,\cdot\right)  \right)  $. Moreover, all the partial
derivatives in $x$ of $\bar{\theta}\left(  \cdot,\cdot\right)  $ up to order
$m$ are uniformly bounded.

Since the terminal condition of (\ref{Eq13**}) is unbounded, we can not
directly apply Theorem \ref{result1} to obtain the well-posedness of
(\ref{Eq13**}). Nonetheless, we state the following result which by means of
easy manipulations shows that the well-posedness of (\ref{Eq13**}) holds indeed.

\begin{theorem}
Let Assumption \textbf{(H}$_{m}$\textbf{)} be satisfied with $m>2+\frac{n}{2}%
$. Then for each $1\leq i\leq n$, the BSPDE (\ref{Eq13**}) admits a unique
adapted classical solution $\left(  \bar{g}^{i}\left(  \cdot,\cdot\right)
,\bar{\eta}^{i}\left(  \cdot,\cdot\right)  \right)  .$ Moreover, all the
partial derivatives in $x$ of $\bar{g}^{i}\left(  \cdot,\cdot\right)  $ up to
order $m$ are uniformly bounded.
\end{theorem}

\bop First define for each $1\leq i\leq n$,%
\begin{align*}
\tilde{g}^{i}\left(  t,x\right)   &  :=\bar{g}^{i}\left(  t,x\right)
-x_{i},\\
\tilde{\eta}^{i}\left(  t,x\right)   &  :=\bar{\eta}^{i}\left(  t,x\right)  .
\end{align*}
Accordingly one has, for each $1\leq i\leq n$,%
\begin{align*}
\tilde{g}_{x}^{i}\left(  t,x\right)   &  \equiv\bar{g}_{x}^{i}\left(
t,x\right)  -\mathbf{e}_{i},\\
\tilde{g}_{xx}^{i}\left(  t,x\right)   &  \equiv\bar{g}_{xx}^{i}\left(
t,x\right)  ,\\
\tilde{\eta}_{x}^{i}\left(  t,x\right)   &  \equiv\bar{\eta}_{x}^{i}\left(
t,x\right)  ,
\end{align*}
where $\mathbf{e}_{i}=\left(  \delta_{i,1},\delta_{i,2},...,\delta
_{i,n}\right)  ^{\top}$, with $\delta_{i,j}=\left\{
\begin{array}
[c]{l}%
1\text{ if }i=j,\\
0\text{ if }i\neq j.
\end{array}
\right.  .$

Then it is not difficult to see that $\left(  \bar{g}^{i}\left(  \cdot
,\cdot\right)  ,\bar{\eta}^{i}\left(  \cdot,\cdot\right)  \right)  $ is a
classical solution of (\ref{Eq13**}), if and only if, $\left(  \tilde{g}%
^{i}\left(  \cdot,\cdot\right)  ,\tilde{\eta}^{i}\left(  \cdot,\cdot\right)
\right)  $ is a classical solution of the following BSPDE%
\begin{equation}
\left\{
\begin{array}
[c]{l}%
d\tilde{g}^{i}\left(  t,x\right)  =-\left\{  \left\langle \tilde{g}_{x}%
^{i}\left(  t,x\right)  ,b^{\bar{u}\left(  \cdot\right)  }\left(  t,x\right)
\right\rangle +\frac{1}{2}\text{\textbf{tr}}\left[  \left(  \sigma\sigma
^{\top}\right)  ^{u\left(  \cdot\right)  }\left(  t,x\right)  \tilde{g}%
_{xx}^{i}\left(  t,x\right)  \right]  \right. \\
\text{ \ \ \ \ \ \ \ \ \ \ \ \ \ \ \ \ \ }+\left.  \text{\textbf{tr}}\left[
\tilde{\eta}_{x}^{i}\left(  t,x\right)  \sigma\left(  t,x\right)  \right]
+\left(  b^{i}\right)  ^{\bar{u}\left(  \cdot\right)  }\left(  t,x\right)
\right\}  dt\\
\text{ \ \ \ \ \ \ \ \ \ \ \ \ \ \ \ \ \ }+\tilde{\eta}^{i}\left(  t,x\right)
^{\top}dW\left(  t\right)  \text{, \ for }\left(  t,x\right)  \in\left[
0,T\right]  \times\mathbb{R}^{n}\text{,}\\
\tilde{g}^{i}\left(  T,x\right)  =0\text{, for }x\in\mathbb{R}^{n}\text{,}%
\end{array}
\right.  \label{Eq13***}%
\end{equation}
where $\left(  b^{i}\right)  ^{\bar{u}\left(  \cdot\right)  }\left(
t,x\right)  =\left\langle \mathbf{e}_{i},b^{\bar{u}\left(  \cdot\right)
}\left(  t,x\right)  \right\rangle $ is the i-th coordinate of $b^{\bar
{u}\left(  \cdot\right)  }\left(  t,x\right)  .$

According to Theorem \ref{result1}, under Assumption \textbf{(H}$_{m}%
$\textbf{)} with $m>2+\frac{n}{2}$, the BSPDE (\ref{Eq13***}) admits a unique
adapted classical solution $\left(  \tilde{g}^{i}\left(  \cdot,\cdot\right)
,\tilde{\eta}^{i}\left(  \cdot,\cdot\right)  \right)  $ and all the partial
derivatives in $x$ of $\tilde{g}^{i}\left(  \cdot,\cdot\right)  $ up to order
$m$ are uniformly bounded. Therefor, the BSPDE (\ref{Eq13**}) admits a unique
adapted classical solution $\left(  \bar{g}^{i}\left(  \cdot,\cdot\right)
,\bar{\eta}^{i}\left(  \cdot,\cdot\right)  \right)  $.\eop

In the next, we define for any $\bar{x}\in%
\mathbb{R}
^{n},$%
\[
G_{\bar{x}_{i}}\left(  \bar{x}\right)  :=\frac{\partial}{\partial\bar{x}_{i}%
}G\left(  \bar{x}_{1},...,\bar{x}_{n}\right)  .
\]

The following theorem constitutes the first main contribution of the paper.

\begin{theorem}
[Stochastic Maximum Principle]\label{result04}Let Assumption \textbf{(H}$_{m}%
$\textbf{)} be satisfied with $m>2+\frac{n}{2}$. If $\left(  \bar{u}\left(
\cdot\right)  ,\bar{X}\left(  \cdot\right)  \right)  $ is an optimal pair of
Problem (S), then the BSPDEs (\ref{Eq13}) and (\ref{Eq13**}), for $1\leq i\leq
n$, admit the classical solutions $\left(  \bar{\theta}\left(  \cdot
,\cdot\right)  ,\bar{\psi}\left(  \cdot,\cdot\right)  \right)  $ and $\left(
\bar{g}^{i}\left(  \cdot,\cdot\right)  ,\bar{\eta}^{i}\left(  \cdot
,\cdot\right)  \right)  $, respectively, such that the following holds%
\begin{align}
0  &  \leq\left\langle b\left(  \tau,\bar{X}\left(  \tau\right)  ,u\right)
-b\left(  \tau,\bar{X}\left(  \tau\right)  ,\bar{u}\left(  \tau\right)
\right)  ,\bar{\theta}_{x}\left(  \tau,\bar{X}\left(  \tau\right)  \right)
\right\rangle \nonumber\\
&  +%
{\textstyle\sum_{i=1}^{n}}
G_{\bar{x}_{i}}\left(  \mathbb{E}\left[  \bar{X}\left(  T\right)  \right]
\right)  \left\langle b\left(  \tau,\bar{X}\left(  \tau\right)  ,u\right)
-b\left(  \tau,\bar{X}\left(  \tau\right)  ,\bar{u}\left(  \tau\right)
\right)  ,\bar{g}_{x}^{i}\left(  \tau,\bar{X}\left(  \tau\right)  \right)
\right\rangle \nonumber\\
&  +%
\begin{array}
[c]{c}%
f\left(  \tau,\bar{X}\left(  \tau\right)  ,u\right)  -f\left(  \tau,\bar
{X}\left(  \tau\right)  ,\bar{u}\left(  \tau\right)  \right)  ,\text{ for all
}u\in U,\text{ a.s., a.e. }\tau\in\left[  0,T\right]  ,
\end{array}
\label{Eq14}%
\end{align}
or, equivalently,%
\begin{align}
&  \left\langle b\left(  \tau,\bar{X}\left(  \tau\right)  ,\bar{u}\left(
\tau\right)  \right)  ,\bar{\theta}_{x}\left(  \tau,\bar{X}\left(
\tau\right)  \right)  +%
{\textstyle\sum_{i=1}^{n}}
G_{\bar{x}_{i}}\left(  \mathbb{E}\left[  \bar{X}\left(  T\right)  \right]
\right)  \bar{g}_{x}^{i}\left(  \tau,\bar{X}\left(  \tau\right)  \right)
\right\rangle \nonumber\\
&  +f\left(  \tau,\bar{X}\left(  \tau\right)  ,\bar{u}\left(  \tau\right)
\right) \nonumber\\
&  =\min\limits_{u\in U}\left\{  \left\langle b\left(  \tau,\bar{X}\left(
\tau\right)  ,u\right)  ,\bar{\theta}_{x}\left(  \tau,\bar{X}\left(
\tau\right)  \right)  +%
{\textstyle\sum_{i=1}^{n}}
G_{\bar{x}_{i}}\left(  \mathbb{E}\left[  \bar{X}\left(  T\right)  \right]
\right)  \bar{g}_{x}^{i}\left(  \tau,\bar{X}\left(  \tau\right)  \right)
\right\rangle \right. \nonumber\\
&
\begin{array}
[c]{c}%
\left.  +f\left(  \tau,\bar{X}\left(  \tau\right)  ,u\right)  \right\}
,\text{ a.e. }\tau\in\left[  0,T\right]  \text{, a.s,.}%
\end{array}
\label{Eq15}%
\end{align}
Furthermore the objective value of $\bar{u}\left(  \cdot\right)  $ is given by%
\begin{equation}
\mathbf{J}\left(  \bar{u}\left(  \cdot\right)  \right)  =\bar{\theta}\left(
0,x_{0}\right)  +G\left(  \bar{g}\left(  0,x_{0}\right)  \right)  .
\label{Eq16}%
\end{equation}

\end{theorem}

A proof of the above theorem will be carried out in the following section.
Analogous to the classical stochastic maximum principle, the optimal state
equation, the corresponding BSPDEs (\ref{Eq13})- (\ref{Eq13**}), along with
the minimum condition (\ref{Eq15}), can be written as the following system of
FBSPDEs:%
\begin{equation}
\left\{
\begin{array}
[c]{l}%
d\bar{X}\left(  t\right)  =b\left(  t,\bar{X}\left(  t\right)  ,\bar{u}\left(
t\right)  \right)  dt+\sigma\left(  t,\bar{X}\left(  t\right)  \right)
dW\left(  t\right)  ,\ \text{for }t\in\left[  0,T\right]  ,\\
d\bar{\theta}\left(  t,x\right)  =-\left\{  \left\langle \bar{\theta}%
_{x}\left(  t,x\right)  ,b\left(  t,x,\bar{u}\left(  t\right)  \right)
\right\rangle \right.  +\frac{1}{2}\text{\textbf{tr}}\left[  \sigma\left(
t,x\right)  \sigma\left(  t,x\right)  ^{\top}\bar{\theta}_{xx}\left(
t,x\right)  \right] \\
\ \ \ \ \ \ \ \ \ \ \ \ \ \ \ +\left.  \text{\textbf{tr}}\left[  \bar{\psi
}_{x}\left(  t,x\right)  \sigma\left(  t,x\right)  \right]  +f\left(
t,x,\bar{u}\left(  t\right)  \right)  \right\}  dt\\
\ \ \ \ \ \ \ \ \ \ \ \ \ \ \ +\bar{\psi}\left(  t,x\right)  ^{\top}dW\left(
t\right)  ,\text{ for }\left(  t,x\right)  \in\left[  0,T\right]
\times\mathbb{R}^{n}\text{,}\\
d\bar{g}^{i}\left(  t,x\right)  =-\left\{  \left\langle \bar{g}_{x}^{i}\left(
t,x\right)  ,b^{\hat{u}\left(  \cdot\right)  }\left(  t,x\right)
\right\rangle +\text{\textbf{tr}}\left[  \sigma\left(  t,x\right)
\sigma\left(  t,x\right)  ^{\top}\bar{g}_{xx}^{i}\left(  t,x\right)  \right]
\right. \\
\text{ \ \ \ \ \ \ \ \ \ \ \ \ }%
\begin{array}
[c]{r}%
+\left.  \text{\textbf{tr}}\left[  \bar{\eta}_{x}^{i}\left(  t,x\right)
\sigma\left(  t,x\right)  \right]  \right\}  dt+\bar{\eta}^{i}\left(
t,x\right)  ^{\top}dW\left(  t\right)  \text{,}\\
\text{for }\left(  t,x\right)  \in\left[  0,T\right]  \times\mathbb{R}%
^{n}\text{, }%
\end{array}
\\
\bar{g}^{i}\left(  T,x\right)  =x_{i}\text{, for }x\in\mathbb{R}^{n}\text{,
}1\leq i\leq n,\\
\bar{X}\left(  0\right)  =x_{0},\text{ }\bar{\theta}\left(  T,x\right)
=h\left(  x\right)  ,\text{ for }x\in\mathbb{R}^{n},\\
\bar{u}\left(  t\right)  \in\arg\min\limits_{u\in U}\left\{  \left\langle
b\left(  t,\bar{X}\left(  t\right)  ,u\right)  ,\bar{\theta}_{x}\left(
t,\bar{X}\left(  t\right)  \right)  \right\rangle +f\left(  t,\bar{X}\left(
t\right)  ,u\right)  \right. \\
\text{ \ \ \ \ \ \ \ \ }%
\begin{array}
[c]{r}%
\left.  +%
{\textstyle\sum_{i=1}^{n}}
\left\langle b\left(  t,\bar{X}\left(  t\right)  ,u\right)  ,G_{\bar{x}_{i}%
}\left(  \mathbb{E}\left[  \bar{X}\left(  T\right)  \right]  \right)  \bar
{g}_{x}^{i}\left(  t,\bar{X}\left(  t\right)  \right)  \right\rangle \right\}
,\text{ }\\
\text{a.e. t}\in\left[  0,T\right]  \text{, a.s,}%
\end{array}
\end{array}
\right.  \label{Eq25}%
\end{equation}

It should be noted that, the above system of FBSPDEs is quite different from
the traditional stochastic (extended) Hamiltonian system (see e.g. Yong and
Zhou \cite{YongZhou}), since two backward stochastic partial differential
equations are involved. We can rephrase Theorem \ref{result04} as the following.

\begin{theorem}
\label{result05}Let Assumption \textbf{(H}$_{m}$\textbf{)} be satisfied with
$m>2+\frac{n}{2}$. If $\left(  \bar{u}\left(  \cdot\right)  ,\bar{X}\left(
\cdot\right)  \right)  $ is an optimal pair of Problem (S), then the family of
processes (resp. random fields) $\left(  \bar{X}\left(  \cdot\right)
\text{,}\bar{u}\left(  \cdot\right)  \text{,}\left\{  \bar{g}^{i}\left(
\cdot,\cdot\right)  ,\bar{\eta}^{i}\left(  \cdot,\cdot\right)  \right\}
_{i=1,...,n}\text{,}\bar{\theta}\left(  \cdot,\cdot\right)  \text{,}\bar{\psi
}\left(  \cdot,\cdot\right)  \right)  $ satisfies the system of FBSPDEs
(\ref{Eq25}). Furthermore the objective value of $\bar{u}\left(  \cdot\right)
$ is given by (\ref{Eq16}).
\end{theorem}

\begin{remark}
On comparing between the SMP in Theorem \ref{result04} and the traditional SMP
of Bensoussan \cite{Ben} we find the following facts:

\begin{enumerate}
\item[(i)] The advantage: The new version of the SMP enables us to derive
simultaneously an optimal solution $\bar{u}\left(  \cdot\right)  $ as well as
its objective value $\mathbf{J}\left(  \bar{u}\left(  \cdot\right)  \right)
$\ in terms of the solution of FBSPDEs (\ref{Eq25}), while the classical SMP
approach permits us to derive the optimal solution $\bar{u}\left(
\cdot\right)  $ only.\newline

\item[(ii)] The disadvantage: Of course, the disadvantage of the new version
of the SMP lies in the fact of assuming differentiability of higher orders in
$x$ on the coefficients, which unfortunately limits the scope of problems
applicable to the theorem.
\end{enumerate}
\end{remark}

\begin{remark}
If the function $U\ni u\rightarrow\left\langle b\left(  t,x,u\right)
,p\right\rangle +f\left(  t,x,u\right)  \in%
\mathbb{R}
$ has a unique minimizer $\mathbf{\bar{u}}\left(  t,x,p\right)  $ for each
$t\in\left[  0,T\right]  $ and $\left(  x,p\right)  \in%
\mathbb{R}
^{n}\times%
\mathbb{R}
^{n}$ such that $\mathbf{\bar{u}}\left(  t,x,p\right)  $ satisfies an
appropriate regularity condition, then the optimal control is characterized by%
\[
\bar{u}\left(  t\right)  \equiv\mathbf{\bar{u}}\left(  t,\bar{X}\left(
t\right)  ,\bar{\theta}_{x}\left(  t,\bar{X}\left(  t\right)  \right)  +%
{\textstyle\sum_{i=1}^{n}}
G_{\bar{x}_{i}}\left(  \mathbb{E}\left[  \bar{X}\left(  T\right)  \right]
\right)  \bar{g}_{x}^{i}\left(  t,\bar{X}\left(  t\right)  \right)  \right)
\text{,}%
\]
where $\left(  \bar{X}\left(  \cdot\right)  ,\left\{  \bar{g}^{i}\left(
\cdot,\cdot\right)  ,\bar{\eta}^{i}\left(  \cdot,\cdot\right)  \right\}
_{i=1,...,n},\bar{\theta}\left(  \cdot,\cdot\right)  ,\bar{\psi}\left(
\cdot,\cdot\right)  \right)  $ is family of processes (resp. random fields)
that satisfies the following system of coupled FBSPDEs:%
\begin{equation}
\left\{
\begin{array}
[c]{l}%
d\bar{X}\left(  t\right)  =b\left(  t,\bar{X}\left(  t\right)  ,\mathrm{\tilde
{u}}\left(  t,\bar{X}\left(  t\right)  \right)  \right)  dt+\sigma\left(
t,\bar{X}\left(  t\right)  \right)  dW\left(  t\right)  ,\ \text{for }%
t\in\left[  0,T\right]  ,\\
d\bar{\theta}\left(  t,x\right)  =-\left\{  \left\langle \bar{\theta}%
_{x}\left(  t,x\right)  ,b\left(  t,x,\mathrm{\tilde{u}}\left(  t,\bar
{X}\left(  t\right)  \right)  \right)  \right\rangle +\frac{1}{2}%
\text{\textbf{tr}}\left[  \sigma\left(  t,x\right)  \sigma\left(  t,x\right)
^{\top}\bar{\theta}_{xx}\left(  t,x\right)  \right]  \right. \\
\ \ \ \ \ \ \ \ \ \ \ \ \ +\left.  \text{\textbf{tr}}\left[  \bar{\psi}%
_{x}\left(  t,x\right)  \sigma\left(  t,x\right)  \right]  +f\left(
t,x,\mathrm{\tilde{u}}\left(  t,\bar{X}\left(  t\right)  \right)  \right)
\right\}  dt\\
\ \ \ \ \ \ \ \ \ \ \ \ \ \ +\bar{\psi}\left(  t,x\right)  ^{\top}dW\left(
t\right)  ,\text{ for }\left(  t,x\right)  \in\left[  0,T\right]
\times\mathbb{R}^{n}\text{,}\\
d\bar{g}^{i}\left(  t,x\right)  =-\left\{  \left\langle \bar{g}_{x}^{i}\left(
t,x\right)  ,b\left(  t,x,\mathrm{\tilde{u}}\left(  t,\bar{X}\left(  t\right)
\right)  \right)  \right\rangle +\text{\textbf{tr}}\left[  \sigma\left(
t,x\right)  \sigma\left(  t,x\right)  ^{\top}\bar{g}_{xx}^{i}\left(
t,x\right)  \right]  \right. \\
\text{ \ \ \ \ \ \ \ \ \ \ \ \ }%
\begin{array}
[c]{r}%
+\left.  \text{\textbf{tr}}\left[  \bar{\eta}_{x}^{i}\left(  t,x\right)
\sigma\left(  t,x\right)  \right]  \right\}  dt+\bar{\eta}^{i}\left(
t,x\right)  ^{\top}dW\left(  t\right)  \text{,}\\
\text{for }\left(  t,x\right)  \in\left[  0,T\right]  \times\mathbb{R}%
^{n}\text{, }%
\end{array}
\\
\bar{g}^{i}\left(  T,x\right)  =x_{i}\text{, for }x\in\mathbb{R}^{n}\text{,
}1\leq i\leq n,\\
\bar{X}\left(  0\right)  =x_{0},\text{ }\bar{\theta}\left(  T,x\right)
=h\left(  x\right)  ,\text{ for }x\in\mathbb{R}^{n},\\
\mathrm{\tilde{u}}\left(  t,\bar{X}\left(  t\right)  \right)  :=\mathbf{\bar
{u}}\left(  t,\bar{X}\left(  t\right)  ,\bar{\theta}_{x}\left(  t,\bar
{X}\left(  t\right)  \right)  +%
{\textstyle\sum_{i=1}^{n}}
G_{\bar{x}_{i}}\left(  \mathbb{E}\left[  \bar{X}\left(  T\right)  \right]
\right)  \bar{g}_{x}^{i}\left(  t,\bar{X}\left(  t\right)  \right)  \right)  .
\end{array}
\right.  \label{Eq326}%
\end{equation}

\end{remark}

It is well known that optimal control theory can be used to solve stochastic
(extended) Hamiltonian systems, which are actually coupled forward-backward
stochastic differential equation (FBSDE, for short); see e.g. Yong and Zhou
\cite{YongZhou}. However, it is seen from Theorem \ref{result05} that
stochastic control theory can also be used to solve coupled equations of the
type (\ref{Eq326}). In Section \ref{section7} we will study more about this
type of equations.

\section{A Proof of Theorem 4.2\label{section5}}

This section is devoted to the proof of Theorem \ref{result04}. Let $\left(
\bar{u}\left(  \cdot\right)  ,\bar{X}\left(  \cdot\right)  \right)  $ be the
given optimal pair. Then the following is satisfied%
\[
\left\{
\begin{array}
[c]{l}%
d\bar{X}\left(  t\right)  =b^{\bar{u}\left(  \cdot\right)  }\left(  t,\bar
{X}\left(  t\right)  \right)  dt+\sigma\left(  t,\bar{X}\left(  t\right)
\right)  dW\left(  t\right)  \text{, }t\in\left[  0,T\right]  \text{,}\\
\bar{X}\left(  0\right)  =x_{0}\text{.}%
\end{array}
\right.
\]

For any $\tau\in\left[  0,T\right)  $, $u\left(  \cdot\right)  \in
\mathcal{U}\left[  0,T\right]  $ and for any $\varepsilon\in\left[
0,T-\tau\right)  $, define%
\begin{equation}
u^{\varepsilon}\left(  t\right)  =\left\{
\begin{array}
[c]{l}%
u\left(  t\right)  ,\text{ for }t\in\left[  \tau,\tau+\varepsilon\right)
\text{,}\\
\bar{u}\left(  t\right)  ,\text{ for }t\in\left[  0,T\right]  /\left[
\tau,\tau+\varepsilon\right)  \text{.}%
\end{array}
\right.  \label{Eq20}%
\end{equation}

Let $\left(  u^{\varepsilon}\left(  \cdot\right)  ,X^{\varepsilon}\left(
\cdot\right)  \right)  $ satisfy the following%
\begin{equation}
\left\{
\begin{array}
[c]{l}%
dX^{\varepsilon}\left(  t\right)  =b^{u^{\varepsilon}\left(  \cdot\right)
}\left(  t,X^{\varepsilon}\left(  t\right)  \right)  dt+\sigma\left(
t,X^{\varepsilon}\left(  t\right)  \right)  dW\left(  t\right)  \text{, }%
t\in\left[  0,T\right]  \text{,}\\
X^{\varepsilon}\left(  0\right)  =x_{0}\text{.}%
\end{array}
\right.  \label{Eq21}%
\end{equation}

To prove Theorem \ref{result04} we need to some preliminary results given in
the following two lemmas.

\begin{lemma}
\label{result06}Suppose that Assumption \textbf{(H}$_{m}$\textbf{)} be
satisfied with $m>2+\frac{n}{2}.$ Let $\left(  \bar{\theta}\left(  \cdot
,\cdot\right)  ,\bar{\psi}\left(  \cdot,\cdot\right)  \right)  $ and
\newline$\left(  \bar{g}^{i}\left(  \cdot,\cdot\right)  ,\bar{\eta}^{i}\left(
\cdot,\cdot\right)  \right)  $, for $1\leq i\leq n,$ be the unique classical
solutions of (\ref{Eq13}) and (\ref{Eq13**}), respectively. Then for each
$t\in\left[  0,T\right]  $ and $1\leq i\leq n$, $\bar{\theta}\left(  t,\bar
{X}\left(  t\right)  \right)  $ and $\bar{g}^{i}\left(  t,\bar{X}\left(
t\right)  \right)  $ have the following probabilistic representations:%
\begin{equation}
\bar{\theta}\left(  t,\bar{X}\left(  t\right)  \right)  =\mathbb{E}_{t}\left[
\int_{t}^{T}f\left(  s,\bar{X}\left(  s\right)  ,\bar{u}\left(  s\right)
\right)  ds+F\left(  \bar{X}\left(  T\right)  \right)  \right]  \text{,}
\label{Eq17}%
\end{equation}
and%
\begin{equation}
\bar{g}^{i}\left(  t,\bar{X}\left(  t\right)  \right)  =\mathbb{E}_{t}\left[
\bar{X}_{i}\left(  T\right)  \right]  \label{Eq17*}%
\end{equation}
where $\bar{X}_{i}\left(  T\right)  $ denotes the i-th coordinate of $\bar
{X}\left(  T\right)  =\left(  \bar{X}_{1}\left(  T\right)  ,...,\bar{X}%
_{n}\left(  T\right)  \right)  ^{\top}$. Furthermore, for any $\tau\in\left[
0,T\right]  $, $u\left(  \cdot\right)  \in\mathcal{U}\left[  0,T\right]  $ and
for any $\varepsilon\in\left[  0,T-\tau\right)  $, the following equality
holds%
\begin{align}
&  \mathbf{J}\left(  u^{\varepsilon}\left(  \cdot\right)  \right)
-\mathbf{J}\left(  \bar{u}\left(  \cdot\right)  \right) \nonumber\\
&  =\mathbb{E}\left[  \int_{\tau}^{\tau+\varepsilon}\left\{  \left\langle
b^{u\left(  \cdot\right)  }\left(  t,X^{\varepsilon}\left(  t\right)  \right)
-b^{\bar{u}\left(  \cdot\right)  }\left(  t,X^{\varepsilon}\left(  t\right)
\right)  ,\bar{\theta}_{x}\left(  t,X^{\varepsilon}\left(  t\right)  \right)
\right\rangle \right.  \right. \nonumber\\
&  +%
{\textstyle\sum_{i=1}^{n}}
\left\langle b^{u\left(  \cdot\right)  }\left(  t,X^{\varepsilon}\left(
t\right)  \right)  -b^{\bar{u}\left(  \cdot\right)  }\left(  t,X^{\varepsilon
}\left(  t\right)  \right)  ,G_{\bar{x}_{i}}\left(  \mathbb{E}\left[  \bar
{g}\left(  t,X^{\varepsilon}\left(  t\right)  \right)  \right]  \right)
\bar{g}_{x}^{i}\left(  t,X^{\varepsilon}\left(  t\right)  \right)
\right\rangle \nonumber\\
&  \left.  \left.  +f^{u\left(  \cdot\right)  }\left(  t,X^{\varepsilon
}\left(  t\right)  \right)  -f^{\bar{u}\left(  \cdot\right)  }\left(
t,X^{\varepsilon}\left(  t\right)  \right)  \right\}  dt\right]  .
\label{Eq18}%
\end{align}

\end{lemma}

\bop Equalities (\ref{Eq17})-(\ref{Eq17*}) can be easily obtained by applying
Theorem \ref{result02}. So we only need to show (\ref{Eq18}). Consider the
difference%
\begin{align*}
&  \mathbf{J}\left(  u^{\varepsilon}\left(  \cdot\right)  \right)
-\mathbf{J}\left(  \bar{u}\left(  \cdot\right)  \right) \\
&  =\mathbb{E}\left[  \int_{0}^{T}f^{u^{\varepsilon}\left(  \cdot\right)
}\left(  t,X^{\varepsilon}\left(  t\right)  \right)  dt+h\left(
X^{\varepsilon}\left(  T\right)  \right)  +G\left(  \mathbb{E}\left[
X^{\varepsilon}\left(  T\right)  \right]  \right)  \right] \\
&  -\mathbb{E}\left[  \int_{0}^{T}f^{\bar{u}\left(  \cdot\right)  }\left(
t,\bar{X}\left(  t\right)  \right)  dt+h\left(  \bar{X}\left(  T\right)
\right)  +G\left(  \mathbb{E}\left[  \bar{X}\left(  T\right)  \right]
\right)  \right]  .
\end{align*}
By the terminal conditions in BSPDEs (\ref{Eq13})-(\ref{Eq13**}) and by
setting $t=0$ in (\ref{Eq17})-(\ref{Eq17*}), we obtain that%
\begin{align*}
\mathbf{J}\left(  u^{\varepsilon}\left(  \cdot\right)  \right)   &
=\mathbb{E}\left[  \int_{0}^{T}f^{u^{\varepsilon}\left(  \cdot\right)
}\left(  t,X^{\varepsilon}\left(  t\right)  \right)  dt+\bar{\theta}\left(
T,X^{\varepsilon}\left(  T\right)  \right)  \right] \\
&  +G\left(  \mathbb{E}\left[  \bar{g}^{1}\left(  T,X^{\varepsilon}\left(
T\right)  \right)  \right]  ,...,\mathbb{E}\left[  \bar{g}^{n}\left(
T,X^{\varepsilon}\left(  T\right)  \right)  \right]  \right)
\end{align*}
and%
\[
\mathbf{J}\left(  \bar{u}\left(  \cdot\right)  \right)  =\bar{\theta}\left(
0,x_{0}\right)  +G\left(  \mathbb{E}\left[  \bar{g}^{1}\left(  0,x_{0}\right)
\right]  ,...,\mathbb{E}\left[  \bar{g}^{n}\left(  0,x_{0}\right)  \right]
\right)  .
\]
Accordingly, we have%
\begin{align}
&  \mathbf{J}\left(  u^{\varepsilon}\left(  \cdot\right)  \right)
-\mathbf{J}\left(  \bar{u}\left(  \cdot\right)  \right) \nonumber\\
&  =\mathbb{E}\left[  \int_{0}^{T}f^{u^{\varepsilon}\left(  \cdot\right)
}\left(  t,X^{\varepsilon}\left(  t\right)  \right)  dt+\bar{\theta}\left(
T,X^{\varepsilon}\left(  T\right)  \right)  -\bar{\theta}\left(
0,x_{0}\right)  \right] \nonumber\\
&  +G\left(  \mathbb{E}\left[  \bar{g}^{1}\left(  T,X^{\varepsilon}\left(
T\right)  \right)  \right]  ,...,\mathbb{E}\left[  \bar{g}^{n}\left(
T,X^{\varepsilon}\left(  T\right)  \right)  \right]  \right) \nonumber\\
&  \left.  -G\left(  \mathbb{E}\left[  \bar{g}^{1}\left(  0,x_{0}\right)
\right]  ,...,\mathbb{E}\left[  \bar{g}^{n}\left(  0,x_{0}\right)  \right]
\right)  \right]  . \label{Eq23}%
\end{align}
Recall that $X^{\varepsilon}\left(  0\right)  =x_{0}.$ Then by applying
It\^{o}-Wentzell formula to $\bar{\theta}\left(  t,X^{\varepsilon}\left(
t\right)  \right)  $ on time interval $\left[  0,T\right]  $, we get%
\begin{align}
&  \mathbb{E}\left[  \bar{\theta}\left(  T,X^{\varepsilon}\left(  T\right)
\right)  -\bar{\theta}\left(  0,x_{0}\right)  \right] \nonumber\\
&  =\mathbb{E}\left[  \bar{\theta}\left(  T,X^{\varepsilon}\left(  T\right)
\right)  -\bar{\theta}\left(  0,X^{\varepsilon}\left(  0\right)  \right)
\right] \nonumber\\
&  =\mathbb{E}\left[  \int_{0}^{T}\left\{  \left\langle \bar{\theta}%
_{x}\left(  t,X^{\varepsilon}\left(  t\right)  \right)  ,\bar{b}%
^{u^{\varepsilon}\left(  \cdot\right)  }\left(  t,X^{\varepsilon}\left(
t\right)  \right)  -\bar{b}^{\bar{u}\left(  \cdot\right)  }\left(
t,X^{\varepsilon}\left(  t\right)  \right)  \right\rangle \right.  \right.
\nonumber\\
&  \left.  \left.  -f^{\bar{u}\left(  \cdot\right)  }\left(  t,X^{\varepsilon
}\left(  t\right)  \right)  \right\}  dt\right]  . \label{Eq24}%
\end{align}
and it follows from the chain rule applied to $G\left(  \mathbb{E}\left[
\bar{g}^{1}\left(  t,X^{\varepsilon}\left(  t\right)  \right)  \right]
,...,\mathbb{E}\left[  \bar{g}^{n}\left(  t,X^{\varepsilon}\left(  t\right)
\right)  \right]  \right)  $ that%
\begin{align}
&  G\left(  \mathbb{E}\left[  \bar{g}^{1}\left(  T,X^{\varepsilon}\left(
T\right)  \right)  \right]  ,...,\mathbb{E}\left[  \bar{g}^{n}\left(
T,X^{\varepsilon}\left(  T\right)  \right)  \right]  \right) \nonumber\\
&  -G\left(  \mathbb{E}\left[  \bar{g}^{1}\left(  0,x_{0}\right)  \right]
,...,\mathbb{E}\left[  \bar{g}^{n}\left(  0,x_{0}\right)  \right]  \right)
\nonumber\\
&  =%
{\textstyle\sum_{i=1}^{n}}
\int_{0}^{T}G_{\bar{x}_{i}}\left(  \mathbb{E}\left[  \bar{g}^{1}\left(
t,X^{\varepsilon}\left(  t\right)  \right)  \right]  ,...,\mathbb{E}\left[
\bar{g}^{n}\left(  t,X^{\varepsilon}\left(  t\right)  \right)  \right]
\right)  d\mathbb{E}\left[  \bar{g}^{i}\left(  t,X^{\varepsilon}\left(
t\right)  \right)  \right] \nonumber\\
&  =%
{\textstyle\sum_{i=1}^{n}}
\int_{0}^{T}G_{\bar{x}_{i}}\left(  \mathbb{E}\left[  \bar{g}\left(
t,X^{\varepsilon}\left(  t\right)  \right)  \right]  \right)  d\mathbb{E}%
\left[  \bar{g}^{i}\left(  t,X^{\varepsilon}\left(  t\right)  \right)
\right]  , \label{Eq24*}%
\end{align}
where%
\begin{align*}
&  d\mathbb{E}\left[  \bar{g}^{i}\left(  t,X^{\varepsilon}\left(  t\right)
\right)  \right] \\
&  =\mathbb{E}\left[  \left\langle \bar{g}_{x}^{i}\left(  t,X^{\varepsilon
}\left(  t\right)  \right)  ,b^{u^{\varepsilon}\left(  \cdot\right)  }\left(
t,X^{\varepsilon}\left(  t\right)  \right)  -b^{\bar{u}\left(  \cdot\right)
}\left(  t,X^{\varepsilon}\left(  t\right)  \right)  \right\rangle \right]
dt.
\end{align*}
Combining (\ref{Eq23}) together with (\ref{Eq24}) and (\ref{Eq24*}), it
follows that%
\begin{align*}
&  \mathbf{J}\left(  u^{\varepsilon}\left(  \cdot\right)  \right)
-\mathbf{J}\left(  \bar{u}\left(  \cdot\right)  \right) \\
&  =\mathbb{E}\left[  \int_{0}^{T}\left\{  \left\langle \bar{\theta}%
_{x}\left(  t,X^{\varepsilon}\left(  t\right)  \right)  ,b^{u^{\varepsilon
}\left(  \cdot\right)  }\left(  t,X^{\varepsilon}\left(  t\right)  \right)
-b^{\bar{u}\left(  \cdot\right)  }\left(  t,X^{\varepsilon}\left(  t\right)
\right)  \right\rangle \right.  \right. \\
&  +\sum\limits_{i=1}^{n}G_{\bar{x}_{i}}\left(  \mathbb{E}\left[  \bar
{g}\left(  t,X^{\varepsilon}\left(  t\right)  \right)  \right]  \right)
\left\langle \bar{g}_{x}^{i}\left(  t,X^{\varepsilon}\left(  t\right)
\right)  ,b^{u^{\varepsilon}\left(  \cdot\right)  }\left(  t,X^{\varepsilon
}\left(  t\right)  \right)  -b^{\bar{u}\left(  \cdot\right)  }\left(
t,X^{\varepsilon}\left(  t\right)  \right)  \right\rangle \\
&  \left.  \left.  +f^{u^{\varepsilon}\left(  \cdot\right)  }\left(
t,X^{\varepsilon}\left(  t\right)  \right)  -f^{\bar{u}\left(  \cdot\right)
}\left(  t,X^{\varepsilon}\left(  t\right)  \right)  \right\}  dt\right]  .
\end{align*}
Since $u^{\varepsilon}\left(  t\right)  =u\left(  t\right)  \chi_{\left[
t,t+\varepsilon\right)  }\left(  t\right)  +\bar{u}\left(  t\right)
\chi_{\left[  0,T\right]  /\left[  \tau,\tau+\varepsilon\right)  }\left(
t\right)  ,$ we obtain that%
\begin{align*}
&  \mathbf{J}\left(  u^{\varepsilon}\left(  \cdot\right)  \right)
-\mathbf{J}\left(  \bar{u}\left(  \cdot\right)  \right) \\
&  =\mathbb{E}\left[  \int_{\tau}^{\tau+\varepsilon}\left\{  \left\langle
\bar{\theta}_{x}\left(  t,X^{\varepsilon}\left(  t\right)  \right)
,b^{u\left(  \cdot\right)  }\left(  t,X^{\varepsilon}\left(  t\right)
\right)  -b^{\bar{u}\left(  \cdot\right)  }\left(  t,X^{\varepsilon}\left(
t\right)  \right)  \right\rangle \right.  \right. \\
&  +\sum\limits_{i=1}^{n}G_{\bar{x}_{i}}\left(  \mathbb{E}\left[  \bar
{g}\left(  t,X^{\varepsilon}\left(  t\right)  \right)  \right]  \right)
\left\langle \bar{g}_{x}^{i}\left(  t,X^{\varepsilon}\left(  t\right)
\right)  ,b^{u\left(  \cdot\right)  }\left(  t,X^{\varepsilon}\left(
t\right)  \right)  -b^{\bar{u}\left(  \cdot\right)  }\left(  t,X^{\varepsilon
}\left(  t\right)  \right)  \right\rangle \\
&  \left.  \left.  +f^{u\left(  \cdot\right)  }\left(  t,X^{\varepsilon
}\left(  t\right)  \right)  -f^{\bar{u}\left(  \cdot\right)  }\left(
t,X^{\varepsilon}\left(  t\right)  \right)  \right\}  dt\right]  .
\end{align*}
This completes the proof.\eop

\begin{lemma}
\label{result07}Suppose that Assumption \textbf{(H}$_{m}$\textbf{)} be
satisfied with $m>2+\frac{n}{2}.$ Let $\left(  \bar{\theta}\left(  \cdot
,\cdot\right)  ,\bar{\psi}\left(  \cdot,\cdot\right)  \right)  $ and
\newline$\left(  \bar{g}^{i}\left(  \cdot,\cdot\right)  ,\bar{\eta}^{i}\left(
\cdot,\cdot\right)  \right)  $, for $1\leq i\leq n,$ be the unique classical
solutions of (\ref{Eq13}) and (\ref{Eq13**}), respectively. Then for any
$\tau\in\left[  0,T\right]  $, $u\left(  \cdot\right)  \in\mathcal{U}\left[
0,T\right]  $ and for any $\varepsilon\in\left[  0,T-\tau\right)  $, the
following equality holds%
\begin{align*}
&  \mathbf{J}\left(  u^{\varepsilon}\left(  \cdot\right)  \right)
-\mathbf{J}\left(  \bar{u}\left(  \cdot\right)  \right) \\
&  =\mathbb{E}\left[  \int_{\tau}^{\tau+\varepsilon}\left\{  \left\langle
b^{u\left(  \cdot\right)  }\left(  t,\bar{X}\left(  t\right)  \right)
-b^{\bar{u}\left(  \cdot\right)  }\left(  t,\bar{X}\left(  t\right)  \right)
,\bar{\theta}_{x}\left(  t,\bar{X}\left(  t\right)  \right)  \right\rangle
\right.  \right. \\
&  +%
{\textstyle\sum_{i=1}^{n}}
\left\langle b^{u\left(  \cdot\right)  }\left(  t,\bar{X}\left(  t\right)
\right)  -b^{\bar{u}\left(  \cdot\right)  }\left(  t,\bar{X}\left(  t\right)
\right)  ,G_{\bar{x}_{i}}\left(  \mathbb{E}\left[  \bar{g}\left(  t,\bar
{X}\left(  t\right)  \right)  \right]  \right)  \bar{g}_{x}^{i}\left(
t,\bar{X}\left(  t\right)  \right)  \right\rangle \\
&  \left.  \left.  f^{u\left(  \cdot\right)  }\left(  t,\bar{X}\left(
t\right)  \right)  -f^{\bar{u}\left(  \cdot\right)  }\left(  t,\bar{X}\left(
t\right)  \right)  \right\}  dt\right]  +o\left(  \varepsilon\right)  .
\end{align*}

\end{lemma}

\bop Let $\xi^{\varepsilon}\left(  \cdot\right)  :=X^{\varepsilon}\left(
\cdot\right)  -\bar{X}\left(  \cdot\right)  $. Then we have $\xi^{\varepsilon
}\left(  \cdot\right)  $ satisfies the following SDE,%
\[
\left\{
\begin{array}
[c]{l}%
d\xi^{\varepsilon}\left(  t\right)  =\left\{  \tilde{b}^{\varepsilon}\left(
t\right)  \xi^{\varepsilon}\left(  t\right)  +\left\{  b^{u^{\varepsilon
}\left(  \cdot\right)  }\left(  t,\bar{X}\left(  t\right)  \right)
-b^{\bar{u}\left(  \cdot\right)  }\left(  t,\bar{X}\left(  t\right)  \right)
\right\}  \chi_{\left[  \tau,\tau+\varepsilon\right)  }\left(  t\right)
\right\}  dt\\
\text{ \ \ \ \ \ \ \ \ \ \ \ \ \ \ \ }+%
{\textstyle\sum\limits_{j=1}^{d}}
\left\{  \tilde{\sigma}^{j}\left(  t\right)  \xi^{\varepsilon}\left(
t\right)  \right\}  dW_{j}\left(  t\right)  \text{, }t\in\left[  0,T\right]
\text{,}\\
\xi^{\varepsilon}\left(  0\right)  =0\text{,}%
\end{array}
\right.
\]
where%
\[
\tilde{b}\left(  t\right)  :=\int_{0}^{1}b_{x}^{u^{\varepsilon}\left(
\cdot\right)  }\left(  t,\bar{X}\left(  t\right)  +r\xi^{\varepsilon}\left(
t\right)  \right)  dr
\]
and%
\[
\tilde{\sigma}^{j}\left(  t\right)  :=\int_{0}^{1}\sigma_{x}^{j}\left(
t,\bar{X}\left(  t\right)  +r\xi^{\varepsilon}\left(  t\right)  \right)  dr.
\]

By Lemma 4.2 in [\cite{YongZhou}, p. 124], we obtain for any $k\geq1$%
\begin{align*}
\sup_{t\in\left[  0,T\right]  }\mathbb{E}\left[  \left\vert \xi^{\varepsilon
}\left(  t\right)  \right\vert ^{2k}\right]   &  \leq K\left(  \int_{\tau
}^{\tau+\varepsilon}\mathbb{E}\left[  \left\vert b^{u\left(  \cdot\right)
}\left(  t,\bar{X}\left(  t\right)  \right)  -b^{\bar{u}\left(  \cdot\right)
}\left(  t,\bar{X}\left(  t\right)  \right)  \right\vert ^{2k}\right]
^{\frac{1}{2k}}dt\right)  ^{2k}\\
&  \leq K\varepsilon^{2k},
\end{align*}
for some $K>0.$

We now consider the difference%
\begin{align*}
\rho\left(  \varepsilon\right)   &  =\left\{  \mathbf{J}\left(  u^{\varepsilon
}\left(  \cdot\right)  \right)  -\mathbf{J}\left(  \bar{u}\left(
\cdot\right)  \right)  \right\} \\
&  -\mathbb{E}\left[  \int_{\tau}^{\tau+\varepsilon}\left\{  \left\langle
b^{u\left(  \cdot\right)  }\left(  t,\bar{X}\left(  t\right)  \right)
-b^{\bar{u}\left(  \cdot\right)  }\left(  t,\bar{X}\left(  t\right)  \right)
,\bar{\theta}_{x}\left(  t,\bar{X}\left(  t\right)  \right)  \right\rangle
\right.  \right. \\
&  +\left\langle b\left(  t,\bar{X}\left(  t\right)  ,u\left(  t\right)
\right)  -b\left(  t,\bar{X}\left(  t\right)  ,\bar{u}\left(  t\right)
\right)  ,%
{\textstyle\sum_{i=1}^{n}}
G_{\bar{x}_{i}}\left(  \mathbb{E}\left[  \bar{g}\left(  \tau,\bar{X}\left(
\tau\right)  \right)  \right]  \right)  \bar{g}_{x}^{i}\left(  \tau,\bar
{X}\left(  \tau\right)  \right)  \right\rangle \\
&  \left.  \left.  +f^{u\left(  \cdot\right)  }\left(  t,\bar{X}\left(
t\right)  \right)  -f^{\bar{u}\left(  \cdot\right)  }\left(  t,\bar{X}\left(
t\right)  \right)  \right\}  dt\right] \\
&  =\mathbb{E}\left[  \int_{\tau}^{\tau+\varepsilon}\left\{  \left\langle
b^{u\left(  \cdot\right)  }\left(  t,X^{\varepsilon}\left(  t\right)  \right)
-b^{\bar{u}\left(  \cdot\right)  }\left(  t,X^{\varepsilon}\left(  t\right)
\right)  ,\bar{\theta}_{x}\left(  t,X^{\varepsilon}\left(  t\right)  \right)
\right\rangle \right.  \right. \\
&  -\left\langle b^{u\left(  \cdot\right)  }\left(  t,\bar{X}\left(  t\right)
\right)  -b^{\bar{u}\left(  \cdot\right)  }\left(  t,\bar{X}\left(  t\right)
\right)  ,\bar{\theta}_{x}\left(  t,\bar{X}\left(  t\right)  \right)
\right\rangle \\
&  +\left\langle b^{u\left(  \cdot\right)  }\left(  t,X^{\varepsilon}\left(
t\right)  \right)  -b^{\bar{u}\left(  \cdot\right)  }\left(  t,X^{\varepsilon
}\left(  t\right)  \right)  ,%
{\textstyle\sum_{i=1}^{n}}
G_{\bar{x}_{i}}\left(  \mathbb{E}\left[  \bar{g}\left(  t,X^{\varepsilon
}\left(  t\right)  \right)  \right]  \right)  \bar{g}_{x}^{i}\left(
t,X^{\varepsilon}\left(  t\right)  \right)  \right\rangle \\
&  -\left\langle b^{u\left(  \cdot\right)  }\left(  t,\bar{X}\left(  t\right)
\right)  -b^{\bar{u}\left(  \cdot\right)  }\left(  t,\bar{X}\left(  t\right)
\right)  ,%
{\textstyle\sum_{i=1}^{n}}
G_{\bar{x}_{i}}\left(  \mathbb{E}\left[  \bar{g}\left(  t,\bar{X}\left(
t\right)  \right)  \right]  \right)  \bar{g}_{x}^{i}\left(  t,\bar{X}\left(
t\right)  \right)  \right\rangle \\
&  \left.  \left.  +f^{u\left(  \cdot\right)  }\left(  t,X^{\varepsilon
}\left(  t\right)  \right)  -f^{\bar{u}\left(  \cdot\right)  }\left(
t,X^{\varepsilon}\left(  t\right)  \right)  -f^{u\left(  \cdot\right)
}\left(  t,\bar{X}\left(  t\right)  \right)  +f^{\bar{u}\left(  \cdot\right)
}\left(  t,\bar{X}\left(  t\right)  \right)  \right\}  dt\right]  .
\end{align*}
Thus%
\begin{align*}
\left\vert \rho\left(  \varepsilon\right)  \right\vert  &  \leq\mathbb{E}%
\left[  \int_{\tau}^{\tau+\varepsilon}\left\{  \left\vert b^{u\left(
\cdot\right)  }\left(  t,X^{\varepsilon}\left(  t\right)  \right)  -b^{\bar
{u}\left(  \cdot\right)  }\left(  t,X^{\varepsilon}\left(  t\right)  \right)
\right\vert \left\vert \bar{\theta}_{x}\left(  t,X^{\varepsilon}\left(
t\right)  \right)  -\bar{\theta}_{x}\left(  t,\bar{X}\left(  t\right)
\right)  \right\vert \right.  \right. \\
&  +\left\vert \bar{\theta}_{x}\left(  t,\bar{X}\left(  t\right)  \right)
\right\vert \left\vert b^{u\left(  \cdot\right)  }\left(  t,X^{\varepsilon
}\left(  t\right)  \right)  -b^{u\left(  \cdot\right)  }\left(  t,\bar
{X}\left(  t\right)  \right)  \right\vert \\
&  +\left\vert \bar{\theta}_{x}\left(  t,\bar{X}\left(  t\right)  \right)
\right\vert \left\vert b^{\bar{u}\left(  \cdot\right)  }\left(
t,X^{\varepsilon}\left(  t\right)  \right)  -b^{\bar{u}\left(  \cdot\right)
}\left(  t,\bar{X}\left(  t\right)  \right)  \right\vert \\
&  +\left\vert f^{u\left(  \cdot\right)  }\left(  t,X^{\varepsilon}\left(
t\right)  \right)  -f^{u\left(  \cdot\right)  }\left(  t,\bar{X}\left(
t\right)  \right)  \right\vert +\left\vert f^{\bar{u}\left(  \cdot\right)
}\left(  t,X^{\varepsilon}\left(  t\right)  \right)  -f^{\bar{u}\left(
\cdot\right)  }\left(  t,\bar{X}\left(  t\right)  \right)  \right\vert \\
&  +\sum_{i=1}^{n}\left\vert G_{\bar{x}_{i}}\left(  \mathbb{E}\left[  \bar
{g}\left(  t,X^{\varepsilon}\left(  t\right)  \right)  \right]  \right)
\right\vert \left\vert g_{x}^{i}\left(  t,X^{\varepsilon}\left(  t\right)
\right)  -g_{x}^{i}\left(  t,\bar{X}\left(  t\right)  \right)  \right\vert
\left\vert \sigma^{u\left(  \cdot\right)  }\left(  t,X^{\varepsilon}\left(
t\right)  \right)  -\sigma^{\bar{u}\left(  \cdot\right)  }\left(
t,X^{\varepsilon}\left(  t\right)  \right)  \right\vert \\
&  +\sum_{i=1}^{n}\left\vert \bar{g}_{x}^{i}\left(  t,\bar{X}\left(  t\right)
\right)  \right\vert \left\vert G_{\bar{x}_{i}}\left(  \mathbb{E}\left[
\bar{g}\left(  t,X^{\varepsilon}\left(  t\right)  \right)  \right]  \right)
-G_{\bar{x}_{i}}\left(  \mathbb{E}\left[  \bar{g}\left(  t,\bar{X}\left(
t\right)  \right)  \right]  \right)  \right\vert \left\vert \sigma^{u\left(
\cdot\right)  }\left(  t,X^{\varepsilon}\left(  t\right)  \right)
-\sigma^{\bar{u}\left(  \cdot\right)  }\left(  t,X^{\varepsilon}\left(
t\right)  \right)  \right\vert \\
&  +\sum_{i=1}^{n}\left\vert G_{\bar{x}_{i}}\left(  \mathbb{E}\left[  \bar
{g}\left(  t,\bar{X}\left(  t\right)  \right)  \right]  \right)  \right\vert
\left\vert \bar{g}_{x}^{i}\left(  t,\bar{X}\left(  t\right)  \right)
\right\vert \left\vert \sigma^{u\left(  \cdot\right)  }\left(
t,X^{\varepsilon}\left(  t\right)  \right)  -\sigma^{u\left(  \cdot\right)
}\left(  t,\bar{X}\left(  t\right)  \right)  \right\vert \\
&  \left.  \left.  +\sum_{i=1}^{n}\left\vert G_{\bar{x}_{i}}\left(
\mathbb{E}\left[  \bar{g}\left(  t,\bar{X}\left(  t\right)  \right)  \right]
\right)  \right\vert \left\vert \bar{g}_{x}^{i}\left(  t,\bar{X}\left(
t\right)  \right)  \right\vert \left\vert \sigma^{u\left(  \cdot\right)
}\left(  t,X^{\varepsilon}\left(  t\right)  \right)  -\sigma^{u\left(
\cdot\right)  }\left(  t,\bar{X}\left(  t\right)  \right)  \right\vert
\right\}  dt\right]  \text{.}%
\end{align*}

Using the fact that the random fields $b_{x}^{u\left(  \cdot\right)  }$,
$b_{x}^{\bar{u}\left(  \cdot\right)  }$, $f_{x}^{u\left(  \cdot\right)  }$,
$f_{x}^{\bar{u}\left(  \cdot\right)  }$, $G_{\bar{x}}$, $G_{\bar{x}\bar{x}}$,
$\bar{\theta}_{x}$, $\bar{\theta}_{xx}$, $\bar{g}_{x}^{i}$, $\bar{g}_{xx}^{i}$
are uniformly bounded, it is not difficult to show that there exists a
constant $K>0$ such that%
\begin{align*}
\left\vert \rho\left(  \varepsilon\right)  \right\vert  &  \leq K\int_{\tau
}^{\tau+\varepsilon}\mathbb{E}\left[  \left\vert \bar{X}\left(  t\right)
-X^{\varepsilon}\left(  t\right)  \right\vert \right]  dt\\
&  =K\mathbb{E}\left[  \int_{\tau}^{\tau+\varepsilon}\mathbb{E}\left[
\left\vert \xi^{\varepsilon}\left(  t\right)  \right\vert \right]  dt\right]
\\
&  \leq K\varepsilon\sup_{t\in\left[  0,T\right]  }\mathbb{E}\left[
\left\vert \xi^{\varepsilon}\left(  t\right)  \right\vert ^{2k}\right]
^{\frac{1}{2k}}\\
&  \leq K\varepsilon^{2}.
\end{align*}

Hence%
\[
\left\vert \rho\left(  \varepsilon\right)  \right\vert \leq o\left(
\varepsilon\right)  \text{.}%
\]
This completes the proof.\eop

Now, we are ready to give a proof of Theorem \ref{result04}.

\textbf{Proof of Theorem \ref{result04}. }From Lemma \ref{result07}, we deduce
that%
\begin{align*}
0  &  \leq\lim_{\varepsilon\downarrow0}\frac{1}{\varepsilon}\left\{
\mathbf{J}\left(  u^{\varepsilon}\left(  \cdot\right)  \right)  -\mathbf{J}%
\left(  \bar{u}\left(  \cdot\right)  \right)  \right\} \\
&  =\lim_{\varepsilon\downarrow0}\frac{1}{\varepsilon}\int_{\tau}%
^{\tau+\varepsilon}\mathbb{E}\left[  \left\langle b\left(  t,\bar{X}\left(
t\right)  ,u\left(  t\right)  \right)  -b\left(  t,\bar{X}\left(  t\right)
,\bar{u}\left(  t\right)  \right)  ,\bar{\theta}_{x}\left(  t,\bar{X}\left(
t\right)  \right)  \right\rangle \right. \\
&  +%
{\textstyle\sum_{i=1}^{n}}
\left\langle b\left(  t,\bar{X}\left(  t\right)  ,u\left(  t\right)  \right)
-b\left(  t,\bar{X}\left(  t\right)  ,\bar{u}\left(  t\right)  \right)
,G_{\bar{x}_{i}}\left(  \mathbb{E}\left[  \bar{g}\left(  t,\bar{X}\left(
t\right)  \right)  \right]  \right)  \bar{g}_{x}^{i}\left(  t,\bar{X}\left(
t\right)  \right)  \right\rangle \\
&  +\left.  \left(  f\left(  t,\bar{X}\left(  t\right)  ,u\left(  t\right)
\right)  -f\left(  t,\bar{X}\left(  t\right)  ,\bar{u}\left(  t\right)
\right)  \right)  \right]  dt\\
&  =\lim_{\varepsilon\downarrow0}\frac{1}{\varepsilon}\int_{\tau}%
^{\tau+\varepsilon}\mathbb{E}\left[  \left\langle b\left(  t,\bar{X}\left(
t\right)  ,u\left(  t\right)  \right)  -b\left(  t,\bar{X}\left(  t\right)
,\bar{u}\left(  t\right)  \right)  ,\bar{\theta}_{x}\left(  t,\bar{X}\left(
t\right)  \right)  \right\rangle \right. \\
&  +%
{\textstyle\sum_{i=1}^{n}}
\left\langle b\left(  t,\bar{X}\left(  t\right)  ,u\left(  t\right)  \right)
-b\left(  t,\bar{X}\left(  t\right)  ,\bar{u}\left(  t\right)  \right)
,G_{\bar{x}_{i}}\left(  \mathbb{E}\left[  \bar{X}\left(  T\right)  \right]
\right)  \bar{g}_{x}^{i}\left(  t,\bar{X}\left(  t\right)  \right)
\right\rangle \\
&  \left.  +f\left(  t,\bar{X}\left(  t\right)  ,u\left(  t\right)  \right)
-f\left(  t,\bar{X}\left(  t\right)  ,\bar{u}\left(  t\right)  \right)
\right]  dt,
\end{align*}
where we have used the following equality,%
\begin{align*}
\mathbb{E}\left[  \bar{g}\left(  t,\bar{X}\left(  t\right)  \right)  \right]
&  =\left(  \mathbb{E}\left[  \bar{g}^{1}\left(  t,\bar{X}\left(  t\right)
\right)  \right]  ,...,\mathbb{E}\left[  \bar{g}^{n}\left(  t,\bar{X}\left(
t\right)  \right)  \right]  \right)  ^{\top}\\
&  =\left(  \mathbb{E}\left[  \mathbb{E}_{t}\left[  \bar{X}_{1}\left(
T\right)  \right]  \right]  ,...,\mathbb{E}\left[  \mathbb{E}_{t}\left[
\bar{X}_{n}\left(  T\right)  \right]  \right]  \right)  ^{\top}\\
&  =\mathbb{E}\left[  \bar{X}\left(  T\right)  \right]  .
\end{align*}

Thus by Lebesgue differentiation theorem%
\begin{align}
0  &  \leq\mathbb{E}\left[  \left\langle b\left(  \tau,\bar{X}\left(
\tau\right)  ,u\left(  \tau\right)  \right)  -b\left(  \tau,\bar{X}\left(
\tau\right)  ,\bar{u}\left(  \tau\right)  \right)  ,\theta_{x}\left(
\tau,\hat{X}\left(  \tau\right)  \right)  \right\rangle \right. \nonumber\\
&  +%
{\textstyle\sum_{i=1}^{n}}
\left\langle b\left(  \tau,\bar{X}\left(  \tau\right)  ,u\left(  \tau\right)
\right)  -b\left(  \tau,\bar{X}\left(  \tau\right)  ,\bar{u}\left(
\tau\right)  \right)  ,G_{\bar{x}_{i}}\left(  \mathbb{E}\left[  \bar{X}\left(
T\right)  \right]  \right)  \bar{g}_{x}^{i}\left(  \tau,\bar{X}\left(
\tau\right)  \right)  \right\rangle \nonumber\\
&  \left.  +f\left(  \tau,\bar{X}\left(  \tau\right)  ,u\left(  \tau\right)
\right)  -f\left(  \tau,\bar{X}\left(  \tau\right)  ,\bar{u}\left(
\tau\right)  \right)  \right]  \text{, a.e. }\tau\in\left[  0,T\right]
\text{.} \label{Eq28}%
\end{align}

Now let%
\[
u\left(  \tau\right)  =u\chi_{A}+\bar{u}\left(  \tau\right)  \chi
_{\Omega\backslash A},
\]
where $A$ is an arbitrarily $\mathcal{F}_{\tau}$-measurable set and $u\in U.$

Hence, in view of (\ref{Eq28}), we have%
\begin{align*}
0  &  \leq\mathbb{E}\left[  \left\{  \left\langle b\left(  \tau,\bar{X}\left(
\tau\right)  ,u\right)  -b\left(  \tau,\bar{X}\left(  \tau\right)  ,\bar
{u}\left(  \tau\right)  \right)  ,\theta_{x}\left(  \tau,\hat{X}\left(
\tau\right)  \right)  \right\rangle \right.  \right. \\
&  +%
{\textstyle\sum_{i=1}^{n}}
\left\langle b\left(  \tau,\bar{X}\left(  \tau\right)  ,u\right)  -b\left(
\tau,\bar{X}\left(  \tau\right)  ,\bar{u}\left(  \tau\right)  \right)
,G_{\bar{x}_{i}}\left(  \mathbb{E}\left[  \bar{X}\left(  T\right)  \right]
\right)  \bar{g}_{x}^{i}\left(  \tau,\bar{X}\left(  \tau\right)  \right)
\right\rangle \\
&  \left.  \left.  +f\left(  \tau,\bar{X}\left(  \tau\right)  ,u\right)
-f\left(  \tau,\bar{X}\left(  \tau\right)  ,\bar{u}\left(  \tau\right)
\right)  \right\}  \chi_{A}\right]  .
\end{align*}
which in turn yields inequality (\ref{Eq14}) since $u\in U$ and the set
$A\in\mathcal{F}_{\tau}$ are arbitrary. This completes the proof.\eop

\section{ Sufficient Condition of Optimality\label{section6}}

In this section, we focus on proving that provided some convexity assumptions
are satisfied, the necessary condition (\ref{Eq15}) turns out to be sufficient.

The following fact concerning the differentiability of stochastic integrals
with parameter is important for our purpose. Let $\phi\left(  \cdot
,\cdot\right)  $ be a bounded random field such that $\phi\left(  \cdot
,\cdot\right)  \in\mathcal{L}_{\mathcal{F}}^{2}\left(  0,T;C_{b}^{2}\left(
\mathbb{R};\mathbb{R}\right)  \right)  $. Then it was shown that (see, e.g.
[Kunita \cite{Kunita2}, Proposition 2.3.1, Page 56] or [Kunita \cite{Kunita1},
Exercise 3.1.5, Page 78]) the stochastic integral with parameter: $\int
_{0}^{\cdot}\phi\left(  t,\cdot\right)  dW\left(  t\right)  $ has a
modification that belongs to $\mathcal{L}_{\mathcal{F}}^{2}\left(
0,T;C^{1}\left(  \mathbb{R}^{n};\mathbb{R}\right)  \right)  $ and it satisfies%
\begin{equation}
\partial_{x}\int_{0}^{t}\phi\left(  t,x\right)  dW\left(  t\right)  =\int
_{0}^{t}\partial_{x}\phi\left(  t,x\right)  dW\left(  t\right)  ,
\label{Eq300}%
\end{equation}
where $\partial_{x}\phi\left(  t,x\right)  $ denote the first derivative of
$\phi\left(  t,x\right)  $ with respect to the variable $x.$

Define the usual Hamiltonian as a map from $\left[  0,T\right]  \times%
\mathbb{R}
^{n}\times U\times%
\mathbb{R}
^{n}\times%
\mathbb{R}
^{n\times d}$ into $%
\mathbb{R}
$ by%
\begin{equation}
\mathbb{H}\left(  t,x,u,p,q\right)  :=\left\langle b\left(  t,x,u\right)
,p\right\rangle +\text{\textbf{tr}}\left[  \sigma\left(  t,x\right)  ^{\top
}q\right]  +f\left(  t,x,u\right)  \label{Eq61}%
\end{equation}
and let us introduce two additional assumptions.

\begin{enumerate}
\item[\textbf{(H1)}] The control domain $U\subset%
\mathbb{R}
^{l}$ is a convex subset and the maps $b$ and $f$ are continuously
differentiable with respect to $u$.

\item[\textbf{(H2)}] The Hamiltonian $\mathbb{H}\left(  t,x,u,p,q\right)  $ is
convex with respect to $\left(  x,u\right)  $, $h\left(  \cdot\right)  $ is
convex with respect to $x$ and $G\left(  \cdot\right)  $ is convex with
respect to $\bar{x}$.
\end{enumerate}

\begin{theorem}
[Sufficient Condition of Optimality]\label{result08}Let Assumption
\textbf{(H1), (H2) }and\textbf{ (H}$_{m}$\textbf{)} be satisfied with
$m>2+\frac{n}{2}$. Let $\left(  \bar{u}\left(  \cdot\right)  ,\bar{X}\left(
\cdot\right)  \right)  $ be an admissible pair and $\left(  \bar{\theta
}\left(  \cdot,\cdot\right)  ,\bar{\psi}\left(  \cdot,\cdot\right)  \right)
,$ $\left(  \bar{g}^{i}\left(  \cdot,\cdot\right)  ,\bar{\eta}^{i}\left(
\cdot,\cdot\right)  \right)  $, for $1\leq i\leq n,$ be the solutions of the
BSPDEs (\ref{Eq13}) and (\ref{Eq13**}), respectively. Suppose that%
\[
\bar{\psi}\left(  \cdot,\cdot\right)  \in\mathcal{L}_{\mathcal{F}}^{2}\left(
0,T;C_{b}^{2}\left(  \mathbb{R}^{n};\mathbb{R}^{d}\right)  \right)  \text{ and
}\bar{\eta}^{i}\left(  \cdot,\cdot\right)  \in\mathcal{L}_{\mathcal{F}}%
^{2}\left(  0,T;C_{b}^{2}\left(  \mathbb{R}^{n};\mathbb{R}^{d}\right)
\right)  .
\]
Suppose further, for a.e. $t\in\left[  0,T\right]  $,%
\begin{equation}%
\begin{array}
[c]{l}%
\bar{u}\left(  t\right)  \in\arg\min\limits_{u\in U}\left\{  \left\langle
\bar{\theta}_{x}\left(  t,\bar{X}\left(  t\right)  \right)  ,b\left(
t,\bar{X}\left(  t\right)  ,\cdot\right)  \right\rangle +f\left(  t,\bar
{X}\left(  t\right)  ,\cdot\right)  \right. \\
\text{ \ \ \ \ \ \ \ \ }%
\begin{array}
[c]{r}%
\left.  +%
{\textstyle\sum_{i=1}^{n}}
G_{\bar{x}_{i}}\left(  \mathbb{E}\left[  \bar{X}\left(  T\right)  \right]
\right)  \left\langle \bar{g}_{x}^{i}\left(  \tau,\bar{X}\left(  \tau\right)
\right)  ,b\left(  t,\bar{X}\left(  t\right)  ,\cdot\right)  \right\rangle
\right\}  ,\text{ }\\
\text{a.e. t}\in\left[  0,T\right]  \text{, a.s,}%
\end{array}
\end{array}
\label{Eq60}%
\end{equation}
Then $\left(  \bar{u}\left(  \cdot\right)  ,\bar{X}\left(  \cdot\right)
\right)  $ is an optimal pair of Problem (S). Furthermore the objective value
of $\bar{u}\left(  \cdot\right)  $ is given by%
\[
\mathbf{J}\left(  \bar{u}\left(  \cdot\right)  \right)  =\inf_{u\left(
\cdot\right)  \in\mathcal{U}\left[  0,T\right]  }\mathbf{J}\left(  u\left(
\cdot\right)  \right)  =\bar{\theta}\left(  0,x_{0}\right)  +G\left(  \bar
{g}\left(  0,x_{0}\right)  \right)  .
\]

\end{theorem}

\bop First, let us begin by some preparations. Let $\left(  \bar{u}\left(
\cdot\right)  ,\bar{X}\left(  \cdot\right)  \right)  $ be an admissible pair
and $\left(  \bar{\theta}\left(  \cdot,\cdot\right)  ,\bar{\psi}\left(
\cdot,\cdot\right)  \right)  $, $\left(  \bar{g}_{i}\left(  \cdot
,\cdot\right)  ,\bar{\eta}_{i}\left(  \cdot,\cdot\right)  \right)  $, for
$1\leq i\leq n,$ be the solutions of the following BSPDEs,%
\[
\left\{
\begin{array}
[c]{l}%
d\bar{\theta}\left(  t,x\right)  =-\left\{  \left\langle \bar{\theta}%
_{x}\left(  t,x\right)  ,b^{\bar{u}\left(  \cdot\right)  }\left(  t,x\right)
\right\rangle +\frac{1}{2}\text{\textbf{tr}}\left[  \left(  \sigma\sigma
^{\top}\right)  \left(  t,x\right)  \bar{\theta}_{xx}\left(  t,x\right)
\right]  \right. \\
\text{ \ \ \ \ \ \ \ \ \ \ \ \ \ \ \ }+\left.  \text{\textbf{tr}}\left[
\bar{\psi}_{x}\left(  t,x\right)  \sigma\left(  t,x\right)  \right]
+f^{\bar{u}\left(  \cdot\right)  }\left(  t,x\right)  \right\}  dt\\
\text{ \ \ \ \ \ \ \ \ \ \ \ \ \ \ \ }+\bar{\psi}\left(  t,x\right)  ^{\top
}dW\left(  t\right)  \text{, }\left(  t,x\right)  \in\left[  0,T\right]
\times\mathbb{R}^{n}\text{,}\\
\bar{\theta}\left(  T,x\right)  =h\left(  x\right)  \text{, for }%
x\in\mathbb{R}^{n}.
\end{array}
\right.
\]
and%
\[
\left\{
\begin{array}
[c]{l}%
d\bar{g}^{i}\left(  t,x\right)  =-\left\{  \left\langle \bar{g}_{x}^{i}\left(
t,x\right)  ,b^{\bar{u}\left(  \cdot\right)  }\left(  t,x\right)
\right\rangle +\frac{1}{2}\text{\textbf{tr}}\left[  \left(  \sigma\sigma
^{\top}\right)  \left(  t,x\right)  \bar{g}_{xx}^{i}\left(  t,x\right)
\right]  \right. \\
\text{ \ \ \ \ \ \ \ \ \ \ \ \ }%
\begin{array}
[c]{r}%
+\left.  \text{\textbf{tr}}\left[  \bar{\eta}_{x}^{i}\left(  t,x\right)
^{\top}\sigma\left(  t,x\right)  \right]  \right\}  dt+\bar{\eta}^{i}\left(
t,x\right)  ^{\top}dW\left(  t\right)  \text{,}\\
\text{for }\left(  t,x\right)  \in\left[  0,T\right]  \times\mathbb{R}%
^{n}\text{,}%
\end{array}
\\
\bar{g}^{i}\left(  T,x\right)  =\bar{x}_{i}\text{,\ for }x\in\mathbb{R}%
^{n}\text{,}%
\end{array}
\right.
\]
respectively, which also can be written as%
\begin{equation}
\left\{
\begin{array}
[c]{l}%
d\bar{\theta}\left(  t,x\right)  =-\left\{  \left\langle b^{\bar{u}\left(
\cdot\right)  }\left(  t,x\right)  ,\bar{\theta}_{x}\left(  t,x\right)
\right\rangle +\frac{1}{2}\sum\limits_{j=1}^{d}\sigma^{j}\left(  t,x\right)
^{\top}\bar{\theta}_{xx}\left(  t,x\right)  \sigma^{j}\left(  t,x\right)
\right. \\
\text{ \ \ \ \ \ \ \ \ \ \ \ \ \ \ \ }+\left.  \sum\limits_{j=1}%
^{d}\left\langle \sigma^{j}\left(  t,x\right)  ,\bar{\psi}_{x}^{j}\left(
t,x\right)  \right\rangle +f^{\bar{u}\left(  \cdot\right)  }\left(
t,x\right)  \right\}  dt\\
\text{ \ \ \ \ \ \ \ \ \ \ \ \ \ \ \ }+\sum\limits_{j=1}^{d}\bar{\psi}%
^{j}\left(  t,x\right)  dW_{j}\left(  t\right)  \text{, }\left(  t,x\right)
\in\left[  0,T\right]  \times\mathbb{R}^{n}\text{,}\\
\bar{\theta}\left(  T,x\right)  =h\left(  x\right)  \text{, for }%
x\in\mathbb{R}^{n}%
\end{array}
\right.  \label{Eq29}%
\end{equation}
and%
\begin{equation}
\left\{
\begin{array}
[c]{l}%
d\bar{g}^{i}\left(  t,x\right)  =-\left\{  \left\langle \bar{g}_{x}^{i}\left(
t,x\right)  ,b^{\bar{u}\left(  \cdot\right)  }\left(  t,x\right)
\right\rangle +\frac{1}{2}\sum\limits_{j=1}^{d}\sigma^{j}\left(  t,x\right)
^{\top}\bar{g}_{xx}^{i}\left(  t,x\right)  \sigma^{j}\left(  t,x\right)
\right. \\
\text{ \ \ \ \ \ \ \ \ \ \ \ \ }%
\begin{array}
[c]{r}%
+\left.  \sum\limits_{j=1}^{d}\left\langle \sigma^{j}\left(  t,x\right)
,\bar{\eta}_{x}^{ij}\left(  t,x\right)  \right\rangle \right\}  dt+\sum
\limits_{j=1}^{d}\bar{\eta}^{ij}\left(  t,x\right)  dW_{j}\left(  t\right)
\text{,}\\
\text{for }\left(  t,x\right)  \in\left[  0,T\right]  \times\mathbb{R}%
^{n}\text{,}%
\end{array}
\\
\bar{g}^{i}\left(  T,x\right)  =\bar{x}_{i}\text{, for }x\in\mathbb{R}%
^{n}\text{.}%
\end{array}
\right.  \label{Eq29*}%
\end{equation}

By using the fact (\ref{Eq300}), we differentiate (\ref{Eq29})-(\ref{Eq29*})
in $x$ to obtain that%
\[
\left\{
\begin{array}
[c]{l}%
d\bar{\theta}_{x}\left(  t,x\right)  =-\left\{  \left\langle b_{x}^{\bar
{u}\left(  \cdot\right)  }\left(  t,x\right)  ,\bar{\theta}_{x}\left(
t,x\right)  \right\rangle +\bar{\theta}_{xx}\left(  t,x\right)  ^{\top}%
b_{x}^{\bar{u}\left(  \cdot\right)  }\left(  t,x\right)  \right. \\
\text{ \ \ \ \ \ \ \ \ \ \ \ \ \ \ \ }+\sum\limits_{j=1}^{d}\sigma_{x}%
^{j}\left(  t,x\right)  ^{\top}\bar{\theta}_{xx}\left(  t,x\right)  \sigma
^{j}\left(  t,x\right)  +\frac{1}{2}\sum\limits_{j=1}^{d}\sigma_{x}^{j}\left(
t,x\right)  ^{\top}\bar{\theta}_{xxx}\left(  t,x\right)  \sigma^{j}\left(
t,x\right) \\
\text{ \ \ \ \ \ \ \ \ \ \ \ \ \ \ \ }+\sum\limits_{j=1}^{d}\left[
\left\langle \sigma_{x}^{j}\left(  t,x\right)  ,\bar{\psi}_{x}^{j}\left(
t,x\right)  \right\rangle +\bar{\psi}_{xx}^{j}\left(  t,x\right)  ^{\top
}\sigma_{x}^{j}\left(  t,x\right)  \right]  +\left.  f_{x}^{\bar{u}\left(
\cdot\right)  }\left(  t,x\right)  \right\}  dt\text{\ \ \ \ \ \ \ \ \ \ }\\
\text{ \ \ \ \ \ \ \ \ \ \ \ \ \ \ \ }+\sum\limits_{j=1}^{d}\bar{\psi}_{x}%
^{j}\left(  t,x\right)  dW_{j}\left(  t\right)  \text{, }\left(  t,x\right)
\in\left[  0,T\right]  \times\mathbb{R}^{n}\text{,}\\
\bar{\theta}_{x}\left(  T,x\right)  =h_{x}\left(  x\right)  \text{, for }%
x\in\mathbb{R}^{n}%
\end{array}
\right.
\]
and%
\[
\left\{
\begin{array}
[c]{l}%
d\bar{g}_{x}^{i}\left(  t,x\right)  =-\left\{  \left\langle b_{x}^{\bar
{u}\left(  \cdot\right)  }\left(  t,x\right)  ,\bar{g}_{x}^{i}\left(
t,x\right)  \right\rangle +\bar{g}_{xx}^{i}\left(  t,x\right)  ^{\top}%
b_{x}^{\bar{u}\left(  \cdot\right)  }\left(  t,x\right)  \right. \\
\text{ \ \ \ \ \ \ \ \ \ \ \ \ }+\sum\limits_{j=1}^{d}\sigma_{x}^{j}\left(
t,x\right)  ^{\top}\bar{g}_{xx}^{i}\left(  t,x\right)  \sigma^{j}\left(
t,x\right)  +\frac{1}{2}\sum\limits_{j=1}^{d}\sigma_{x}^{j}\left(  t,x\right)
^{\top}\bar{g}_{xxx}^{i}\left(  t,x\right)  \sigma^{j}\left(  t,x\right) \\
\text{ \ \ \ \ \ \ \ \ \ \ \ }+\left.  \sum\limits_{j=1}^{d}\left[
\left\langle \sigma_{x}^{j,\bar{u}\left(  \cdot\right)  }\left(  t,x\right)
,\bar{\eta}_{x}^{ij}\left(  t,x\right)  \right\rangle +\bar{\eta}_{xx}%
^{ij}\left(  t,x\right)  ^{\top}\sigma_{x}^{j}\left(  t,x\right)  \right]
\right\}  dt\text{\ \ \ \ \ \ \ \ \ \ }\\
\text{ \ \ \ \ \ \ \ \ \ \ \ \ }+\sum\limits_{j=1}^{d}\bar{\eta}_{x}%
^{ij}\left(  t,x\right)  dW_{j}\left(  t\right)  \text{, }\left(  t,x\right)
\in\left[  0,T\right]  \times\mathbb{R}^{n}\text{,}\\
\bar{g}_{x}^{i}\left(  T,x\right)  =\mathbf{e}_{i}\text{, for }x\in
\mathbb{R}^{n},
\end{array}
\right.
\]
where for $\varrho=\bar{\theta},$ $\bar{g}^{i}$%
\[
\left(  \sigma_{x}^{j}\right)  ^{\top}\varrho_{xxx}\left(  \sigma^{j}\right)
:=\left(  \left(  \sigma_{x}^{j}\right)  ^{\top}\left(  \varrho_{x}\right)
_{xx}^{1}\left(  \sigma^{j}\right)  ,...,\left(  \sigma_{x}^{j}\right)
^{\top}\left(  \varrho_{x}\right)  _{xx}^{n}\sigma^{j}\right)  ,
\]
with%
\[
\varrho_{x}:=\left(  \left(  \varrho_{x}\right)  ^{1},...,\left(  \varrho
_{x}\right)  ^{n}\right)  ^{\top}.
\]

Accordingly, applying It\^{o}-Wentzell formula to $\bar{\theta}_{x}\left(
t,\bar{X}\left(  t\right)  \right)  $ and $\bar{g}_{x}^{i}\left(  t,\bar
{X}\left(  t\right)  \right)  $, we get%
\[
\left\{
\begin{array}
[c]{l}%
d\bar{\theta}_{x}\left(  t,\bar{X}\left(  t\right)  \right)  =-\left\{
\left\langle b_{x}^{\bar{u}\left(  \cdot\right)  }\left(  t,\bar{X}\left(
t\right)  \right)  ,\bar{\theta}_{x}\left(  t,\bar{X}\left(  t\right)
\right)  \right\rangle +f_{x}\left(  t,\bar{X}\left(  t\right)  ,\bar
{u}\left(  t\right)  \right)  \right. \\
\text{ \ \ \ \ \ \ \ \ \ \ \ \ \ \ \ }+\left.  \sum\limits_{j=1}%
^{d}\left\langle \sigma_{x}^{j}\left(  t,\bar{X}\left(  t\right)  \right)
,\bar{\psi}_{x}^{j}\left(  t,\bar{X}\left(  t\right)  \right)  +\bar{\theta
}_{xx}\left(  t,\bar{X}\left(  t\right)  \right)  \sigma^{j}\left(  t,\bar
{X}\left(  t\right)  \right)  \right\rangle \right\}  dt\\
\text{ \ \ \ \ \ \ \ \ \ \ \ \ }+\sum\limits_{j=1}^{d}\left\{  \bar{\psi}%
_{x}^{j}\left(  t,\bar{X}\left(  t\right)  \right)  +\bar{\theta}_{xx}\left(
t,\bar{X}\left(  t\right)  \right)  \sigma^{j}\left(  t,\bar{X}\left(
t\right)  \right)  \right\}  dW_{j}\left(  t\right)  \text{, }t\in\left[
0,T\right]  \text{,}\\
\bar{\theta}_{x}\left(  T,\bar{X}\left(  T\right)  \right)  =h_{x}\left(
\bar{X}\left(  T\right)  \right)  \text{,}%
\end{array}
\right.
\]
and, for each $1\leq i\leq n,$%
\[
\left\{
\begin{array}
[c]{l}%
d\bar{g}_{x}^{i}\left(  t,\bar{X}\left(  t\right)  \right)  =-\left\{
\left\langle b_{x}^{\bar{u}\left(  \cdot\right)  }\left(  t,\bar{X}\left(
t\right)  \right)  ,\bar{g}_{x}^{i}\left(  t,\bar{X}\left(  t\right)  \right)
\right\rangle \right. \\
\text{ \ \ \ \ \ \ \ \ \ \ \ \ \ \ \ }+\left.  \sum\limits_{j=1}%
^{d}\left\langle \sigma_{x}^{j}\left(  t,\bar{X}\left(  t\right)  \right)
,\bar{\eta}_{x}^{ij}\left(  t,\bar{X}\left(  t\right)  \right)  +\bar{g}%
_{xx}^{i}\left(  t,\bar{X}\left(  t\right)  \right)  \sigma^{j}\left(
t,\bar{X}\left(  t\right)  \right)  \right\rangle \right\}  dt\\
\text{ \ \ \ \ \ \ \ \ \ \ \ \ }+\sum\limits_{j=1}^{d}\left\{  \bar{\eta}%
_{x}^{ij}\left(  t,\bar{X}\left(  t\right)  \right)  +\bar{g}_{xx}^{i}\left(
t,\bar{X}\left(  t\right)  \right)  \sigma^{j}\left(  t,\bar{X}\left(
t\right)  \right)  \right\}  dW_{j}\left(  t\right)  \text{, }t\in\left[
0,T\right]  \text{,}\\
\bar{g}_{x}^{i}\left(  T,\bar{X}\left(  T\right)  \right)  =\mathbf{e}%
_{i}\text{.}%
\end{array}
\right.
\]

Define for $t\in\left[  0,T\right]  ,$%
\begin{equation}
\bar{p}\left(  t\right)  :=\bar{\theta}_{x}\left(  t,\bar{X}\left(  t\right)
\right)  +\sum_{i=1}^{n}G_{\bar{x}_{i}}\left(  \mathbb{E}\left[  \bar
{X}\left(  T\right)  \right]  \right)  \bar{g}_{x}^{i}\left(  t,\bar{X}\left(
t\right)  \right)  \label{Eq2000}%
\end{equation}
and for each $j=1,...,d,$%
\begin{align}
\bar{q}_{j}\left(  t\right)   &  :=\bar{\psi}_{x}^{j}\left(  t,\bar{X}\left(
t\right)  \right)  +\bar{\theta}_{xx}\left(  t,\bar{X}\left(  t\right)
\right)  \sigma^{j}\left(  t,\bar{X}\left(  t\right)  \right) \nonumber\\
&  +\sum_{i=1}^{n}\bar{g}_{\bar{x}_{i}}\left(  \mathbb{E}\left[  \bar
{X}\left(  T\right)  \right]  \right)  \left\{  \bar{\eta}_{x}^{ij}\left(
t,\bar{X}\left(  t\right)  \right)  +\bar{g}_{xx}^{i}\left(  t,\bar{X}\left(
t\right)  \right)  \sigma^{j}\left(  t,\bar{X}\left(  t\right)  \right)
\right\}  . \label{Eq2001}%
\end{align}

Then it is easy to verify that $\bar{p}\left(  \cdot\right)  $ and $\bar
{q}\left(  \cdot\right)  =\left(  \bar{q}_{1}\left(  \cdot\right)
,...,\bar{q}_{d}\left(  \cdot\right)  \right)  $ satisfy the following linear
BSDEs%
\begin{equation}
\left\{
\begin{array}
[c]{l}%
d\bar{p}\left(  t\right)  =-\mathbb{H}_{x}\left(  t,\bar{X}\left(  t\right)
,\bar{u}\left(  t\right)  ,\bar{p}\left(  t\right)  ,\bar{q}\left(  t\right)
\right)  dt+\sum\limits_{j=1}^{d}\bar{q}^{j}\left(  t\right)  dW_{j}\left(
t\right)  \text{, }t\in\left[  0,T\right]  \text{,}\\
\bar{p}\left(  T\right)  =h_{x}\left(  \bar{X}\left(  T\right)  \right)
+\bar{g}_{\bar{x}}\left(  \mathbb{E}\left[  \bar{X}\left(  T\right)  \right]
\right)  \text{.}%
\end{array}
\right.  \label{Eq30}%
\end{equation}
where $\mathbb{H}$ is given by (\ref{Eq61}).

\bigskip Recall that%
\[%
\begin{array}
[c]{l}%
\bar{u}\left(  t\right)  \in\arg\min\limits_{u\in U}\left\{  \left\langle
\bar{\theta}_{x}\left(  t,\bar{X}\left(  t\right)  \right)  ,b\left(
t,\bar{X}\left(  t\right)  ,\cdot\right)  \right\rangle +f\left(  t,\bar
{X}\left(  t\right)  ,\cdot\right)  \right. \\
\text{ \ \ \ \ \ \ \ \ }%
\begin{array}
[c]{r}%
\left.  +%
{\textstyle\sum_{i=1}^{n}}
G_{\bar{x}_{i}}\left(  \mathbb{E}\left[  \bar{X}\left(  T\right)  \right]
\right)  \left\langle \bar{g}_{x}^{i}\left(  \tau,\bar{X}\left(  \tau\right)
\right)  ,b\left(  t,\bar{X}\left(  t\right)  ,\cdot\right)  \right\rangle
\right\}  ,\text{ }\\
\text{a.e. t}\in\left[  0,T\right]  \text{, a.s.}%
\end{array}
\end{array}
\]
Accordingly, the first order optimality condition yields%
\begin{equation}%
\begin{array}
[c]{l}%
0\leq\left\langle u-\bar{u}\left(  t\right)  ,b_{u}\left(  t,\bar{X}\left(
t\right)  ,\bar{u}\left(  t\right)  \right)  ^{\top}\bar{\theta}_{x}\left(
t,\bar{X}\left(  t\right)  \right)  +f_{u}\left(  t,\bar{X}\left(  t\right)
,\bar{u}\left(  t\right)  \right)  \right\rangle \\
\text{ \ \ }%
\begin{array}
[c]{r}%
+%
{\textstyle\sum_{i=1}^{n}}
G_{\bar{x}_{i}}\left(  \mathbb{E}\left[  \bar{X}\left(  T\right)  \right]
\right)  \left\langle u-\bar{u}\left(  t\right)  ,b_{u}\left(  t,\bar
{X}\left(  t\right)  ,\bar{u}\left(  t\right)  \right)  ^{\top}\bar{g}_{x}%
^{i}\left(  t,\bar{X}\left(  t\right)  \right)  \right\rangle \\
\text{for all }u\in U,\text{ a.s., a.e. }t\in\left[  0,T\right]  .
\end{array}
\end{array}
\label{Eq36}%
\end{equation}
Combining this with Equalities (\ref{Eq2000})-(\ref{Eq2001}), we get%
\begin{equation}
0\leq\left(  u-\bar{u}\left(  t\right)  \right)  ^{\top}\mathbb{H}_{u}\left(
t,\hat{X}^{x_{0}}\left(  t\right)  ,\bar{u}\left(  t\right)  ,p\left(
t\right)  ,q\left(  t\right)  ,r\left(  t,\cdot\right)  \right)  ,\text{
}\mathbb{P-}\text{a.s., }\forall u\in U\text{, a.e. }t\in\left[  0,T\right]
\text{.} \label{4.9}%
\end{equation}

Now for an arbitrarily admissible state-control pair $\left(  u\left(
\cdot\right)  ,X\left(  \cdot\right)  \right)  $, consider the difference%
\begin{align*}
&  \mathbf{J}\left(  u\left(  \cdot\right)  \right)  -\mathbf{J}\left(
\bar{u}\left(  \cdot\right)  \right) \\
&  =\mathbb{E}\left[  \int_{0}^{T}\left\{  f^{u\left(  \cdot\right)  }\left(
t,X\left(  t\right)  \right)  -f^{\bar{u}\left(  \cdot\right)  }\left(
t,\bar{X}\left(  t\right)  \right)  \right\}  dt+h\left(  X\left(  T\right)
\right)  -h\left(  \bar{X}\left(  T\right)  \right)  \right. \\
&  +\bar{g}\left(  \mathbb{E}\left[  X\left(  T\right)  \right]  \right)
-\bar{g}\left(  \mathbb{E}\left[  \bar{X}\left(  T\right)  \right]  \right)
\end{align*}
By the convexity of $h\left(  \cdot\right)  $ and $\bar{g}\left(
\cdot\right)  $, we have%
\[
\mathbb{E}\left[  h\left(  X\left(  T\right)  \right)  -h\left(  \bar
{X}\left(  T\right)  \right)  \right]  \geq\mathbb{E}\left[  \left(  X\left(
T\right)  -\bar{X}\left(  T\right)  \right)  ^{\top}h_{x}\left(  \bar
{X}\left(  T\right)  \right)  \right]
\]
and%
\[
G\left(  \mathbb{E}\left[  X\left(  T\right)  \right]  \right)  -G\left(
\mathbb{E}\left[  \bar{X}\left(  T\right)  \right]  \right)  \geq
\mathbb{E}\left[  \left(  X\left(  T\right)  -\bar{X}\left(  T\right)
\right)  ^{\top}G_{\bar{x}}\left(  \mathbb{E}\left[  \bar{X}\left(  T\right)
\right]  \right)  \right]  .
\]

Accordingly, by the terminal condition in the BSDE (\ref{Eq30}) we obtain that%
\begin{align}
&  \mathbf{J}\left(  u\left(  \cdot\right)  \right)  -\mathbf{J}\left(
\bar{u}\left(  \cdot\right)  \right) \nonumber\\
&  \geq\mathbb{E}\left[  \int_{0}^{T}\left\{  f^{u\left(  \cdot\right)
}\left(  t,X\left(  t\right)  \right)  -f^{\bar{u}\left(  \cdot\right)
}\left(  t,\bar{X}\left(  t\right)  \right)  \right\}  dt+\left(  X\left(
T\right)  -\bar{X}\left(  T\right)  \right)  ^{\top}\bar{p}\left(  T\right)
\right]  \text{.} \label{Eq31}%
\end{align}

From It\^{o}'s lemma, applied to $t\mapsto\left(  X\left(  t\right)  -\bar
{X}\left(  t\right)  \right)  ^{\top}\bar{p}\left(  T\right)  $ on time
interval $\left[  0,T\right]  $, it follows that%
\begin{align}
&  \mathbb{E}\left[  \left(  X\left(  T\right)  -\bar{X}\left(  T\right)
\right)  ^{\top}\bar{p}\left(  T\right)  \right] \nonumber\\
&  =\mathbb{E}\left[  \int_{0}^{T}\left\{  \left\langle b^{u\left(
\cdot\right)  }\left(  t,X\left(  t\right)  \right)  -b^{\bar{u}\left(
\cdot\right)  }\left(  t,\bar{X}\left(  t\right)  \right)  ,\bar{p}\left(
t\right)  \right\rangle \right.  \right. \nonumber\\
&  +\sum\limits_{j=1}^{d}\left\langle \sigma^{j}\left(  t,X\left(  t\right)
\right)  -\sigma^{j}\left(  t,\bar{X}\left(  t\right)  \right)  ,\bar{q}%
_{j}\left(  t\right)  \right\rangle \nonumber\\
&  \left.  \left.  -\left\langle X\left(  t\right)  -\bar{X}\left(  t\right)
,\mathbb{H}_{x}\left(  t,\bar{X}\left(  t\right)  ,\bar{u}\left(  t\right)
,\bar{p}\left(  t\right)  ,\bar{q}\left(  t\right)  \right)  \right\rangle
\right\}  dt\right]  . \label{Eq32}%
\end{align}

On the other hand, by the definition of the Hamiltonian $\mathbb{H}$, we have%
\begin{align}
&  \mathbb{E}\left[  \int_{0}^{T}\left\{  f\left(  t,X\left(  t\right)
,u\left(  t\right)  \right)  -f\left(  t,\bar{X}\left(  t\right)  ,\bar
{u}\left(  t\right)  \right)  \right\}  dt\right] \nonumber\\
&  =\mathbb{E}\left[  \int_{0}^{T}\left\{  \mathbb{H}\left(  t,X\left(
t\right)  ,u\left(  t\right)  ,\bar{p}\left(  t\right)  ,\bar{q}\left(
t\right)  \right)  -\mathbb{H}\left(  t,\bar{X}\left(  t\right)  ,\bar
{u}\left(  t\right)  ,\bar{p}\left(  t\right)  ,\bar{q}\left(  t\right)
\right)  \right.  \right. \nonumber\\
&  -\sum\limits_{j=1}^{d}\left\langle \sigma^{j}\left(  t,X\left(  t\right)
\right)  -\sigma^{j}\left(  t,\bar{X}\left(  t\right)  \right)  ,\bar{q}%
_{j}\left(  t\right)  \right\rangle \nonumber\\
&  \left.  \left.  -\left\langle b^{u\left(  \cdot\right)  }\left(  t,X\left(
t\right)  \right)  -b^{\bar{u}\left(  \cdot\right)  }\left(  t,\bar{X}\left(
t\right)  \right)  ,\bar{p}\left(  t\right)  \right\rangle \right\}
dt\right]  \text{.} \label{Eq33}%
\end{align}

Invoking (\ref{Eq32}) and (\ref{Eq33}) into (\ref{Eq31}), we obtain that%
\begin{align}
&  \mathbf{J}\left(  u\left(  \cdot\right)  \right)  -\mathbf{J}\left(
\bar{u}\left(  \cdot\right)  \right) \nonumber\\
&  \geq\mathbb{E}\left[  \int_{0}^{T}\left\{  \mathbb{H}\left(  t,X\left(
t\right)  ,u\left(  t\right)  ,\bar{p}\left(  t\right)  ,\bar{q}\left(
t\right)  \right)  -\mathbb{H}\left(  t,\bar{X}\left(  t\right)  ,\bar
{u}\left(  t\right)  ,\bar{p}\left(  t\right)  ,\bar{q}\left(  t\right)
\right)  \right.  \right. \nonumber\\
&  \left.  \left.  -\left\langle X\left(  t\right)  -\bar{X}\left(  t\right)
,\mathbb{H}_{x}\left(  t,\bar{X}\left(  t\right)  ,\bar{u}\left(  t\right)
,\bar{p}\left(  t\right)  ,\bar{q}\left(  t\right)  \right)  \right\rangle
\right\}  dt\right]  \text{.} \label{Eq34}%
\end{align}
By the convexity of the Hamiltonian $\mathbb{H}$, we have%
\begin{align}
&  \mathbb{H}\left(  t,X\left(  t\right)  ,u\left(  t\right)  ,\bar{p}\left(
t\right)  ,\bar{q}\left(  t\right)  \right)  -\mathbb{H}\left(  t,\bar
{X}\left(  t\right)  ,\bar{u}\left(  t\right)  ,\bar{p}\left(  t\right)
,\bar{q}\left(  t\right)  \right) \nonumber\\
&  \geq\left\langle X\left(  t\right)  -\bar{X}\left(  t\right)
,\mathbb{H}_{x}\left(  t,\bar{X}\left(  t\right)  ,\bar{u}\left(  t\right)
,\bar{p}\left(  t\right)  ,\bar{q}\left(  t\right)  \right)  \right\rangle
\nonumber\\
&  +\left\langle u\left(  t\right)  -\bar{u}\left(  t\right)  ,\mathbb{H}%
_{u}\left(  t,\bar{X}\left(  t\right)  ,\bar{u}\left(  t\right)  ,\bar
{p}\left(  t\right)  ,\bar{q}\left(  t\right)  \right)  \right\rangle \text{.}
\label{Eq35}%
\end{align}
Combining (\ref{Eq34})-(\ref{Eq35}) together with (\ref{4.9}), it follows that%
\begin{align*}
&  \mathbf{J}\left(  u\left(  \cdot\right)  \right)  -\mathbf{J}\left(
\bar{u}\left(  \cdot\right)  \right) \\
&  \geq\mathbb{E}\left[  \int_{0}^{T}\left\langle u\left(  t\right)  -\bar
{u}\left(  t\right)  ,\mathbb{H}_{u}\left(  t,\bar{X}\left(  t\right)
,\bar{u}\left(  t\right)  ,\bar{p}\left(  t\right)  ,\bar{q}\left(  t\right)
\right)  \right\rangle dt\right] \\
&  \geq0\text{,}%
\end{align*}
which means that $\bar{u}\left(  \cdot\right)  $ is an optimal control for
Problem (S).\eop

\begin{remark}
It should be noted that the assumptions in Theorem \ref{result08} can be made
weaker, but we do not focus on this here.\newline
\end{remark}

\subsection{Example}

In this subsection, we consider a simple example to illustrate our results.
For $n=d=1,$ consider the following controlled system:%
\begin{equation}
\left\{
\begin{array}
[c]{l}%
dX\left(  t\right)  =u\left(  t\right)  dt+dW\left(  t\right)  ,\text{ }%
t\in\left[  0,1\right]  ,\\
X\left(  0\right)  =0,
\end{array}
\right.  \label{Eq80**}%
\end{equation}
with the control domain being $U=\left[  -2,2\right]  $ and the cost
functional being%
\begin{equation}
\mathbf{J}\left(  u\left(  \cdot\right)  \right)  :=\frac{1}{2}\mathbb{E}%
\left[  \int_{0}^{1}\left(  u\left(  t\right)  +1\right)  ^{2}%
dt-e^{-\mathbb{E}\left[  X\left(  T\right)  \right]  ^{2}}\right]  .
\label{Eq80*}%
\end{equation}

We want to address the following stochastic control problem.

\subparagraph{Problem (E).}

\textit{Minimize (\ref{Eq80*}) over }$\mathcal{U}\left[  0,1\right]  .$

Suppose that $\left(  \bar{u}\left(  \cdot\right)  ,\bar{X}\left(
\cdot\right)  \right)  $ is an optimal pair (which we are going to identify),
then according to Theorem \ref{result05} the optimal 6-tuple $\left(  \bar
{X}\left(  \cdot\right)  ,\bar{u}\left(  \cdot\right)  ,\bar{g}\left(
\cdot,\cdot\right)  ,\bar{\eta}\left(  \cdot,\cdot\right)  ,\bar{\theta
}\left(  \cdot,\cdot\right)  ,\bar{\psi}\left(  \cdot,\cdot\right)  \right)  $
satisfies the system of FBSPDEs%
\begin{equation}
\left\{
\begin{array}
[c]{l}%
d\bar{X}\left(  t\right)  =\bar{u}\left(  t\right)  dt+dW\left(  t\right)
,\text{ }t\in\left[  0,T\right]  ,\\
d\bar{\theta}\left(  t,x\right)  =-\left\{  \bar{u}\left(  t\right)
\bar{\theta}_{x}\left(  t,x\right)  +\frac{1}{2}\bar{\theta}_{xx}\left(
t,x\right)  +\bar{\psi}_{x}\left(  t,x\right)  +\frac{1}{2}\left(  \bar
{u}\left(  t\right)  +1\right)  ^{2}\right\}  dt\\
\text{ \ \ \ \ \ \ \ \ \ \ \ \ \ \ \ }+\bar{\psi}\left(  t,x\right)  ^{\top
}dW\left(  t\right)  \text{, }\left(  t,x\right)  \in\left[  0,T\right]
\times\mathbb{R}^{n}\text{,}\\
d\bar{g}\left(  t,x\right)  =-\left\{  \bar{u}\left(  t\right)  \bar{g}%
_{x}\left(  t,x\right)  +\frac{1}{2}\bar{g}_{xx}\left(  t,x\right)  +\bar
{\eta}_{x}\left(  t,x\right)  \right\}  dt\\
\text{ \ \ \ \ \ \ \ \ \ \ \ }+\bar{\eta}\left(  t,x\right)  ^{\top}dW\left(
t\right)  \text{, }\left(  t,x\right)  \in\left[  0,T\right]  \times
\mathbb{R}^{n}\text{,}\\
\bar{X}\left(  0\right)  =0,\text{ }\bar{\theta}\left(  T,x\right)  =0,\text{
}\bar{g}\left(  T,x\right)  =x,\text{ for }x\in%
\mathbb{R}
,
\end{array}
\right.  \label{Eq66}%
\end{equation}
with the minimum condition: For a.e. $t\in\left[  0,1\right]  ,$%
\begin{equation}
\bar{u}\left(  t\right)  \in\arg\min\limits_{u\in\left[  -2,2\right]
}\left\{  u\left(  \bar{\theta}_{x}\left(  t,\bar{X}\left(  t\right)  \right)
+\mathbb{E}\left[  \bar{X}\left(  1\right)  \right]  e^{-\mathbb{E}\left[
\bar{X}\left(  1\right)  \right]  ^{2}}\bar{g}_{x}\left(  t,\bar{X}\left(
t\right)  \right)  \right)  +\frac{1}{2}\left(  u+1\right)  ^{2}\right\}
\text{.} \label{Eq67}%
\end{equation}
Moreover, the objective value of $\bar{u}\left(  \cdot\right)  $ is given by%
\begin{equation}
\mathbf{J}\left(  \bar{u}\left(  \cdot\right)  \right)  =\bar{\theta}\left(
0,0\right)  -\frac{1}{2}e^{-g\left(  0,0\right)  ^{2}}. \label{Eq68}%
\end{equation}

To solve the above system, we consider the following ansatz: For all $\left(
t,x\right)  \in\left[  0,T\right]  \times%
\mathbb{R}
,$%
\begin{equation}
\bar{\theta}\left(  t,x\right)  =M\left(  t\right)  \text{ } \label{Eq69}%
\end{equation}
and%
\begin{equation}
\bar{g}\left(  t,x\right)  =N\left(  t\right)  x+L\left(  t\right)  ,
\label{Eq69*}%
\end{equation}
where $\left(  M\left(  \cdot\right)  ,V_{M}\left(  \cdot\right)  \right)  $,
$\left(  N\left(  \cdot\right)  ,V_{N}\left(  \cdot\right)  \right)  $ and
$\left(  L\left(  \cdot\right)  ,V_{L}\left(  \cdot\right)  \right)  $ are the
pairs of adapted processes which are assumed to satisfy the BSDEs,%
\[
\left\{
\begin{array}
[c]{l}%
dM\left(  t\right)  =-U_{M}\left(  t\right)  dt+V_{M}\left(  t\right)
dW\left(  t\right)  ,\\
M\left(  T\right)  =0,
\end{array}
\right.
\]%
\[
\left\{
\begin{array}
[c]{l}%
dN\left(  t\right)  =-U_{N}\left(  t\right)  dt+V_{N}\left(  t\right)
dW\left(  t\right)  ,\\
N\left(  T\right)  =1,
\end{array}
\right.
\]
and%
\[
\left\{
\begin{array}
[c]{l}%
dL\left(  t\right)  =-U_{L}\left(  t\right)  dt+V_{L}\left(  t\right)
dW\left(  t\right)  ,\\
L\left(  T\right)  =0.
\end{array}
\right.
\]
In this case, the partial derivatives of $\bar{\theta}\left(  t,x\right)  $
and $\bar{g}\left(  t,x\right)  $ are%
\begin{equation}
\bar{\theta}_{x}\left(  t,x\right)  \equiv0\text{, }\bar{\theta}_{xx}\left(
t,x\right)  \equiv0,\text{ }\bar{g}_{x}\left(  t,x\right)  \equiv N\left(
t\right)  \text{ and }\bar{g}_{xx}\left(  t,x\right)  \equiv0. \label{Eq70}%
\end{equation}
We would like to determine the equation that $\left(  M\left(  \cdot\right)
,V_{M}\left(  \cdot\right)  \right)  $, $\left(  N\left(  \cdot\right)
,V_{N}\left(  \cdot\right)  \right)  $ and $\left(  L\left(  \cdot\right)
,V_{L}\left(  \cdot\right)  \right)  $ should satisfy. To this end, we
differentiate (\ref{Eq69})-(\ref{Eq69*}) and compare them with (\ref{Eq66}),
we obtain that%
\begin{align}
&  -U_{M}\left(  t\right)  dt+V_{M}\left(  t\right)  dW\left(  t\right)
\nonumber\\
&  =-\left\{  \bar{\psi}_{x}\left(  t,x\right)  +\frac{1}{2}\left(  \bar
{u}\left(  t\right)  +1\right)  ^{2}\right\}  dt+\bar{\psi}\left(  t,x\right)
dW\left(  t\right)  \text{.} \label{Eq71}%
\end{align}
and%
\begin{align}
&  -\left\{  U_{N}\left(  t\right)  x+U_{L}\left(  t\right)  \right\}
dt+\left\{  V_{N}\left(  t\right)  x+V_{L}\left(  t\right)  \right\}
dW\left(  t\right)  ,\nonumber\\
&  =-\left\{  N\left(  t\right)  \bar{u}\left(  t\right)  +\bar{\eta}%
_{x}\left(  t,x\right)  \right\}  dt+\bar{\eta}\left(  t,x\right)  dW\left(
t\right)  \text{.} \label{Eq71*}%
\end{align}
\bigskip Thus for all $\left(  t,x\right)  \in\left[  0,T\right]  \times%
\mathbb{R}
$,%
\begin{equation}
\bar{\psi}\left(  t,x\right)  =V_{M}\left(  t\right)  \label{Eq72}%
\end{equation}
and%
\begin{equation}
\bar{\eta}\left(  t,x\right)  =\left\{  V_{N}\left(  t\right)  x+V_{L}\left(
t\right)  \right\}  . \label{Eq72*}%
\end{equation}
Consequently,%
\begin{equation}
\bar{\psi}_{x}\left(  t,x\right)  \equiv0 \label{Eq73}%
\end{equation}
and%
\begin{equation}
\bar{\eta}_{x}\left(  t,x\right)  \equiv V_{N}\left(  t\right)  .
\label{Eq73*}%
\end{equation}

Next, comparing the $dt$ terms in (\ref{Eq71})-(\ref{Eq71*}) and using the
equalities (\ref{Eq73})- (\ref{Eq73*}), we obtain that%
\[
U_{M}\left(  t\right)  =\frac{1}{2}\left(  \bar{u}\left(  t\right)  +1\right)
^{2}%
\]
and%
\[
U_{N}\left(  t\right)  x+U_{L}\left(  t\right)  =\left\{  N\left(  t\right)
\bar{u}\left(  t\right)  +V_{N}\left(  t\right)  \right\}
\]
which leads to the following BSDEs%
\begin{equation}
\left\{
\begin{array}
[c]{l}%
dM\left(  t\right)  =-\frac{1}{2}\left(  \bar{u}\left(  t\right)  +1\right)
^{2}dt+V_{M}\left(  t\right)  dW\left(  t\right)  ,\\
M\left(  T\right)  =0,
\end{array}
\right.  \label{Eq75}%
\end{equation}%
\begin{equation}
\left\{
\begin{array}
[c]{l}%
dN\left(  t\right)  =V_{N}\left(  t\right)  dW\left(  t\right)  ,\\
N\left(  T\right)  =1,
\end{array}
\right.  \label{Eq75*}%
\end{equation}
and%
\begin{equation}
\left\{
\begin{array}
[c]{l}%
dL\left(  t\right)  =-\left\{  N\left(  t\right)  \bar{u}\left(  t\right)
+V_{N}\left(  t\right)  \right\}  dt+V_{L}\left(  t\right)  dW\left(
t\right)  ,\\
L\left(  T\right)  =0.
\end{array}
\right.  \label{Eq75**}%
\end{equation}
Equation (\ref{Eq70}) can be easily solved, whose solution is given by%
\[
\left(  N\left(  t\right)  ,V_{N}\left(  t\right)  \right)  \equiv\left(
1,0\right)  \text{, }\forall t\in\left[  0,1\right]  .
\]
Moreover by taking (\ref{Eq70}) into (\ref{Eq67}) we obtain that,%
\[
\bar{u}\left(  t\right)  \in\arg\min\limits_{u\in\left[  -2,2\right]
}\left\{  u\mathbb{E}\left[  \bar{X}\left(  1\right)  \right]  e^{-\mathbb{E}%
\left[  \bar{X}\left(  1\right)  \right]  ^{2}}+\frac{1}{2}\left(  u+1\right)
^{2}\right\}
\]
which suggests that%
\begin{equation}
\bar{u}\left(  t\right)  \equiv-\mathbb{E}\left[  \bar{X}\left(  1\right)
\right]  e^{-\mathbb{E}\left[  \bar{X}\left(  1\right)  \right]  ^{2}}%
-1\in\left]  -2,+2\right[  \text{, }\forall t\in\left[  0,1\right]  .
\label{Eq74}%
\end{equation}

In order to see that $\bar{u}\left(  t\right)  \equiv-\mathbb{E}\left[
\bar{X}\left(  T\right)  \right]  e^{-\mathbb{E}\left[  \bar{X}\left(
T\right)  \right]  ^{2}}-1$ is indeed optimal, we note that%
\[
\mathbb{H}\left(  t,x,u,p,q\right)  =up+q+\frac{1}{2}\left(  u+1\right)  ^{2}%
\]
is convex with respect to $\left(  x,u\right)  $. Moreover $h\left(  x\right)
=0$ and $G\left(  \bar{x}\right)  =-e^{-\bar{x}^{2}}$ are convex. So the
optimality follows from the sufficient condition (Theorem \ref{result08}%
).\newline

Taking the optimal control (\ref{Eq74}) into the state equation (\ref{Eq80**}%
), we obtain that%
\begin{equation}
\left\{
\begin{array}
[c]{l}%
d\bar{X}\left(  t\right)  =-\left\{  \mathbb{E}\left[  \bar{X}\left(
1\right)  \right]  e^{-\mathbb{E}\left[  \bar{X}\left(  1\right)  \right]
^{2}}+1\right\}  dt+dW\left(  t\right)  ,\text{ }t\in\left[  0,1\right]  ,\\
X\left(  0\right)  =0,
\end{array}
\right.  \label{Eq77}%
\end{equation}
Accordingly we have,%
\[
\left\{
\begin{array}
[c]{l}%
\mathbb{E}\left[  \bar{X}\left(  t\right)  \right]  =-\left\{  \mathbb{E}%
\left[  \bar{X}\left(  T\right)  \right]  e^{-\mathbb{E}\left[  \bar{X}\left(
T\right)  \right]  ^{2}}+1\right\}  dt,\text{ }t\in\left[  0,1\right]  ,\\
\mathbb{E}\left[  \bar{X}\left(  0\right)  \right]  =0,
\end{array}
\right.
\]
which leads to
\[
\mathbb{E}\left[  \bar{X}\left(  t\right)  \right]  =-\left\{  \mathbb{E}%
\left[  \bar{X}\left(  T\right)  \right]  e^{-\mathbb{E}\left[  \bar{X}\left(
T\right)  \right]  ^{2}}+1\right\}  t,\text{ }\forall t\in\left[  0,1\right]
,
\]
Setting $t=1$ in the above, we obtain that $\mathbb{E}\left[  \bar{X}\left(
1\right)  \right]  =\mathbf{X}_{\ast}$ coincides with the unique solution of
the following equation,%
\[
\mathbf{X}_{\ast}+\mathbf{X}_{\ast}e^{-\mathbf{X}_{\ast}^{2}}+1=0.
\]

Now invoking (\ref{Eq74}) into (\ref{Eq75})-(\ref{Eq75**}), the BSDEs
satisfied by $\left(  M\left(  \cdot\right)  ,V_{M}\left(  \cdot\right)
\right)  $ and $\left(  L\left(  \cdot\right)  ,V_{L}\left(  \cdot\right)
\right)  $ reduce to%
\[
\left\{
\begin{array}
[c]{l}%
dM\left(  t\right)  =-\frac{1}{2}\mathbf{X}_{\ast}^{2}e^{-2\mathbf{X}_{\ast
}^{2}}dt+V_{M}\left(  t\right)  dW\left(  t\right)  ,\\
M\left(  1\right)  =0;
\end{array}
\right.
\]
and%
\[
\left\{
\begin{array}
[c]{l}%
dL\left(  t\right)  =\left\{  \mathbf{X}_{\ast}e^{-\mathbf{X}_{\ast}^{2}%
}+1\right\}  dt+V_{L}\left(  t\right)  dW\left(  t\right)  ,\\
L\left(  1\right)  =0.
\end{array}
\right.
\]
These equations can be easily solved, whose solutions are%
\begin{equation}
\text{ }\left(  M\left(  t\right)  ,V_{M}\left(  t\right)  \right)
\equiv\left(  \frac{1}{2}\mathbf{X}_{\ast}^{2}e^{-2\mathbf{X}_{\ast}^{2}%
}\left(  1-t\right)  ,0\right)  \text{, }\forall t\in\left[  0,1\right]  .
\label{Eq76}%
\end{equation}
and%
\begin{equation}
\text{ }\left(  L\left(  t\right)  ,V_{L}\left(  t\right)  \right)
\equiv\left(  -\left(  \mathbf{X}_{\ast}e^{-\mathbf{X}_{\ast}^{2}}+1\right)
\left(  1-t\right)  ,0\right)  \text{, }\forall t\in\left[  0,1\right]  .
\label{Eq76*}%
\end{equation}
Taking (\ref{Eq76})-(\ref{Eq76*}) into (\ref{Eq69})-(\ref{Eq69*}) we obtain
that the objective value of $\bar{u}\left(  \cdot\right)  $ is given by%
\[
\mathbf{J}\left(  \bar{u}\left(  \cdot\right)  \right)  =\frac{1}{2}%
\mathbf{X}_{\ast}^{2}e^{-2\mathbf{X}_{\ast}^{2}}-\frac{1}{2}e^{-\left(
\mathbf{X}_{\ast}e^{-\mathbf{X}_{\ast}^{2}}+1\right)  ^{2}}.
\]

\begin{remark}
Summarizing the preceding analysis, we have the following results:\newline(i)
The optimal control of Problem (E) is given by%
\[
\bar{u}\left(  t\right)  \equiv-\mathbf{X}_{\ast}e^{-\mathbf{X}_{\ast}^{2}}-1
\]
with the corresponding state%
\[
\bar{X}\left(  t\right)  \equiv W\left(  t\right)  -t\left(  \mathbf{X}_{\ast
}e^{-\mathbf{X}_{\ast}^{2}}+1\right)  ,
\]
where $\mathbf{X}_{\ast}$ satisfies the following equation
\[
\mathbf{X}_{\ast}+\mathbf{X}_{\ast}e^{-\mathbf{X}_{\ast}}+1=0.
\]
(ii) The corresponding solution of the BSPDEs are given by,%
\begin{align*}
\bar{\theta}\left(  t,x\right)   &  \equiv\frac{1}{2}\mathbf{X}_{\ast}%
^{2}e^{-2\mathbf{X}_{\ast}^{2}}\left(  1-t\right)  ,\\
\bar{\psi}\left(  t,x\right)   &  \equiv0.
\end{align*}
and%
\begin{align*}
\bar{g}\left(  t,x\right)   &  \equiv x-\left(  \mathbf{X}_{\ast
}e^{-\mathbf{X}_{\ast}^{2}}+1\right)  \left(  1-t\right)  ,\\
\bar{\eta}\left(  t,x\right)   &  \equiv0.
\end{align*}
(iii) The objective value of $\bar{u}\left(  \cdot\right)  $ is given by%
\begin{align*}
\mathbf{J}\left(  \bar{u}\left(  \cdot\right)  \right)   &  =\inf_{u\left(
\cdot\right)  \in\mathcal{U}\left[  0,1\right]  }\mathbf{J}\left(  u\left(
\cdot\right)  \right) \\
&  =\frac{1}{2}\mathbf{X}_{\ast}^{2}e^{-2\mathbf{X}_{\ast}^{2}}-\frac{1}%
{2}e^{-\left(  \mathbf{X}_{\ast}e^{-\mathbf{X}_{\ast}^{2}}+1\right)  ^{2}}.
\end{align*}
Again, we emphasize that the new version of the SMP permits us to derive the
objective value $\mathbf{J}\left(  \bar{u}\left(  \cdot\right)  \right)  $,
which is different from the traditional SMP approach.
\end{remark}

\begin{remark}
Approximately, we have
\begin{align*}
\mathbb{E}\left[  \bar{X}\left(  1\right)  \right]   &  =\mathbf{X}_{\ast
}\simeq-0.584\,62;\\
\bar{u}\left(  t\right)   &  \simeq0.58462,\text{ }\forall t\in\left[
0,1\right]  ;\\
\mathbf{J}\left(  \bar{u}\left(  \cdot\right)  \right)   &  \simeq-0.29.
\end{align*}

\end{remark}

\section{A Class of forward-backward stochastic partial differential
equations\label{section7}}

In this section, motivated by the system of FBSDEs (\ref{Eq326}), we are going
to study the solvability of the following class of forward-backward stochastic
partial differential equations:%
\begin{equation}
\left\{
\begin{array}
[c]{l}%
dX\left(  t\right)  =\bar{b}\left(  t,X\left(  t\right)  ,p_{x}\left(
t,X\left(  t\right)  \right)  \right)  dt+\bar{\sigma}\left(  t,X\left(
t\right)  \right)  dW\left(  t\right)  \text{,\ for }t\in\left[  0,T\right]
\text{,}\\
dp\left(  t,x\right)  =-\left\{  \left\langle p_{x}\left(  t,x\right)
,\bar{b}\left(  t,x,p_{x}\left(  t,X\left(  t\right)  \right)  \right)
\right\rangle +\frac{1}{2}\text{\textbf{tr}}\left[  \bar{\sigma}\left(
t,x\right)  \bar{\sigma}\left(  t,x\right)  ^{\top}p_{xx}\left(  t,X\left(
t\right)  \right)  \right]  \right.  \text{\ \ \ }\\
\ \ \ \ \ \ \ \ \ \ \ \ \ \ \ +\left.  \text{\textbf{tr}}\left[  q_{x}\left(
t,x\right)  \bar{\sigma}\left(  t,x\right)  \right]  +\bar{f}\left(
t,x,p_{x}\left(  t,X\left(  t\right)  \right)  \right)  \right\}  dt\\
\ \ \ \ \ \ \ \ \ \ \ \ \ \ \ +\bar{q}\left(  t,x\right)  ^{\top}dW\left(
t\right)  \text{, for }\left(  t,x\right)  \in\left[  0,T\right]
\times\mathbb{R}^{n}\text{,}\\
X\left(  0\right)  =x_{0},\text{ }p\left(  T,x\right)  =\bar{F}\left(
x\right)  \text{, for }x\in\mathbb{R}^{n}\text{.}%
\end{array}
\right.  \label{Eq200}%
\end{equation}

In the above, $X\left(  \cdot\right)  $ is the unknown process, $\left(
p\left(  \cdot,\cdot\right)  ,q\left(  \cdot,\cdot\right)  \right)  $ are the
unknown random fields, and they are required to be $\left(  \mathcal{F}%
_{t}\right)  _{t\in\left[  0,T\right]  }$-adapted; $\bar{b}:\left[
0,T\right]  \times\mathbb{R}^{n}\times\mathbb{R}^{n}\rightarrow\mathbb{R}^{n}%
$, $\bar{\sigma}:\left[  0,T\right]  \times\mathbb{R}^{n}\rightarrow
\mathbb{R}^{n\times d}$, $\bar{f}:\left[  0,T\right]  \times\mathbb{R}%
^{n}\times\mathbb{R}^{n}\rightarrow\mathbb{R}$ and $\bar{F}:\mathbb{R}%
^{n}\rightarrow\mathbb{R}$ are given deterministic measurable functions. The
main feature of the above system is that $X\left(  \cdot\right)  $ satisfies a
forward SDE and the pair of random fields $\left(  p\left(  \cdot
,\cdot\right)  ,q\left(  \cdot,\cdot\right)  \right)  $ satisfies a nonlinear BSPDE.

We introduce the following definition.

\begin{definition}
A triple $\left(  X\left(  \cdot\right)  ,p\left(  \cdot,\cdot\right)
,q\left(  \cdot,\cdot\right)  \right)  $ is called a \textit{solution} of
(\ref{Eq200}) if the following holds%
\[
\left\{
\begin{array}
[c]{l}%
X\left(  \cdot\right)  \in\mathcal{C}_{\mathcal{F}}^{2}\left(  0,T;\mathbb{R}%
^{n}\right) \\
p\left(  \cdot,\cdot\right)  \in\mathcal{C}_{\mathcal{F}}^{2}\left(
0,T;C^{2}\left(  \bar{B}_{R};\mathbb{R}\right)  \right)  ,\\
q\left(  \cdot,\cdot\right)  \in\mathcal{L}_{\mathcal{F}}^{2}\left(
0,T;C^{1}\left(  \bar{B}_{R};\mathbb{R}^{d}\right)  \right)  ,
\end{array}
\right.  \forall R>0\text{,}%
\]
such that the following holds, for all $\left(  t,x\right)  \in\left[
0,T\right]  \times\mathbb{R}^{n}$, almost surely:%
\[
\left\{
\begin{array}
[c]{l}%
X\left(  t\right)  =x_{0}+\int_{0}^{t}\bar{b}\left(  \tau,X\left(
\tau\right)  ,p_{x}\left(  \tau,X\left(  \tau\right)  \right)  \right)
d\tau+\int_{0}^{t}\bar{\sigma}\left(  \tau,X\left(  \tau\right)  \right)
dW\left(  \tau\right)  ,\\
p\left(  t,x\right)  =\bar{F}\left(  x\right)  +\int_{t}^{T}\left\{
\left\langle p_{x}\left(  \tau,x\right)  ,\bar{b}\left(  \tau,x,p_{x}\left(
\tau,X\left(  \tau\right)  \right)  \right)  \right\rangle \right. \\
\text{ \ \ \ \ \ \ \ \ \ \ \ \ \ }+\frac{1}{2}\text{\textbf{tr}}\left[
\bar{\sigma}\left(  \tau,x\right)  \bar{\sigma}\left(  \tau,x\right)  ^{\top
}p_{xx}\left(  \tau,x\right)  \right] \\
\ \ \ \ \ \ \ \ \ \ \ \ \ \ +\left.  \text{\textbf{tr}}\left[  q_{x}\left(
\tau,x\right)  \bar{\sigma}\left(  \tau,x\right)  \right]  +\bar{f}\left(
\tau,x,p_{x}\left(  \tau,X\left(  \tau\right)  \right)  \right)  \right\}
d\tau\\
\ \ \ \ \ \ \ \ \ \ \ \ \ \ -\int_{t}^{T}\bar{q}\left(  \tau,x\right)  ^{\top
}dW\left(  \tau\right)  .
\end{array}
\right.
\]

\end{definition}

\subsection{A version of three-step scheme, a heuristic derivation}

Now, essentially inspired by the idea of the classical four-step scheme
introduced in Ma et al. \cite{MPY} (see also e.g. Ma and Yong \cite{MY} or
Yong and Zhou \cite{YongZhou}), we are going to introduce a method for solving
the FBSPDE (\ref{Eq200}) over any time duration $\left[  0,T\right]  $.

As in~[\cite{YongZhou}, Section 7.5.2] let us give a heuristic derivation
first. Suppose that $\left(  X\left(  \cdot\right)  ,p\left(  \cdot
,\cdot\right)  ,q\left(  \cdot,\cdot\right)  \right)  $ is an adapted solution
to (\ref{Eq200}). We assume that $p\left(  t,x\right)  $ and $X\left(
t\right)  $ are related by%
\begin{equation}
p\left(  t,x\right)  =\theta\left(  t,x,X\left(  t\right)  \right)  ,\text{
for any }\left(  t,x\right)  \in\left[  0,T\right]  \times%
\mathbb{R}
^{n}\text{,} \label{Eq201}%
\end{equation}
where $\theta\left(  t,x,y\right)  $ is some function to be determined such
that for any $\left(  x,y\right)  \in%
\mathbb{R}
^{n}\times%
\mathbb{R}
^{n}$,
\[
\theta\left(  T,x,y\right)  =F\left(  x\right)  .
\]

We suppose that $\theta\left(  t,x,y\right)  $ is $C^{1}$ in $t$ and \ $C^{2}$
in $\left(  x,y\right)  $. Accordingly, we have for any $\left(  t,x\right)
\in\left[  0,T\right]  \times%
\mathbb{R}
^{n}$,%
\[
p_{x}\left(  t,x\right)  =\theta_{x}\left(  t,x,X\left(  t\right)  \right)
\text{,}%
\]
and%
\[
p_{xx}\left(  t,x\right)  =\theta_{xx}\left(  t,x,X\left(  t\right)  \right)
\text{.}%
\]

Then by It\^{o}'s formula, we differentiate (\ref{Eq201}) and compare it with
(\ref{Eq200}) to obtain that%
\begin{align}
dp\left(  t,x\right)   &  =d\theta\left(  t,x,X\left(  t\right)  \right)
\nonumber\\
&  =\left\{  \theta_{t}\left(  t,x,X\left(  t\right)  \right)  +\left\langle
\theta_{y}\left(  t,x,X\left(  t\right)  \right)  ,\bar{b}\left(  t,X\left(
t\right)  ,\theta_{x}\left(  t,X\left(  t\right)  ,X\left(  t\right)  \right)
\right)  \right\rangle \right. \nonumber\\
&  \left.  +\frac{1}{2}\text{\textbf{tr}}\left[  \bar{\sigma}\left(
t,X\left(  t\right)  \right)  \bar{\sigma}\left(  t,X\left(  t\right)
\right)  ^{\top}\theta_{yy}\left(  t,x,X\left(  t\right)  \right)  \right]
\right\}  dt\nonumber\\
&  +\left\langle \theta_{y}\left(  t,x,X\left(  t\right)  \right)
,\bar{\sigma}\left(  t,X\left(  t\right)  \right)  dW\left(  t\right)
\right\rangle \nonumber\\
&  =-\left\{  \left\langle \theta_{x}\left(  t,x,X\left(  t\right)  \right)
,\bar{b}\left(  t,x,\theta_{x}\left(  t,X\left(  t\right)  ,X\left(  t\right)
\right)  \right)  \right\rangle \right. \nonumber\\
&  +\frac{1}{2}\text{\textbf{tr}}\left[  \bar{\sigma}\left(  t,x\right)
\bar{\sigma}\left(  t,x\right)  ^{\top}\theta_{xx}\left(  t,x,X\left(
t\right)  \right)  \right] \nonumber\\
&  +\left.  \text{\textbf{tr}}\left[  q_{x}\left(  t,x\right)  \bar{\sigma
}\left(  t,x\right)  \right]  +\bar{f}\left(  t,x,\theta_{x}\left(  t,X\left(
t\right)  ,X\left(  t\right)  \right)  \right)  \right\}  dt\ \nonumber\\
&  +q\left(  t,x\right)  ^{\top}dW\left(  t\right)  , \label{Eq203}%
\end{align}
Accordingly, we have%
\[
q\left(  t,x\right)  \equiv\bar{\sigma}\left(  t,X\left(  t\right)  \right)
^{\top}\theta_{y}\left(  t,x,X\left(  t\right)  \right)  .
\]
Thus
\begin{equation}
q_{x}\left(  t,x\right)  \equiv\bar{\sigma}\left(  t,X\left(  t\right)
\right)  ^{\top}\theta_{yx}\left(  t,x,X\left(  t\right)  \right)  .
\label{Eq204}%
\end{equation}

Now compare the dt terms in (\ref{Eq203}) and using (\ref{Eq204}), we obtain
that%
\begin{align*}
0  &  =\theta_{t}\left(  t,x,X\left(  t\right)  \right)  +\left\langle
\theta_{y}\left(  t,x,X\left(  t\right)  \right)  ,\bar{b}\left(  t,X\left(
t\right)  ,\theta_{x}\left(  t,X\left(  t\right)  ,X\left(  t\right)  \right)
\right)  \right\rangle \\
&  +\frac{1}{2}\text{\textbf{tr}}\left[  \bar{\sigma}\left(  t,X\left(
t\right)  \right)  \bar{\sigma}\left(  t,X\left(  t\right)  \right)  ^{\top
}\theta_{yy}\left(  t,x,X\left(  t\right)  \right)  \right]  +\left\langle
\theta_{y}\left(  t,x,X\left(  t\right)  \right)  ,\bar{\sigma}\left(
t,X\left(  t\right)  \right)  dW\left(  t\right)  \right\rangle \\
&  +\left\{  \left\langle \theta_{x}\left(  t,x,X\left(  t\right)  \right)
,\bar{b}\left(  t,x,\theta_{x}\left(  t,X\left(  t\right)  ,X\left(  t\right)
\right)  \right)  \right\rangle +\frac{1}{2}\text{\textbf{tr}}\left[
\bar{\sigma}\left(  t,x\right)  \bar{\sigma}\left(  t,x\right)  ^{\top}%
\theta_{xx}\left(  t,x,X\left(  t\right)  \right)  \right]  \right. \\
&  +\text{\textbf{tr}}\left[  \bar{\sigma}\left(  t,X\left(  t\right)
\right)  ^{\top}\theta_{yx}\left(  t,x,X\left(  t\right)  \right)  \bar
{\sigma}\left(  t,x\right)  \right]  +\bar{f}\left(  t,x,\theta_{x}\left(
t,X\left(  t\right)  ,X\left(  t\right)  \right)  \right)  \text{.}%
\end{align*}

The above argument suggests that we design the following three-step scheme to
solve the FBSDE (\ref{Eq200}).

\subparagraph{A Three-step scheme:}

\begin{enumerate}
\item[\textbf{Step 1.}] Solve the following nonlinear parabolic system for
$\theta\left(  t,x,y\right)  $:%
\begin{equation}
\left\{
\begin{array}
[c]{l}%
0=\theta_{t}\left(  t,x,y\right)  +\left\langle \theta_{y}\left(
t,x,y\right)  ,\bar{b}\left(  t,y,\theta_{x}\left(  t,y,y\right)  \right)
\right\rangle +\frac{1}{2}\text{\textbf{tr}}\left[  \bar{\sigma}\left(
t,y\right)  \bar{\sigma}\left(  t,y\right)  ^{\top}\theta_{yy}\left(
t,x,y\right)  \right] \\
\text{ \ \ }+\left\langle \theta_{x}\left(  t,x,y\right)  ,\bar{b}\left(
t,x,\theta_{x}\left(  t,y,y\right)  \right)  \right\rangle +\frac{1}%
{2}\text{\textbf{tr}}\left[  \bar{\sigma}\left(  t,x\right)  \bar{\sigma
}\left(  t,x\right)  ^{\top}\theta_{xx}\left(  t,x,y\right)  \right] \\
\text{ \ \ }+\text{\textbf{tr}}\left[  \bar{\sigma}\left(  t,y\right)  ^{\top
}\theta_{yx}\left(  t,x,y\right)  \bar{\sigma}\left(  t,x\right)  \right]
\text{\ }+\bar{f}\left(  t,x,\theta_{x}\left(  t,y,y\right)  \right)  \text{,
for }\left(  t,x,y\right)  \in\left[  0,T\right]  \times%
\mathbb{R}
^{n}\times%
\mathbb{R}
^{n},\\
\theta\left(  T,x,y\right)  =\bar{F}\left(  x\right)  ,\text{ for }\left(
t,x,y\right)  \in%
\mathbb{R}
^{n}\times%
\mathbb{R}
^{n}.
\end{array}
\right.  \label{Eq205}%
\end{equation}

\item[\textbf{Step 2.}] Use $\theta$ obtained in Steps 1 to solve the
following forward SDE:%
\begin{equation}
\left\{
\begin{array}
[c]{l}%
dX\left(  t\right)  =\tilde{b}\left(  t,X\left(  t\right)  \right)
dt+\bar{\sigma}\left(  t,X\left(  t\right)  \right)  dW\left(  t\right)
,\ \text{for }t\in\left[  0,T\right]  ,\\
X\left(  0\right)  =x_{0},
\end{array}
\right.  \label{Eq2014}%
\end{equation}
where%
\begin{equation}
\tilde{b}\left(  t,y\right)  :=\bar{b}\left(  t,y,\theta_{x}\left(
t,y,y\right)  \right)  \label{Eq210}%
\end{equation}

\item[\textbf{Step 3.}] Set%
\begin{equation}
\left\{
\begin{array}
[c]{l}%
p\left(  t,x\right)  :=\theta\left(  t,x,X\left(  t\right)  \right)  ,\text{
}\\
q\left(  t,x\right)  :=\bar{\sigma}\left(  t,X\left(  t\right)  \right)
^{\top}\theta_{y}\left(  t,x,X\left(  t\right)  \right)  ,\text{ for any
}\left(  t,x\right)  \in\left[  0,T\right]  \times%
\mathbb{R}
^{n}.
\end{array}
\right.  \label{Eq211}%
\end{equation}

\end{enumerate}

Should this scheme be realizable, $\left(  X\left(  \cdot\right)  ,p\left(
\cdot,\cdot\right)  ,q\left(  \cdot,\cdot\right)  \right)  $ would give an
adapted solution of (\ref{Eq200}). We have the following result.

\begin{theorem}
Assume that (\ref{Eq205}) admits a unique classical solution $\theta\left(
t,x,y\right)  $ with bounded $\theta_{x}\left(  t,x,y\right)  $, $\theta
_{y}\left(  t,x,y\right)  $, $\theta_{xx}\left(  t,x,y\right)  $, $\theta
_{yy}\left(  t,x,y\right)  $ and $\theta_{xy}\left(  t,x,y\right)  .$ Assume
further that the functions $\bar{b}\left(  t,y,p\right)  $ and $\bar{\sigma
}\left(  t,y\right)  $ are uniformly Lipschitz continuous in $\left(
y,p\right)  $ with $\bar{b}\left(  t,0,0\right)  $ and $\bar{\sigma}\left(
t,0\right)  $ being bounded. Then the triple $\left(  X\left(  \cdot\right)
,p\left(  \cdot,\cdot\right)  ,q\left(  \cdot,\cdot\right)  \right)  $
determined by (\ref{Eq2014}) and (\ref{Eq211}) is an adapted solution to
(\ref{Eq200}).
\end{theorem}

\bop Under our conditions both $\bar{\sigma}\left(  t,y\right)  $ and
$\tilde{b}\left(  t,y\right)  $ defined by (\ref{Eq210}) are uniformly
Lipschitz continuous in $y$. Therefore, for any $x_{0}\in%
\mathbb{R}
^{n}$, (\ref{Eq2014}) has a unique strong solution. Then, by defining
$p\left(  t,x\right)  $ and $q\left(  t,x\right)  $ via (\ref{Eq211}) and
applying It\^{o}'s formula, we can easily check that (\ref{Eq200}) is
satisfied. Hence, $\left(  X\left(  \cdot\right)  ,p\left(  \cdot
,\cdot\right)  ,q\left(  \cdot,\cdot\right)  \right)  $ is a solution of
(\ref{Eq200}).\eop

We now turn our attention on the Cauchy problem (\ref{Eq205}). First, let us
consider the following notations\footnote{Here we adopt the following
convention. For any $x=\left(  x_{1},...,x_{n}\right)  ^{\top}\in%
\mathbb{R}
^{n}$ and $y=\left(  y_{1},...,y_{n}\right)  ^{\top}\in%
\mathbb{R}
^{n},$ we consider the pair $\left(  x,y\right)  $ as an element of $%
\mathbb{R}
^{2n}$ and write:
\[
\left(  x,y\right)  :=\left(  x_{1},...,x_{n},y_{1},...,y_{n}\right)  \in%
\mathbb{R}
^{2n}.
\]
}%
\begin{align*}
\mathbf{X}  &  \mathbf{:}\mathbf{=}\left(  x,y\right)  \in%
\mathbb{R}
^{2n},\\
\varphi\left(  \mathbf{X}\right)   &  =\varphi\left(  x,y\right)  :=\left(
y,y\right)  \in%
\mathbb{R}
^{2n},\\
\theta_{\mathbf{X}}\left(  t,\mathbf{X}\right)   &  :=\left(
\begin{array}
[c]{c}%
\theta_{x}\left(  t,x,y\right) \\
\theta_{y}\left(  t,x,y\right)
\end{array}
\right)  ,\\
\theta_{\mathbf{XX}}\left(  t,\mathbf{X}\right)   &  :=\left(
\begin{array}
[c]{cc}%
\theta_{xx}\left(  t,x,y\right)  & \theta_{xy}\left(  t,x,y\right) \\
\theta_{yx}\left(  t,x,y\right)  & \theta_{yy}\left(  t,x,y\right)
\end{array}
\right)  .\\
\mathbf{b}\left(  t,\mathbf{X,}\theta_{\mathbf{X}}\left(  t,\varphi\left(
\mathbf{X}\right)  \right)  \right)   &  :=\left(
\begin{array}
[c]{c}%
\bar{b}\left(  t,x,\theta_{x}\left(  t,y,y\right)  \right) \\
\bar{b}\left(  t,y,\theta_{x}\left(  t,y,y\right)  \right)
\end{array}
\right)  ,\\
\mathbf{f}\left(  t,\mathbf{X,}\theta_{\mathbf{X}}\left(  t,\varphi\left(
\mathbf{X}\right)  \right)  \right)   &  :=f\left(  t,x,\theta_{x}\left(
t,y,y\right)  \right)  ,\\
\mathbf{F}\left(  \mathbf{X}\right)   &  :=\bar{F}\left(  x\right)  ,\\
\mathbf{a}\left(  t,\mathbf{X}\right)   &  :=\frac{1}{2}\left(
\begin{array}
[c]{cc}%
\bar{\sigma}\left(  t,x\right)  \bar{\sigma}\left(  t,x\right)  ^{\top} &
\bar{\sigma}\left(  t,y\right)  \bar{\sigma}\left(  t,x\right)  ^{\top}\\
\bar{\sigma}\left(  t,x\right)  \sigma\left(  t,y\right)  ^{\top} &
\bar{\sigma}\left(  t,y\right)  \bar{\sigma}\left(  t,y\right)  ^{\top}%
\end{array}
\right)  .
\end{align*}

Then it is not difficult to see that $\theta\left(  t,x,y\right)  $ is a
classical solution to (\ref{Eq205}), if and only if, $\theta\left(
t,\mathbf{X}\right)  $ is a classical solution of the following nonlinear
parabolic system%
\begin{equation}
\left\{
\begin{array}
[c]{l}%
0=\theta_{t}\left(  t,\mathbf{X}\right)  +\left\langle \theta_{\mathbf{X}%
}\left(  t,\mathbf{X}\right)  ,\mathbf{b}\left(  t,\mathbf{X,}\theta
_{\mathbf{X}}\left(  t,\varphi\left(  \mathbf{X}\right)  \right)  \right)
\right\rangle \\
\text{ \ \ \ }+\text{\textbf{tr}}\left[  \mathbf{a}\left(  t,\mathbf{X}%
\right)  \theta_{X\mathbf{X}}\left(  t,\mathbf{X}\right)  \right] \\
\text{ \ \ \ }+\mathbf{f}\left(  t,\mathbf{X},\theta_{\mathbf{X}}\left(
t,\varphi\left(  \mathbf{X}\right)  \right)  \right)  ,\text{ for }\left(
t,\mathbf{X}\right)  \in\left[  0,T\right]  \times%
\mathbb{R}
^{2n},\\
\theta\left(  T,\mathbf{X}\right)  =\mathbf{F}\left(  \mathbf{X}\right)
,\text{ for }\mathbf{X}\in%
\mathbb{R}
^{2n}.
\end{array}
\right.  \label{Eq212}%
\end{equation}

\subsection{Well-posedness of the parabolic PDE}

In this subsections, we discuss the well-posedness for the PDE (\ref{Eq212})
by adopting a fixed point method. Let us make some preparations.

Let $C^{\alpha}\left(
\mathbb{R}
^{2n}\right)  $ be the space of continuous functions $\varphi\left(
\cdot\right)  $ such that%
\[
\left\Vert \varphi\right\Vert _{\alpha}=\left\Vert \varphi\right\Vert
_{0}+\left\lfloor \varphi\right\rfloor _{\alpha}<\infty,
\]
where%
\[
\left\Vert \varphi\right\Vert _{0}=\sup_{\mathbf{X}\in%
\mathbb{R}
^{2n}}\left\vert \varphi\left(  \mathbf{X}\right)  \right\vert ,\text{
}\left\lfloor \varphi\right\rfloor _{\alpha}=\sup_{\mathbf{X,Y}\in%
\mathbb{R}
^{2n}}\frac{\left\vert \varphi\left(  \mathbf{X}\right)  -\varphi\left(
\mathbf{Y}\right)  \right\vert }{\left\Vert \mathbf{X-Y}\right\Vert ^{\alpha}%
}<\infty
\]

Further let $C^{1+\alpha}\left(
\mathbb{R}
^{2n}\right)  $ and $C^{1+\alpha}\left(
\mathbb{R}
^{2n}\right)  $ be the space of continuous functions $\varphi\left(
\cdot\right)  $ such that%
\[
\left\Vert \varphi\right\Vert _{1+\alpha}=\left\Vert \varphi\right\Vert
_{0}+\left\Vert \varphi_{x}\right\Vert _{0}+\left\lfloor \varphi
_{x}\right\rfloor _{\alpha}<\infty
\]
and%
\[
\left\Vert \varphi\right\Vert _{1+\alpha}=\left\Vert \varphi\right\Vert
_{0}+\left\Vert \varphi_{x}\right\Vert _{0}+\left\Vert \varphi_{xx}\right\Vert
_{0}+\left\lfloor \varphi_{xx}\right\rfloor _{\alpha}<\infty,
\]
respectively. Next, let $\mathbb{B}\left(  \left[  0,T\right]  ;C^{\alpha
}\left(
\mathbb{R}
^{2n}\right)  \right)  $, the set of all measurable functions
$f:\left[
0,T\right]  \times%
\mathbb{R}
^{2n}\rightarrow%
\mathbb{R}
$ such that for each $t\in\left[  0,T\right]  ,$ $f\left(  t,\cdot\right)  \in
C^{\alpha}\left(
\mathbb{R}
^{2n}\right)  $ and%
\[
\left\Vert f\right\Vert _{\mathbb{B}\left(  \left[  0,T\right]  ;C^{\alpha
}\left(
\mathbb{R}
^{2n}\right)  \right)  }=\sup_{t\in\left[  0,T\right]  }\left\Vert f\left(
t.\cdot\right)  \right\Vert _{\alpha}<\infty.
\]

Also, we let $\mathbb{C}\left(  \left[  0,T\right]  ;C^{\alpha}\left(
\mathbb{R}
^{2n}\right)  \right)  $ be the set of all continuous functions that are also
in $\mathbb{B}\left(  \left[  0,T\right]  ;C^{\alpha}\left(
\mathbb{R}
^{2n}\right)  \right)  $. \ Similarly, we define $\mathbb{B}\left(  \left[
0,T\right]  ;C^{k+\alpha}\left(
\mathbb{R}
^{2n}\right)  \right)  $ and $\mathbb{C}\left(  \left[  0,T\right]
;C^{k+\alpha}\left(
\mathbb{R}
^{2n}\right)  \right)  ,$ for $k=1,2.$

We now introduce the following assumptions.

\begin{enumerate}
\item[\textbf{(H)}] The maps $\bar{b},$ $\bar{\sigma},$ $\bar{f}$ and $\bar
{h}$ are continuous and bounded. Moreover, there exists a constant $K>0$ such
that $\forall\left(  t,x,p\right)  \in\left[  0,T\right]  \times%
\mathbb{R}
^{n}\times%
\mathbb{R}
^{n},$
\begin{align*}
K  &  \geq\left\vert \bar{\sigma}_{x}\left(  t,x\right)  \right\vert
+\left\vert \bar{b}_{x}\left(  t,x,p\right)  \right\vert +\left\vert \bar
{f}_{x}\left(  t,x,p\right)  \right\vert \\
&  +\left\vert \bar{b}_{p}\left(  t,x,p\right)  \right\vert +\left\vert
\bar{f}_{p}\left(  t,x,p\right)  \right\vert +\left\vert \bar{h}_{x}\left(
t,x\right)  \right\vert .
\end{align*}
Further $\mathbf{a}\left(  t,\mathbf{X}\right)  ^{-1}$ the inverse of
$\mathbf{a}\left(  t,\mathbf{X}\right)  $ exists for all $\left(
t,\mathbf{X}\right)  \in\left[  0,T\right]  \times%
\mathbb{R}
^{2n}$ and there exists constants $\lambda_{0},\lambda_{1}>0$ such that
$\forall\mathbf{X}\in%
\mathbb{R}
^{2n}$%
\[
\lambda_{1}\mathbf{I}_{2n}\leq\mathbf{a}\left(  t,\mathbf{X}\right)  ^{-1}%
\leq\lambda_{2}\mathbf{I}_{2n}\text{,}%
\]
where $\mathbf{I}_{2n}$ denotes the $\left(  2n\times2n\right)  $ identity matrix.
\end{enumerate}

The following theorem grantee the existence and uniqueness of a classical
solution to the PDEs (\ref{Eq212}); its proof follows an argument adapted from
Proof of Theorem 5.2. in Yong \cite{Yong2012}.

\begin{theorem}
Let Assumption \textbf{(H)} hold. Then (\ref{Eq212}) admits a unique solution
$\theta\left(  \cdot,\mathbf{\cdot}\right)  $.
\end{theorem}

\bop For any fixed $v\left(  \cdot,\cdot\right)  \in\mathbb{C}\left(  \left[
0,T\right]  ;C^{1+\alpha}\left(
\mathbb{R}
^{2n}\right)  \right)  $, we consider the following linear parabolic PDE%
\begin{equation}
\left\{
\begin{array}
[c]{l}%
0=\theta_{t}\left(  t,\mathbf{X}\right)  +\left[  \mathcal{L}^{0}\theta\left(
t,\mathbf{\cdot}\right)  \right]  \left(  \mathbf{X}\right)  +\left\langle
\theta_{\mathbf{X}}\left(  t,\mathbf{X}\right)  ,\mathbf{b}\left(
t,\mathbf{X,}v_{\mathbf{X}}\left(  t,\varphi\left(  \mathbf{X}\right)
\right)  \right)  \right\rangle \\
\text{ \ \ \ }+\mathbf{f}\left(  t,\mathbf{X},v_{\mathbf{X}}\left(
t,\varphi\left(  \mathbf{X}\right)  \right)  \right)  ,\text{ for }\left(
t,\mathbf{X}\right)  \in\left[  0,T\right]  \times%
\mathbb{R}
^{2n}\text{,}\\
\theta\left(  T,\mathbf{X}\right)  =\mathbf{F}\left(  \mathbf{X}\right)
,\text{ for }\mathbf{X}\in%
\mathbb{R}
^{2n}\text{.}%
\end{array}
\right.  \label{Eq206}%
\end{equation}
where $\mathcal{L}^{0}\left(  \cdot\right)  $ is the differentiable operator
defined as follows: For any \ $\psi\left(  \cdot\right)  \in C^{2}\left(
\mathbb{R}
^{2n}\right)  ,$%
\[
\left[  \mathcal{L}^{0}\psi\left(  \cdot\right)  \right]  \left(
\mathbf{X}\right)  =\text{\textbf{tr}}\left[  \mathbf{a}\left(  t,\mathbf{X}%
\right)  \psi_{\mathbf{XX}}\left(  \mathbf{X}\right)  \right]
\]

Applying Proposition 5.1 in Yong \cite{Yong2012}, we obtain that%
\begin{align*}
\theta\left(  t,\mathbf{X}\right)   &  =%
{\textstyle\int_{\mathbb{R}^{2n}}}
\Gamma^{0}\left(  t,\mathbf{X;}T,\mathbf{Z}\right)  \mathbf{\bar{F}}\left(
\mathbf{Z}\right)  \text{ }d\mathbf{Z}\\
&  \mathbf{+}%
{\textstyle\int_{t}^{T}}
{\textstyle\int_{\mathbb{R}^{2n}}}
\Gamma^{0}\left(  t,\mathbf{X;}s,\mathbf{Z}\right)  \left\langle
\mathbf{b}\left(  t,\mathbf{Z,}v_{\mathbf{X}}\left(  t,\varphi\left(
\mathbf{Z}\right)  \right)  \right)  ,\theta_{\mathbf{X}}\left(
t,\mathbf{Z}\right)  \right\rangle d\mathbf{Z}ds\\
&  +%
{\textstyle\int_{t}^{T}}
{\textstyle\int_{\mathbb{R}^{2n}}}
\Gamma^{0}\left(  t,\mathbf{X;}s,\mathbf{Z}\right)  \mathbf{f}\left(
t,\mathbf{Z},v_{\mathbf{X}}\left(  t,\varphi\left(  \mathbf{Z}\right)
\right)  \right)  d\mathbf{Z}ds\text{,}%
\end{align*}
where $\Gamma^{0}\left(  t,\mathbf{X;}s,\mathbf{Z}\right)  $ is the
fundamental solution of $\mathcal{L}^{0}\left(  \cdot\right)  $, given
explicitly by%
\[
\Gamma^{0}\left(  t,\mathbf{X;}s,\mathbf{Z}\right)  =\frac{1}{\left(
4\pi\left(  s-t\right)  \right)  ^{n}\left(  \det\left[  \mathbf{a}\left(
t,\mathbf{Z}\right)  \right]  \right)  ^{\frac{1}{2}}}e^{\frac{\left\langle
\mathbf{a}\left(  t,\mathbf{Z}\right)  ^{-1}\left(  \mathbf{X-Z}\right)
,\left(  \mathbf{X-Z}\right)  \right\rangle }{4\left(  s-t\right)  }}.
\]

On the other hand, by some computations one has (see~e.g. \cite{Friedman},
Page 24.)%
\begin{equation}
\left\{
\begin{array}
[c]{l}%
\left\vert \Gamma^{0}\left(  t,\mathbf{X;}s,\mathbf{Z}\right)  \right\vert
\leq\frac{K}{\left(  s-t\right)  ^{n}}e^{-\lambda\frac{\left\vert
\mathbf{X-Z}\right\vert ^{2}}{4\left(  s-t\right)  }},\\
\left\vert \Gamma_{\mathbf{X}}^{0}\left(  t,\mathbf{X;}s,\mathbf{Z}\right)
\right\vert \leq\frac{K}{\left(  s-t\right)  ^{\frac{2n+1}{2}}}e^{-\lambda
\frac{\left\vert \mathbf{X-Z}\right\vert ^{2}}{4\left(  s-t\right)  }}.
\end{array}
\right.  \text{ }\lambda<\lambda_{0}. \label{Eq207*}%
\end{equation}
Moreover, arguing as in \cite{Yong2012} we get,%
\begin{equation}
\Gamma_{\mathbf{Z}}^{0}\left(  t,\mathbf{X;}s,\mathbf{Z}\right)
=-\Gamma_{\mathbf{X}}^{0}\left(  t,\mathbf{X;}s,\mathbf{Z}\right)  -\Gamma
^{0}\left(  t,\mathbf{X;}s,\mathbf{Z}\right)  \rho\left(  t,\mathbf{X;}%
s,\mathbf{Z}\right)  , \label{Eq207}%
\end{equation}
where%
\[
\rho\left(  t,\mathbf{X;}s,\mathbf{Z}\right)  :=\frac{1}{2}\left(  \det\left[
\mathbf{a}\left(  t,\mathbf{Z}\right)  \right]  \right)  _{\mathbf{Z}}%
+\dfrac{\left\langle \left[  a\left(  t,\mathbf{Z}\right)  ^{-1}\right]
_{\mathbf{Z}}\left(  \mathbf{X-Z}\right)  ,\left(  \mathbf{X-Z}\right)
\right\rangle }{4\left(  s-t\right)  },
\]
with%
\[
\left\langle \left[  \mathbf{a}\left(  t,\mathbf{Z}\right)  ^{-1}\right]
_{\mathbf{Z}}\left(  \mathbf{X-Z}\right)  ,\left(  \mathbf{X-Z}\right)
\right\rangle :=\left(
\begin{array}
[c]{c}%
\left\langle \left[  a\left(  t,\mathbf{Z}\right)  ^{-1}\right]  _{z_{1}%
}\left(  \mathbf{X-Z}\right)  ,\left(  \mathbf{X-Z}\right)  \right\rangle \\
\left\langle \left[  a\left(  t,\mathbf{Z}\right)  ^{-1}\right]  _{z_{2}%
}\left(  \mathbf{X-Z}\right)  ,\left(  \mathbf{X-Z}\right)  \right\rangle \\
...\\
\left\langle \left[  a\left(  t,\mathbf{Z}\right)  ^{-1}\right]  _{z_{2n}%
}\left(  \mathbf{X-Z}\right)  ,\left(  \mathbf{X-Z}\right)  \right\rangle
\end{array}
\right)  .
\]
Under Assumption \textbf{(H), }we have%
\begin{equation}
\rho\left(  t,\mathbf{X;}s,\mathbf{Z}\right)  \leq K\left(  1+\dfrac
{\left\vert \mathbf{X-Z}\right\vert ^{2}}{s-t}\right)  . \label{Eq208}%
\end{equation}
Then using (\ref{Eq207}), we have%
\begin{align}
&  \theta_{\mathbf{X}}\left(  t,\mathbf{X}\right) \nonumber\\
&  =%
{\textstyle\int_{\mathbb{R}^{2n}}}
\Gamma_{\mathbf{X}}^{0}\left(  t,\mathbf{X;}T,\mathbf{Z}\right)
\mathbf{F}\left(  \mathbf{Z}\right)  d\mathbf{Z}\nonumber\\
&  \mathbf{+}%
{\textstyle\int_{t}^{T}}
{\textstyle\int_{\mathbb{R}^{2n}}}
\Gamma_{\mathbf{X}}^{0}\left(  t,\mathbf{X;}s,\mathbf{Z}\right)  \left\langle
\mathbf{b}\left(  s,\mathbf{Z,}v_{\mathbf{X}}\left(  s,\varphi\left(
\mathbf{Z}\right)  \right)  \right)  ,\theta_{\mathbf{X}}\left(
s,\mathbf{Z}\right)  \right\rangle d\mathbf{Z}ds\nonumber\\
&  +%
{\textstyle\int_{t}^{T}}
{\textstyle\int_{\mathbb{R}^{2n}}}
\Gamma_{\mathbf{X}}^{0}\left(  t,\mathbf{X;}s,\mathbf{Z}\right)
\mathbf{f}\left(  s,\mathbf{Z,}v_{\mathbf{X}}\left(  s,\varphi\left(
\mathbf{Z}\right)  \right)  \right)  d\mathbf{Z}ds\nonumber\\
&  =-%
{\textstyle\int_{\mathbb{R}^{2n}}}
\Gamma_{\mathbf{Z}}^{0}\left(  t,\mathbf{X;}s,\mathbf{Z}\right)
\mathbf{F}\left(  \mathbf{Z}\right)  d\mathbf{Z}-%
{\textstyle\int_{\mathbb{R}^{2n}}}
\Gamma^{0}\left(  t,\mathbf{X;}s,\mathbf{Z}\right)  \rho\left(  t,\mathbf{X;}%
s,\mathbf{Z}\right)  \mathbf{F}\left(  \mathbf{Z}\right)  d\mathbf{Z}%
\nonumber\\
&  \mathbf{+}%
{\textstyle\int_{t}^{T}}
{\textstyle\int_{\mathbb{R}^{2n}}}
\Gamma_{\mathbf{X}}^{0}\left(  t,\mathbf{X;}s,\mathbf{Z}\right)  \left\langle
\mathbf{b}\left(  s,\mathbf{Z,}v_{\mathbf{X}}\left(  s,\varphi\left(
\mathbf{Z}\right)  \right)  \right)  ,\theta_{\mathbf{X}}\left(
s,\mathbf{Z}\right)  \right\rangle d\mathbf{Z}ds\nonumber\\
&  +%
{\textstyle\int_{t}^{T}}
{\textstyle\int_{\mathbb{R}^{2n}}}
\Gamma_{\mathbf{X}}^{0}\left(  t,\mathbf{X;}s,\mathbf{Z}\right)
\mathbf{f}\left(  s,\mathbf{Z,}v_{\mathbf{X}}\left(  t,\varphi\left(
\mathbf{Z}\right)  \right)  \right)  d\mathbf{Z}ds\nonumber\\
&  =%
{\textstyle\int_{\mathbb{R}^{2n}}}
\Gamma^{0}\left(  t,\mathbf{X;}s,\mathbf{Z}\right)  \mathbf{F}_{\mathbf{Z}%
}\left(  \mathbf{Z}\right)  d\mathbf{Z}-%
{\textstyle\int_{\mathbb{R}^{2n}}}
\Gamma^{0}\left(  t,\mathbf{X;}s,\mathbf{Z}\right)  \rho\left(  t,\mathbf{X;}%
s,\mathbf{Z}\right)  \mathbf{F}\left(  \mathbf{Z}\right)  d\mathbf{Z}%
\nonumber\\
&  \mathbf{+}%
{\textstyle\int_{t}^{T}}
{\textstyle\int_{\mathbb{R}^{2n}}}
\Gamma_{\mathbf{X}}^{0}\left(  t,\mathbf{X;}s,\mathbf{Z}\right)  \left\langle
\mathbf{b}\left(  s,\mathbf{Z,}v_{\mathbf{X}}\left(  s,\varphi\left(
\mathbf{Z}\right)  \right)  \right)  ,\theta_{\mathbf{X}}\left(
s,\mathbf{Z}\right)  \right\rangle d\mathbf{Z}ds\nonumber\\
&  +%
{\textstyle\int_{t}^{T}}
{\textstyle\int_{\mathbb{R}^{2n}}}
\Gamma_{\mathbf{X}}^{0}\left(  t,\mathbf{X;}s,\mathbf{Z}\right)
\mathbf{f}\left(  s,\mathbf{Z,}v_{\mathbf{X}}\left(  s,\varphi\left(
\mathbf{Z}\right)  \right)  \right)  d\mathbf{Z}ds. \label{Eq209}%
\end{align}
Therefor, combining (\ref{Eq207*})-(\ref{Eq208}) together with (\ref{Eq209})
and using the fact that $\mathbf{f}$ and $\mathbf{b}$ are uniformly bounded,
we obtain that%
\begin{align}
\left\vert \theta_{\mathbf{X}}\left(  t,\mathbf{X}\right)  \right\vert  &
\leq%
{\textstyle\int_{\mathbb{R}^{2n}}}
\frac{K}{\left(  T-t\right)  ^{n}}e^{-\lambda\frac{\left\vert \mathbf{X-Z}%
\right\vert ^{2}}{4\left(  T-t\right)  }}\left[  \left\vert \mathbf{F}%
_{\mathbf{Z}}\left(  \mathbf{Z}\right)  \right\vert +\left(  1+\frac
{\left\vert \mathbf{X-Z}\right\vert ^{2}}{T-t}\right)  \left\vert
\mathbf{F}\left(  \mathbf{Z}\right)  \right\vert \right]  d\mathbf{Z}%
\nonumber\\
&  +%
{\textstyle\int_{t}^{T}}
{\textstyle\int_{\mathbb{R}^{2n}}}
\frac{K}{\left(  s-t\right)  ^{\frac{2n+1}{2}}}e^{-\lambda\frac{\left\vert
\mathbf{X-Z}\right\vert ^{2}}{4\left(  s-t\right)  }}\left[  \left\vert
\theta_{\mathbf{X}}\left(  t,\mathbf{Z}\right)  \right\vert +1\right]
d\mathbf{Z}ds\nonumber\\
&  \leq K\left(  \left\Vert F\left(  \cdot\right)  \right\Vert _{C^{1}\left(
\mathbb{R}
^{n}\right)  }+1\right)  +%
{\textstyle\int_{t}^{T}}
{\textstyle\int_{\mathbb{R}^{2n}}}
\frac{K}{\left(  s-t\right)  ^{\frac{2n+1}{2}}}e^{-\lambda\frac{\left\vert
\mathbf{X-Z}\right\vert ^{2}}{4\left(  s-t\right)  }}\left\vert \theta
_{\mathbf{X}}\left(  s,\mathbf{Z}\right)  \right\vert d\mathbf{Z}ds
\label{Eq211*}%
\end{align}
Thus we also have%
\[
\left\vert \theta_{\mathbf{X}}\left(  s,\mathbf{Z}\right)  \right\vert \leq
K\left(  \left\Vert F\left(  \cdot\right)  \right\Vert _{C^{1}\left(
\mathbb{R}
^{n}\right)  }+1\right)  +%
{\textstyle\int_{s}^{T}}
{\textstyle\int_{\mathbb{R}^{2n}}}
\frac{K}{\left(  r-s\right)  ^{\frac{2n+1}{2}}}e^{-\lambda\frac{\left\vert
\mathbf{Z-Y}\right\vert ^{2}}{4\left(  r-s\right)  }}\left\vert \theta
_{\mathbf{X}}\left(  r,\mathbf{Z}\right)  \right\vert d\mathbf{Z}dr.
\]
Accordingly, by using Lemma 3 in (\cite{Friedman}, Page 24), we get%
\begin{align*}
&
{\textstyle\int_{t}^{T}}
{\textstyle\int_{\mathbb{R}^{2n}}}
\frac{K}{\left(  s-t\right)  ^{\frac{2n+1}{2}}}e^{-\lambda\frac{\left\vert
\mathbf{X-Z}\right\vert ^{2}}{4\left(  s-t\right)  }}\left\vert \theta
_{\mathbf{X}}\left(  s,\mathbf{Z}\right)  \right\vert d\mathbf{Z}ds\\
&  \leq%
{\textstyle\int_{t}^{T}}
{\textstyle\int_{\mathbb{R}^{2n}}}
\frac{K}{\left(  s-t\right)  ^{\frac{2n+1}{2}}}e^{-\lambda\frac{\left\vert
\mathbf{X-Z}\right\vert ^{2}}{4\left(  s-t\right)  }}\left\vert \left(
\left\Vert F\left(  \cdot\right)  \right\Vert _{C^{1}\left(
\mathbb{R}
^{n}\right)  }+1\right)  \right\vert d\mathbf{Z}ds\\
&  +%
{\textstyle\int_{t}^{T}}
{\textstyle\int_{\mathbb{R}^{2n}}}
\frac{K}{\left(  s-t\right)  ^{\frac{2n+1}{2}}}e^{-\lambda\frac{\left\vert
\mathbf{X-Z}\right\vert ^{2}}{4\left(  s-t\right)  }}%
{\textstyle\int_{s}^{T}}
{\textstyle\int_{\mathbb{R}^{2n}}}
\frac{K}{\left(  r-s\right)  ^{\frac{2n+1}{2}}}e^{-\lambda\frac{\left\vert
\mathbf{Z-Y}\right\vert ^{2}}{4\left(  r-s\right)  }}\left\vert \theta
_{\mathbf{X}}\left(  r,\mathbf{Y}\right)  \right\vert d\mathbf{Y}%
drd\mathbf{Z}ds\\
&  \leq K\left(  \left\Vert F\left(  \cdot\right)  \right\Vert _{C^{1}\left(
\mathbb{R}
^{n}\right)  }+1\right) \\
&  +%
{\textstyle\int_{t}^{T}}
{\textstyle\int_{\mathbb{R}^{2n}}}
\left(
{\textstyle\int_{t}^{r}}
{\textstyle\int_{\mathbb{R}^{2n}}}
\frac{K}{\left(  s-t\right)  ^{\frac{2n+1}{2}}}e^{-\lambda\frac{\left\vert
\mathbf{X-Z}\right\vert ^{2}}{4\left(  s-t\right)  }}\frac{K}{\left(
r-s\right)  ^{\frac{2n+1}{2}}}e^{-\lambda\frac{\left\vert \mathbf{Z-Y}%
\right\vert ^{2}}{4\left(  r-s\right)  }}d\mathbf{Z}ds\right)  \left\vert
\theta_{\mathbf{X}}\left(  r,\mathbf{Y}\right)  \right\vert d\mathbf{Y}dr\\
&  \leq K\left(  \left\Vert F\left(  \cdot\right)  \right\Vert _{C^{1}\left(
\mathbb{R}
^{n}\right)  }+1\right)  +%
{\textstyle\int_{t}^{T}}
{\textstyle\int_{\mathbb{R}^{2n}}}
\frac{K}{\left(  r-t\right)  ^{\frac{2n}{2}}}e^{-\lambda\frac{\left\vert
\mathbf{X-Y}\right\vert ^{2}}{4\left(  r-t\right)  }}\left\vert \theta
_{\mathbf{X}}\left(  r,\mathbf{Y}\right)  \right\vert d\mathbf{Y}dr.
\end{align*}
Invoking this into (\ref{Eq211*})$,$ we get%
\begin{equation}
\left\vert \theta_{\mathbf{X}}\left(  t,\mathbf{X}\right)  \right\vert \leq
K\left(  \left\Vert F\left(  \cdot\right)  \right\Vert _{C^{1}\left(
\mathbb{R}
^{n}\right)  }+1\right)  +%
{\textstyle\int_{t}^{T}}
{\textstyle\int_{\mathbb{R}^{2n}}}
\frac{K}{\left(  r-t\right)  ^{\frac{2n}{2}}}e^{-\lambda\frac{\left\vert
\mathbf{X-Y}\right\vert ^{2}}{4\left(  r-t\right)  }}\left\vert \theta
_{\mathbf{X}}\left(  r,\mathbf{Y}\right)  \right\vert d\mathbf{Y}dr
\label{Eq213}%
\end{equation}

We can repeat the above procedure $2n$ times and then use Gronwall's
inequality to obtain%
\begin{equation}
\sup_{\mathbf{X\in%
\mathbb{R}
}^{2n}}\left\vert \theta_{\mathbf{X}}\left(  t,\mathbf{X}\right)  \right\vert
\leq K\left(  \left\Vert F\left(  \cdot\right)  \right\Vert _{C^{1}\left(
\mathbb{R}
^{n}\right)  }+1\right)  . \label{Eq214}%
\end{equation}

Now let $v_{1}\left(  \cdot,\cdot\right)  $, $v_{2}\left(  \cdot,\cdot\right)
\in\mathbb{C}\left(  \left[  0,T\right]  ;C^{1}\left(
\mathbb{R}
^{2n}\right)  \right)  $ and we consider $\theta^{1}\left(  \cdot
,\cdot\right)  $, $\theta^{2}\left(  \cdot,\cdot\right)  $ their corresponding
solutions of (\ref{Eq206}). Then%
\begin{align}
&  \theta^{1}\left(  t,\mathbf{X}\right)  -\theta^{2}\left(  t,\mathbf{X}%
\right) \nonumber\\
&  =%
{\textstyle\int_{t}^{T}}
{\textstyle\int_{\mathbb{R}^{2n}}}
\Gamma^{0}\left(  t,\mathbf{X;}s,\mathbf{Z}\right)  \left(  \left\langle
\mathbf{b}\left(  s,\mathbf{Z,}v_{\mathbf{X}}^{1}\left(  s,\varphi\left(
\mathbf{Z}\right)  \right)  \right)  ,\theta_{\mathbf{X}}^{1}\left(
s,\mathbf{Z}\right)  -\theta_{\mathbf{X}}^{2}\left(  s,\mathbf{Z}\right)
\right\rangle \right. \nonumber\\
&  +\left\langle \mathbf{b}\left(  s,\mathbf{Z,}v_{\mathbf{X}}^{1}\left(
s,\varphi\left(  \mathbf{Z}\right)  \right)  \right)  -\mathbf{b}\left(
s,\mathbf{Z,}v_{\mathbf{X}}^{2}\left(  s,\varphi\left(  \mathbf{Z}\right)
\right)  \right)  ,\theta_{\mathbf{X}}^{2}\left(  s,\mathbf{Z}\right)
\right\rangle \nonumber\\
&  \left.  +\mathbf{f}\left(  s,\mathbf{Z,}v_{\mathbf{X}}^{1}\left(
s,\varphi\left(  \mathbf{Z}\right)  \right)  \right)  -\mathbf{f}\left(
s,\mathbf{Z,}v_{\mathbf{X}}^{2}\left(  s,\varphi\left(  \mathbf{Z}\right)
\right)  \right)  \right)  d\mathbf{Z}ds \label{Eq215}%
\end{align}
and%
\begin{align*}
&  \theta_{\mathbf{X}}^{1}\left(  t,\mathbf{X}\right)  -\theta_{\mathbf{X}%
}^{2}\left(  t,\mathbf{X}\right) \\
&  =%
{\textstyle\int_{t}^{T}}
{\textstyle\int_{\mathbb{R}^{2n}}}
\Gamma_{\mathbf{X}}^{0}\left(  t,\mathbf{X;}s,\mathbf{Z}\right)  \left(
\left\langle \mathbf{b}\left(  s,\mathbf{Z,}v_{\mathbf{X}}^{1}\left(
s,\varphi\left(  \mathbf{Z}\right)  \right)  \right)  ,\theta_{\mathbf{X}}%
^{1}\left(  s,\mathbf{Z}\right)  -\theta_{\mathbf{X}}^{2}\left(
s,\mathbf{Z}\right)  \right\rangle \right. \\
&  +\left\langle \mathbf{b}\left(  s,\mathbf{Z,}v_{\mathbf{X}}^{1}\left(
s,\varphi\left(  \mathbf{Z}\right)  \right)  \right)  -\mathbf{b}\left(
s,\mathbf{Z,}v_{\mathbf{X}}^{2}\left(  s,\varphi\left(  \mathbf{Z}\right)
\right)  \right)  ,\theta_{\mathbf{X}}^{2}\left(  s,\mathbf{Z}\right)
\right\rangle \\
&  \left.  +\mathbf{f}\left(  s,\mathbf{Z,}v_{\mathbf{X}}^{1}\left(
s,\varphi\left(  \mathbf{Z}\right)  \right)  \right)  -\mathbf{f}\left(
s,\mathbf{Z,}v_{\mathbf{X}}^{2}\left(  s,\varphi\left(  \mathbf{Z}\right)
\right)  \right)  \right)  d\mathbf{Z}ds.
\end{align*}
Accordingly, by Assumption (\textbf{H) }we have,%
\begin{align*}
&  \left\vert \theta_{\mathbf{X}}^{1}\left(  t,\mathbf{X}\right)
-\theta_{\mathbf{X}}^{2}\left(  t,\mathbf{X}\right)  \right\vert \\
&  \leq%
{\textstyle\int_{t}^{T}}
{\textstyle\int_{\mathbb{R}^{2n}}}
\frac{K}{\left(  s-t\right)  ^{\frac{2n+1}{2}}}e^{-\lambda\frac{\left\vert
\mathbf{X-Z}\right\vert ^{2}}{4\left(  s-t\right)  }}\left(  \left\vert
\theta_{\mathbf{X}}^{1}\left(  s,\mathbf{Z}\right)  -\theta_{\mathbf{X}}%
^{2}\left(  s,\mathbf{Z}\right)  \right\vert \right. \\
&  \left.  +\left(  1+\left\vert \theta_{\mathbf{X}}^{2}\left(  s,\mathbf{Z}%
\right)  \right\vert \right)  \left\vert v_{\mathbf{X}}^{1}\left(
s,\varphi\left(  \mathbf{Z}\right)  \right)  -v_{\mathbf{X}}^{2}\left(
s,\varphi\left(  \mathbf{Z}\right)  \right)  \right\vert \right)
d\mathbf{Z}ds\\
&  \leq%
{\textstyle\int_{t}^{T}}
{\textstyle\int_{\mathbb{R}^{2n}}}
\frac{K}{\left(  s-t\right)  ^{\frac{2n+1}{2}}}e^{-\lambda\frac{\left\vert
\mathbf{X-Z}\right\vert ^{2}}{4\left(  s-t\right)  }}\left\vert \theta
_{\mathbf{X}}^{1}\left(  s,\mathbf{Z}\right)  -\theta_{\mathbf{X}}^{2}\left(
s,\mathbf{Z}\right)  \right\vert d\mathbf{Z}ds\\
&  +K\left(  T-t\right)  ^{\frac{1}{2}}\left(  \left\Vert F\left(
\cdot\right)  \right\Vert _{C^{1}\left(
\mathbb{R}
^{n}\right)  }+1\right)  \left\Vert v_{\mathbf{X}}^{1}\left(  \cdot
,\mathbf{\cdot}\right)  -v_{\mathbf{X}}^{2}\left(  \cdot,\mathbf{\cdot
}\right)  \right\Vert _{\mathbb{C}\left(  \left[  t,T\right]  ;C\left(
\mathbb{R}
^{2n}\right)  \right)  }%
\end{align*}
Then, following an iterative procedure as the one used to obtain
(\ref{Eq214}), we can obtain%
\begin{align*}
\sup_{\mathbf{X\in%
\mathbb{R}
}^{2n}}  &  \left\vert \theta_{\mathbf{X}}^{1}\left(  t,\mathbf{X}\right)
-\theta_{\mathbf{X}}^{2}\left(  t,\mathbf{X}\right)  \right\vert \\
&  \leq K\left(  T-t\right)  ^{\frac{1}{2}}\left(  \left\Vert F\left(
\cdot\right)  \right\Vert _{C^{1}\left(
\mathbb{R}
^{n}\right)  }+1\right)  \left\Vert v_{\mathbf{X}}^{1}\left(  \cdot
,\mathbf{\cdot}\right)  -v_{\mathbf{X}}^{2}\left(  \cdot,\mathbf{\cdot
}\right)  \right\Vert _{\mathbb{C}\left(  \left[  t,T\right]  ;C\left(
\mathbb{R}
^{2n}\right)  \right)  }.
\end{align*}
On the other hand, from (\ref{Eq215}), we have%
\begin{align*}
&  \sup_{\mathbf{X\in%
\mathbb{R}
}^{2n}}\left\vert \theta^{1}\left(  t,\mathbf{X}\right)  -\theta^{2}\left(
t,\mathbf{X}\right)  \right\vert \\
&  \leq%
{\textstyle\int_{t}^{T}}
{\textstyle\int_{\mathbb{R}^{2n}}}
\left(  \frac{K}{\left(  s-t\right)  ^{\frac{2n}{2}}}e^{-\lambda
\frac{\left\vert \mathbf{X-Z}\right\vert ^{2}}{4\left(  s-t\right)  }%
}\left\vert \theta_{\mathbf{X}}^{1}\left(  s,\mathbf{Z}\right)  -\theta
_{\mathbf{X}}^{2}\left(  s,\mathbf{Z}\right)  \right\vert \right. \\
&  \left.  +\left(  1+\left\vert \theta_{\mathbf{X}}^{2}\left(  s,\mathbf{Z}%
\right)  \right\vert \right)  \left\vert v_{\mathbf{X}}^{1}\left(
s,\varphi\left(  \mathbf{Z}\right)  \right)  -v_{\mathbf{X}}^{2}\left(
s,\varphi\left(  \mathbf{Z}\right)  \right)  \right\vert \right)
d\mathbf{Z}ds\\
&  \leq K\left(  T-t\right)  \left(  \left\Vert F\left(  \cdot\right)
\right\Vert _{C^{1}\left(
\mathbb{R}
^{n}\right)  }+1\right)  \left\Vert v_{\mathbf{X}}^{1}\left(  \cdot
,\mathbf{\cdot}\right)  -v_{\mathbf{X}}^{2}\left(  \cdot,\mathbf{\cdot
}\right)  \right\Vert _{\mathbb{C}\left(  \left[  t,T\right]  ;C\left(
\mathbb{R}
^{2n}\right)  \right)  }.
\end{align*}
Hence we obtain%
\begin{align*}
\sup_{\mathbf{X\in%
\mathbb{R}
}^{2n}}  &  \left\vert \theta^{1}\left(  t,\mathbf{X}\right)  -\theta
^{2}\left(  t,\mathbf{X}\right)  \right\vert +\sup_{\mathbf{X\in%
\mathbb{R}
}^{2n}}\left\vert \theta_{X}^{1}\left(  t,\mathbf{X}\right)  -\theta_{X}%
^{2}\left(  t,\mathbf{X}\right)  \right\vert \\
&  \leq K\left(  T-t\right)  ^{\frac{1}{2}}\left(  \left\Vert F\left(
\cdot\right)  \right\Vert _{C^{1}\left(
\mathbb{R}
^{n}\right)  }+1\right)  \left\Vert v_{\mathbf{X}}^{1}\left(  \cdot
,\mathbf{\cdot}\right)  -v_{\mathbf{X}}^{2}\left(  \cdot,\mathbf{\cdot
}\right)  \right\Vert _{\mathbb{C}\left(  \left[  t,T\right]  ;C\left(
\mathbb{R}
^{2n}\right)  \right)  }\\
&  \leq K\left(  T-t\right)  ^{\frac{1}{2}}\left(  \left\Vert F\left(
\cdot\right)  \right\Vert _{C^{1}\left(
\mathbb{R}
^{n}\right)  }+1\right)  \left\Vert v^{1}\left(  \cdot,\mathbf{\cdot}\right)
-v^{2}\left(  \cdot,\mathbf{\cdot}\right)  \right\Vert _{\mathbb{C}\left(
\left[  t,T\right]  ;C^{1}\left(
\mathbb{R}
^{2n}\right)  \right)  }.
\end{align*}
Thus, in particular, we have%
\begin{align*}
&  \left\Vert \theta^{1}\left(  \cdot,\mathbf{\cdot}\right)  -\theta
^{2}\left(  \cdot,\mathbf{\cdot}\right)  \right\Vert _{\mathbb{C}\left(
\left[  T-\delta,T\right]  ;C^{1}\left(
\mathbb{R}
^{2n}\right)  \right)  }\\
&  \leq K\delta^{\frac{1}{2}}\left(  \left\Vert F\left(  \cdot\right)
\right\Vert _{C^{1}\left(
\mathbb{R}
^{n}\right)  }+1\right)  \left\Vert v^{1}\left(  \cdot,\cdot\right)
-v^{1}\left(  \cdot,\cdot\right)  \right\Vert _{\mathbb{C}\left(  \left[
T-\delta,T\right]  ;C^{1}\left(
\mathbb{R}
^{2n}\right)  \right)  }.
\end{align*}
Clearly, by choosing $\delta>0$ small enough, we get a contraction mapping on
$v\left(  \cdot,\cdot\right)  \rightarrow\theta\left(  \cdot,\mathbf{\cdot
}\right)  $ on \newline$\mathbb{C}\left(  \left[  T-\delta,T\right]
;C^{1}\left(
\mathbb{R}
^{2n}\right)  \right)  .$ Therefore, this map admits a unique fixed point.
Since we may obtain similar estimates on $\mathbb{C}\left(  \left[
T-2\delta,T-\delta\right]  ;C^{1}\left(
\mathbb{R}
^{2n}\right)  \right)  $, etc., one sees that the fixed point will exists on
the whole space \newline$\mathbb{C}\left(  \left[  0,T\right]  ;C^{1}\left(
\mathbb{R}
^{2n}\right)  \right)  $ for the map $v\left(  \cdot,\cdot\right)
\rightarrow\theta\left(  \cdot,\mathbf{\cdot}\right)  $. Then we obtain the
well-posedness of the following%
\begin{align*}
&  \theta\left(  t,\mathbf{X}\right) \\
&  =%
{\textstyle\int_{\mathbb{R}^{2n}}}
\Gamma^{0}\left(  t,\mathbf{X;}T,\mathbf{Z}\right)  \mathbf{F}\left(
\mathbf{Z}\right)  \text{ }d\mathbf{Z}\\
&  \mathbf{+}%
{\textstyle\int_{t}^{T}}
\Gamma^{0}\left(  t,\mathbf{X;}s,\mathbf{Z}\right)  \left\langle
\mathbf{b}\left(  t,\mathbf{Z,}\theta_{\mathbf{X}}\left(  t,\varphi\left(
\mathbf{Z}\right)  \right)  \right)  ,\theta_{\mathbf{X}}\left(
t,\mathbf{Z}\right)  \right\rangle d\mathbf{Z}ds\\
&  +%
{\textstyle\int_{t}^{T}}
\Gamma^{0}\left(  t,\mathbf{X;}s,\mathbf{Z}\right)  \mathbf{f}\left(
t,\mathbf{Z},\theta_{\mathbf{X}}\left(  t,\varphi\left(  \mathbf{Z}\right)
\right)  \right)  d\mathbf{Z}ds.
\end{align*}
Finally, by the regularity of $\Gamma^{0}\left(  t,\mathbf{X;}T,\mathbf{Z}%
\right)  $, we know that $\theta\left(  t,\mathbf{X}\right)  $ is
$C^{2+\alpha}$ in $\mathbf{X=}\left(  x,y\right)  \mathbf{,}$ $C^{1+\alpha}$
in $t$ for some $\alpha\in\left(  0,1\right)  $, and the PDE (\ref{Eq212}) is
satisfied.\eop

\bigskip

Department of Mathematics University of Bordj Bou Arreridj 34000 Algeria;
izacalia@yahoo.com, fermasof@yahoo.fr.

\end{document}